\documentclass[12pt]{amsart}

\usepackage[paperheight=11in,paperwidth=8.5in,left=1in,top=1.25in,right=1in,bottom=1.25in,]{geometry}
\usepackage[linktocpage=false,colorlinks=true,linkcolor=black,citecolor=black,
filecolor=black,urlcolor=black,pagebackref=false,pdfstartpage={1},
pdfstartview={FitH},pdftitle={},pdfauthor={},pdfsubject={},pdfcreator={},
pdfproducer={},pdfkeywords={}]{hyperref}

\usepackage[lofdepth,lotdepth]{subfig}

\usepackage{comment}
\usepackage{afterpage}
\usepackage{amsfonts}
\usepackage{amsmath}
\allowdisplaybreaks[4]
\usepackage{amssymb}
\usepackage{amsthm}
\usepackage{bbm}
\usepackage[justification=centering]{caption}
\usepackage[capitalize,nameinlink,noabbrev,nosort]{cleveref}
\usepackage{enumerate}
\usepackage{enumitem}
\usepackage{float}
\restylefloat{figure}
\usepackage{graphicx}
\usepackage{mathrsfs}
\usepackage{mathtools}
\usepackage[framemethod=tikz,ntheorem,xcolor]{mdframed}
\usepackage{parcolumns}
\usepackage{tabularx}
\usepackage{tikz}\pdfpageattr{/Group <</S /Transparency /I true /CS /DeviceRGB>>}
\usetikzlibrary{arrows,calc,intersections,matrix,positioning,through}
\usepackage{tikz-cd}
\tikzset{commutative diagrams/.cd,every label/.append style = {font = \normalsize}}
\usepackage[all]{xy}
\usepackage[aligntableaux=center]{ytableau}
\usepackage{epstopdf}

\numberwithin{equation}{section}

\newtheorem{thm}[equation]{Theorem}
\AfterEndEnvironment{thm}{\noindent\ignorespaces}
\newtheorem{cor}[equation]{Corollary}
\AfterEndEnvironment{cor}{\noindent\ignorespaces}
\newtheorem{lem}[equation]{Lemma}
\AfterEndEnvironment{lem}{\noindent\ignorespaces}
\newtheorem{prop}[equation]{Proposition}
\AfterEndEnvironment{prop}{\noindent\ignorespaces}

\AfterEndEnvironment{conj}{\noindent\ignorespaces}
\newtheorem{prob}[equation]{Problem}
\AfterEndEnvironment{prob}{\noindent\ignorespaces}

\theoremstyle{definition}
\newtheorem{defn}[equation]{Definition}

\AfterEndEnvironment{defn}{\noindent\ignorespaces}
\AfterEndEnvironment{definition}{\noindent\ignorespaces}
\newtheorem*{pf_no_qed}{Proof}
\newenvironment{pf}[1][]{\begin{pf_no_qed}[#1]\pushQED{\qed}}{\popQED\end{pf_no_qed}}
\AfterEndEnvironment{pf}{\noindent\ignorespaces}
\newtheorem{eg_no_qed}[equation]{Example}
\newenvironment{eg}[1][]{\begin{eg_no_qed}[#1]\pushQED{\qed}}{\popQED\end{eg_no_qed}}
\AfterEndEnvironment{eg}{\noindent\ignorespaces}

\AfterEndEnvironment{example}{\noindent\ignorespaces}
\newtheorem{rmk}[equation]{Remark}
\AfterEndEnvironment{rmk}{\noindent\ignorespaces}

\AfterEndEnvironment{remark}{\noindent\ignorespaces}

\theoremstyle{remark}
\newtheorem*{claim}{Claim}
\AfterEndEnvironment{claim}{\noindent\ignorespaces}
\newtheorem*{claimpf_no_qed}{Proof of Claim}
\newenvironment{claimpf}[1][]{\begin{claimpf_no_qed}[#1]\pushQED{\qed}}{\popQED\end{claimpf_no_qed}}
\AfterEndEnvironment{claimpf}{\noindent\ignorespaces}

\newcommand{\Le}{\textup{\protect\scalebox{-1}[1]{L}}}

\font\pipefont=lcircle10
\def\elbow{\smash{\raise3pt\hbox{\pipefont\rlap{\rlap{\char'014}\char'016}}}}
\def\textelbow{\;\;\,\elbow\;\;\,}
\def\halfelbow{\smash{\raise2pt\hbox{\pipefont\rlap{\rlap{\rlap{\char'015}\phantom{\char'017}}}}}}
\def\cross{\smash{\lower5pt\hbox{\rlap{\vrule height16pt}}\raise3pt\hbox{\rlap{\hskip-8pt \vrule height0.4pt depth0pt width16pt}}}}
\def\textcross{\;\;\,\cross\;\;\,}
\ytableausetup{boxsize=16pt}

\def\DS{\Omega_{DS}}
\def\DM{\Omega_{DM}}
\def\O{\mathcal{O}}

\DeclareMathOperator{\alt}{alt}

\newcommand{\R}{\mathbb{R}}
\DeclareMathOperator{\colspan}{colspan}

\DeclareMathOperator{\Fl}{Fl}

\DeclareMathOperator{\Gr}{Gr}

\DeclareMathOperator{\interior}{int}

\newcommand{\PSign}{\mathbb{P}\!\operatorname{Sign}}
\DeclareMathOperator{\rank}{rank}

\DeclareMathOperator{\rowspan}{rowspan}
\newcommand{\rf}[1]{\hyperref[#1]{(\ref*{#1})}}

\DeclareMathOperator{\sign}{sign}
\DeclareMathOperator{\Sign}{Sign}

\DeclareMathOperator{\Slide}{Slide}

\DeclareMathOperator{\spn}{span}

\DeclareMathOperator{\var}{var}


\title{The $m=1$ amplituhedron and cyclic hyperplane arrangements}

\author{Steven N.\ Karp}
\address{Department of Mathematics, University of California, Berkeley}
\email{\href{mailto:skarp@berkeley.edu}{skarp@berkeley.edu}}
\email{\href{mailto:williams@math.berkeley.edu}{williams@math.berkeley.edu}}

\author{Lauren K.\ Williams}

\thanks{LW was partially supported by a Rose Hills Innovator award
and the NSF CAREER award DMS-1049513.  
Both authors were partially supported by the NSF grant DMS-1600447.}

\begin{document}

\begin{abstract}
The \emph{(tree) amplituhedron} $\mathcal{A}_{n,k,m}$ is the image 
in the Grassmannian $\Gr_{k,k+m}$ of the totally nonnegative part of 
$\Gr_{k,n}$, under a (map induced by a) linear map which is totally 
positive. It was introduced by 
Arkani-Hamed and Trnka in 2013 in order to give a geometric basis for the 
computation of scattering amplitudes in $\mathcal{N}=4$
supersymmetric Yang-Mills theory.  When $k+m=n$, the amplituhedron
is isomorphic to the totally nonnegative Grassmannian, and when
$k=1$, the amplituhedron is a cyclic polytope.
While the case $m=4$ is most relevant to physics, the amplituhedron is an interesting mathematical object for any $m$.  In this paper we study it in the case $m=1$. 
We start by taking an orthogonal point of 
view and define a related ``B-amplituhedron'' $\mathcal{B}_{n,k,m}$,
which we show is isomorphic to $\mathcal{A}_{n,k,m}$.  We use this 
reformulation to describe the amplituhedron in terms of sign variation.
We then give a cell decomposition of the amplituhedron 
$\mathcal{A}_{n,k,1}$ 
using the images of a collection of distinguished cells of the totally nonnegative Grassmannian.  
We also show that $\mathcal{A}_{n,k,1}$ can 
be identified with the complex of bounded faces of 
a \emph{cyclic hyperplane arrangement}, and describe how its cells fit together.  We deduce that $\mathcal{A}_{n,k,1}$ is homeomorphic to a ball.
\end{abstract}

\maketitle
\setcounter{tocdepth}{1}
\tableofcontents

\section{Introduction}\label{sec_intro}

\noindent The totally nonnegative Grassmannian 
$\Gr_{k,n}^{\geq 0}$ 
is the subset of the real Grassmannian $\Gr_{k,n}$ 
consisting of points with all Pl\"ucker coordinates nonnegative.  Following seminal
work of Lusztig \cite{lusztig}, as well as by Fomin and Zelevinsky \cite{FZ}, 
Postnikov initiated the combinatorial study of $\Gr_{k,n}^{\geq 0}$ and its cell decomposition \cite{postnikov}.
Since then the totally nonnegative Grassmannian has found applications in 
diverse contexts such as mirror symmetry \cite{MarshRietsch}, 
soliton solutions to the KP equation \cite{KodamaWilliams}, 
and scattering amplitudes for $\mathcal{N}=4$
supersymmetric Yang-Mills theory \cite{abcgpt}.

Building on \cite{abcgpt}, 
Arkani-Hamed and Trnka \cite{arkani-hamed_trnka} recently introduced a beautiful new 
mathematical object called the \emph{(tree) amplituhedron}, which 
is the image of the totally nonnegative Grassmannian under a particular map.

\begin{defn}\label{defn_amplituhedron}
Let $Z$ be a $(k+m) \times n$ real matrix whose maximal minors are all positive, where
$m\ge 0$ is fixed with $k+m \leq n$. 
Then it induces a map
$$\tilde{Z}:\Gr_{k,n}^{\ge 0} \to \Gr_{k,k+m}$$ defined by 
$$\tilde{Z}(\langle v_1,\dots, v_k \rangle) := \langle Z(v_1),\dots, Z(v_k) \rangle,$$
where $\langle v_1,\dots,v_k\rangle$ is an element of 
$\Gr_{k,n}^{\ge 0}$ written as the span of $k$ basis vectors.\footnote{The
fact that $Z$ has positive maximal minors ensures that $\tilde{Z}$
is well defined
\cite{arkani-hamed_trnka}.  
See \cite[Theorem 4.2]{karp} for a characterization of when 
a matrix $Z$ gives rise to a well-defined map $\tilde{Z}$, and also \cref{grassmann_polytope_remark}.}
The \emph{(tree) amplituhedron} $\mathcal{A}_{n,k,m}(Z)$ is defined to be the image
$\tilde{Z}(\Gr_{k,n}^{\ge 0})$ inside $\Gr_{k,k+m}$.  
\end{defn}

In special cases the amplituhedron recovers familiar objects. If $Z$ is a square matrix, i.e.\ 
$k+m=n$, then $\mathcal{A}_{n,k,m}(Z)$ is isomorphic to 
the totally nonnegative Grassmannian. If $k=1$, then $\mathcal{A}_{n,1,m}(Z)$ is a {\itshape cyclic polytope} in projective space \cite{Sturmfels}.

While the amplituhedron $\mathcal{A}_{n,k,m}(Z)$ is an interesting mathematical object for any $m$, the case of immediate relevance to physics is $m=4$. In this case, it provides a geometric basis for the computation of {\itshape scattering amplitudes} in $\mathcal{N}=4$ supersymmetric Yang-Mills theory. These amplitudes are complex numbers related to the probability of observing a certain scattering process of $n$ particles.
It is expected that such amplitudes can be expressed (modulo higher-order terms) as an integral over the amplituhedron $\mathcal{A}_{n,k,4}(Z)$.  This statement would follow
from the conjecture of Arkani-Hamed and Trnka \cite{arkani-hamed_trnka} that the images of 
a certain collection of $4k$-dimensional cells of $\Gr_{k,n}^{\geq 0}$ provide a ``triangulation'' of the amplituhedron
$\mathcal{A}_{n,k,4}(Z)$.
More specifically, 
the BCFW recurrence
\cite{BCF, BCFW}
provides one way to compute scattering amplitudes.
Translated into the Grassmannian formulation of 
\cite{abcgpt}, the terms in the BCFW recurrence can 
be identified with a collection of $4k$-dimensional cells in 
$\Gr_{k,n}^{\geq 0}$. If the images of these {\itshape BCFW cells} in $\mathcal{A}_{n,k,4}(Z)$ fit together in a nice way, then we can combine the contributions from each term into a single integral over $\mathcal{A}_{n,k,4}(Z)$.

In this paper, we study the amplituhedron $\mathcal{A}_{n,k,1}(Z)$ for $m=1$.
We find that this object is already interesting and non-trivial.
Since $\mathcal{A}_{n,k,1}(Z) \subseteq \Gr_{k, k+1}$, it is convenient to 
take orthogonal complements and work with lines rather than $k$-planes in $\R^{k+1}$.
This leads us to define a related ``B-amplituhedron''
$$\mathcal{B}_{n,k,m}(W):= \{V^{\perp} \cap W : V\in \Gr_{k,n}^{\geq 0} \} \subseteq 
\Gr_m(W),$$
which is homeomorphic to $\mathcal{A}_{n,k,m}(Z)$, where $W$ is the subspace of $\mathbb{R}^n$ spanned by the rows of $Z$ (\cref{sec_B_amplituhedron}). In the context of scattering amplitudes ($m=4$), $W$ is the span of $4$ bosonic variables and $k$ fermionic variables. Building on results of \cite{karp}, we use this reformulation to give a description of the amplituhedron $\mathcal{A}_{n,k,m}(Z)$ in terms of sign variation (\cref{sec_intrinsic_description}).

Modeling the $m=4$ case, we define a BCFW-like recursion in the case $m=1$, which we
use to produce a subset of $k$-dimensional ``BCFW cells'' of $\Gr_{k,n}^{\geq 0}$ (\cref{sec_bcfw}).  
The set of all cells of $\Gr_{k,n}^{\geq 0}$ are in bijection with various combinatorial
objects, including \Le-diagrams and  decorated permutations, so we describe our 
$m=1$ BCFW cells in terms of these objects.  We then show that their images
triangulate the $m=1$ amplituhedron;  more specifically, we show that 
$\mathcal{A}_{n,k,1}(Z)$ is homeomorphic to a $k$-dimensional subcomplex of 
the totally nonnegative Grassmannian $\Gr_{k,n}^{\geq 0}$ (\cref{sec_induced_subcomplex}). See \cref{G24} for $\mathcal{A}_{4,2,1}(Z)$ as a subcomplex of $\Gr_{2,4}^{\ge 0}$.

We also show that 
$\mathcal{A}_{n,k,1}(Z)$ can be identified with
\begin{figure}[t]
\begin{center}
\begin{gather*}
\arraycolsep=2pt
\begin{array}{cccc}
\includegraphics{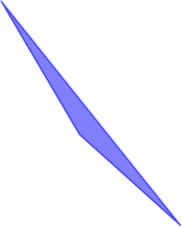} & \includegraphics{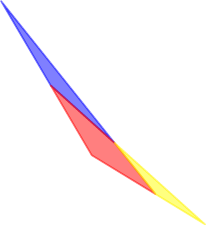} & \includegraphics{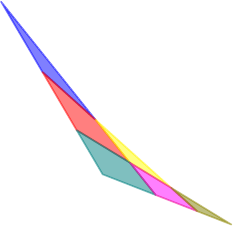} & \includegraphics{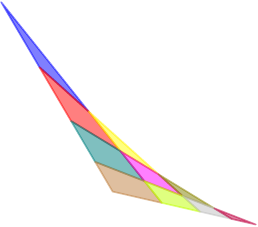} \\
\hspace*{8pt}\mathcal{A}_{3,2,1} & \hspace*{8pt}\mathcal{A}_{4,2,1} & \hspace*{8pt}\mathcal{A}_{5,2,1} & \hspace*{8pt}\mathcal{A}_{6,2,1}
\end{array} \\[2pt]
\arraycolsep=7pt
\begin{array}{ccc}
\hspace*{-4pt}\includegraphics{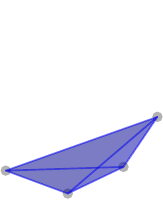} & \includegraphics{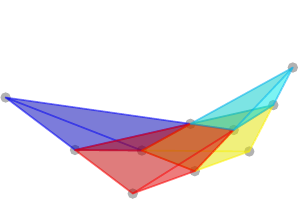} & \includegraphics{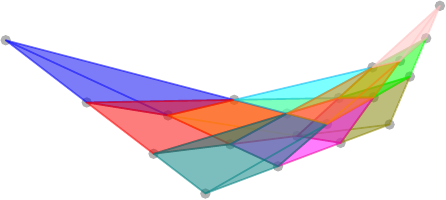} \\
\mathcal{A}_{4,3,1} & \mathcal{A}_{5,3,1} & \mathcal{A}_{6,3,1}
\end{array}
\end{gather*}
\caption{The amplituhedron $\mathcal{A}_{n,k,1}(Z)$ as the complex of bounded faces of a cyclic hyperplane arrangement of $n$ hyperplanes in $\mathbb{R}^k$, for $k=2,3$ and $n\le 6$.}
\label{fig:hyperplanes}
\end{center}
\end{figure}
the complex of bounded faces of a certain hyperplane arrangement of $n$ hyperplanes in $\mathbb{R}^k$, called a {\itshape cyclic hyperplane 
arrangement} (\cref{sec_cyclic_arrangement}). We use this description of the $m=1$ amplituhedron to describe how its cells fit together (\cref{sec_glue}).

It is known that the totally nonnegative Grassmannian has a remarkably simple topology:
it is contractible with boundary
a sphere \cite{RietschWilliams}, and
its poset of cells is Eulerian \cite{Williams}.
While there are not yet any general results in this direction,
calculations of Euler characteristics \cite{anatomy} indicate that the amplituhedron $\mathcal{A}_{n,k,m}(Z)$ is likely also topologically very nice.
Our description of $\mathcal{A}_{n,k,1}(Z)$ as the complex of bounded faces of a hyperplane arrangement, together with a result of Dong \cite{Dong}, implies that the $m=1$
amplituhedron is homeomorphic to a closed ball (\cref{homeomorphic_to_ball}).

Since $k+m \leq n$, the map 
$$\tilde{Z}:\Gr_{k,n}^{\ge 0} \to \Gr_{k,k+m}$$ 
is far from injective in general.
We determine when an arbitrary cell of $\Gr_{k,n}^{\ge 0}$ is mapped injectively by $\tilde{Z}$ into $\mathcal{A}_{n,k,1}(Z)$, and in this case we describe its image in $\mathcal{A}_{n,k,1}(Z)$ (\cref{sec_injectivity}).

Finally, we discuss
to what extent our results hold in the setting of {\itshape Grassmann polytopes} (\cref{sec_grassmann_polytopes}). Grassmann polytopes are generalizations of amplituhedra obtained by relaxing the positivity condition on the matrix $Z$ \cite{lam}.

\textsc{Acknowledgements:}
This paper is an offshoot of a larger ongoing project 
which is joint with Yan Zhang.   We would like to thank him for many
helpful conversations.  We would also like to thank Nima Arkani-Hamed,
Hugh Thomas,
and Jaroslav Trnka
for sharing their results, 
Richard Stanley for providing a reference 
on hyperplane arrangements, Thomas Lam for giving useful comments about \cref{image_dimension}, Pavel Galashin for resolving \cref{lam_problem}, and anonymous referees for their feedback on the paper.

\section{Background on the totally nonnegative Grassmannian}\label{sec_TNN_Grassmannian}

\noindent The {\itshape (real) Grassmannian} $\Gr_{k,n}$ is the space of all
$k$-dimensional subspaces of $\R^n$, for $0\le k \le n$.  An element of
$\Gr_{k,n}$ can be viewed as a $k\times n$ matrix of rank $k$, modulo left
multiplication by invertible $k\times k$ matrices. That is, two
$k\times n$ matrices of rank $k$ represent the same point in $\Gr_{k,n}$ if and only if they
can be obtained from each other by invertible row operations.

Let $[n]$ denote $\{1,\dots,n\}$, and $\binom{[n]}{k}$ the set of all $k$-element subsets of $[n]$. Given $V\in\Gr_{k,n}$ represented by a $k\times n$ matrix $A$, for $I\in \binom{[n]}{k}$ we let $\Delta_I(V)$ be the maximal minor of $A$ located in the column set $I$. The $\Delta_I(V)$ do not depend on our choice of matrix $A$ (up to simultaneous rescaling by a nonzero constant), and are called the {\itshape Pl\"{u}cker coordinates} of $V$.

\begin{defn}[{\cite[Section~3]{postnikov}}]\label{def:positroid}
We say that $V\in\Gr_{k,n}$ is {\itshape totally nonnegative} if $\Delta_I(V)\ge 0$ for all $I\in\binom{[n]}{k}$, and {\itshape totally positive} if $\Delta_I(V) > 0$ for all $I\in\binom{[n]}{k}$. The set of all totally nonnegative $V\in\Gr_{k,n}$ is the {\it totally nonnegative Grassmannian} $\Gr_{k,n}^{\ge 0}$, and the set of all totally positive $V$ is the {\itshape totally positive Grassmannian} $\Gr_{k,n}^{>0}$. For $M\subseteq \binom{[n]}{k}$,
the {\it positroid cell} $S_{M}$ is
the set of $V\in\Gr_{k,n}^{\geq 0}$ with the prescribed collection of Pl\"{u}cker coordinates strictly positive (i.e.\ $\Delta_I(V)>0$ for all $I\in M$), and the remaining Pl\"{u}cker coordinates
equal to zero (i.e.\ $\Delta_J(V)=0$ for all $J\in\binom{[n]}{k}\setminus M$). We call $M$ a \emph{positroid} if $S_M$ is nonempty. We let $Q_{k,n}$ denote the poset on the cells of $\Gr_{k,n}^{\ge 0}$ defined by 
$S_{M} \leq S_{M'}$ if and only if 
$S_M \subseteq \overline{S_{M'}}$.
\end{defn}

\begin{rmk}\label{positive_torus_remark}
There is an action of the ``positive torus'' 
$T_{>0} = \R_{>0}^n$ on $\Gr_{k,n}^{\geq 0}$.  Concretely,
if $A$ is a $k \times n$ matrix representing an element of 
$\Gr_{k,n}^{\geq 0}$, then the positive torus acts on $A$ by rescaling
its columns.  If $\mathbf{t} = (t_1,\dots,t_n) \in T_{>0}$ and 
$A$ represents an element of $S_M$, then 
$\mathbf{t} \cdot A$ also represents an element of $S_M$.
\end{rmk}

The fact that each nonempty $S_{M}$ is a topological cell is due to
Postnikov \cite{postnikov}.  Moreover, it was shown in 
\cite{PSW} that the cells glue together to form a CW 
decomposition of $\Gr_{k,n}^{\geq 0}$.

\subsection{Combinatorial objects parameterizing cells}

In \cite{postnikov}, Postnikov gave several families of combinatorial objects in bijection with 
cells of the totally nonnegative Grassmannian.   In this section we will start by defining
\emph{\Le-diagrams}, \emph{decorated permutations}, 
and equivalence classes of \emph{reduced plabic graphs}, and give
(compatible) bijections among all these objects.  This will give us a 
canonical way to label each positroid by a \Le-diagram, a decorated permutation, and an equivalence class of plabic graphs.

\begin{defn}
A \emph{decorated permutation} of the set $[n]$ is a bijection $\pi : [n] \to [n]$ whose fixed points are colored either black or white. 
We denote a black fixed point by $\pi(i) = \underline{i}$ and 
a white fixed point by $\pi(i) = \overline{i}$.
An \emph{anti-excedance} of the decorated permutation $\pi$ is an element $i \in [n]$ such that either $\pi^{-1}(i) > i$ or 
$\pi(i)=\overline{i}$ (i.e.\ $i$ is a white fixed point). 
\end{defn}

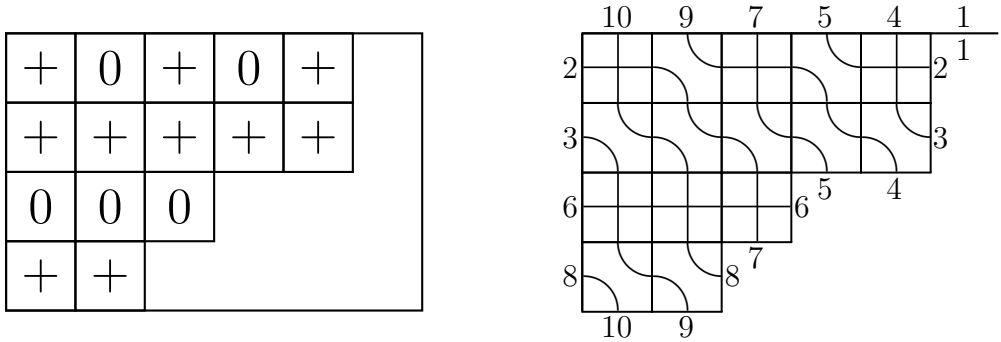
\begin{figure}[ht]
\begin{center}
$$
\begin{tikzpicture}[baseline=(current bounding box.center)]
\pgfmathsetmacro{\scalar}{1.6};
\pgfmathsetmacro{\unit}{\scalar*0.922/1.6};
\draw[thick](0,0)rectangle(6*\unit,-4*\unit);
\foreach \x in {1,...,5}{
\draw[thick](\x*\unit-\unit,0)rectangle(\x*\unit,-\unit);}
\foreach \x in {1,...,5}{
\draw[thick](\x*\unit-\unit,-\unit)rectangle(\x*\unit,-2*\unit);}
\foreach \x in {1,...,3}{
\draw[thick](\x*\unit-\unit,-2*\unit)rectangle(\x*\unit,-3*\unit);}
\foreach \x in {1,...,2}{
\draw[thick](\x*\unit-\unit,-3*\unit)rectangle(\x*\unit,-4*\unit);}
\node[inner sep=0]at(0.5*\unit,-0.5*\unit){\scalebox{\scalar}{$+$}};
\node[inner sep=0]at(1.5*\unit,-0.5*\unit){\scalebox{\scalar}{$0$}};
\node[inner sep=0]at(2.5*\unit,-0.5*\unit){\scalebox{\scalar}{$+$}};
\node[inner sep=0]at(3.5*\unit,-0.5*\unit){\scalebox{\scalar}{$0$}};
\node[inner sep=0]at(4.5*\unit,-0.5*\unit){\scalebox{\scalar}{$+$}};
\node[inner sep=0]at(0.5*\unit,-1.5*\unit){\scalebox{\scalar}{$+$}};
\node[inner sep=0]at(1.5*\unit,-1.5*\unit){\scalebox{\scalar}{$+$}};
\node[inner sep=0]at(2.5*\unit,-1.5*\unit){\scalebox{\scalar}{$+$}};
\node[inner sep=0]at(3.5*\unit,-1.5*\unit){\scalebox{\scalar}{$+$}};
\node[inner sep=0]at(4.5*\unit,-1.5*\unit){\scalebox{\scalar}{$+$}};
\node[inner sep=0]at(0.5*\unit,-2.5*\unit){\scalebox{\scalar}{$0$}};
\node[inner sep=0]at(1.5*\unit,-2.5*\unit){\scalebox{\scalar}{$0$}};
\node[inner sep=0]at(2.5*\unit,-2.5*\unit){\scalebox{\scalar}{$0$}};
\node[inner sep=0]at(0.5*\unit,-3.5*\unit){\scalebox{\scalar}{$+$}};
\node[inner sep=0]at(1.5*\unit,-3.5*\unit){\scalebox{\scalar}{$+$}};
\end{tikzpicture}\qquad\qquad
\begin{tikzpicture}[baseline=(current bounding box.center)]
\pgfmathsetmacro{\unit}{0.922};
\useasboundingbox(-0.5*\unit,0.5*\unit)rectangle(6.5*\unit,-4.5*\unit);
\coordinate (vstep)at(0,-0.24*\unit);
\coordinate (hstep)at(0.17*\unit,0);
\draw[thick](0,0)--(6*\unit,0) (0,0)--(0,-4*\unit);
\node[inner sep=0]at(0,0){\scalebox{1.6}{\begin{ytableau}
\none \\
\none \\
\none \\
\none \\
\none & \none & \none & \none & \none & \cross & \elbow & \cross & \elbow & \cross \\
\none & \none & \none & \none & \none & \elbow & \elbow & \elbow & \elbow & \elbow \\
\none & \none & \none & \none & \none & \cross & \cross & \cross \\
\none & \none & \none & \none & \none & \elbow & \elbow
\end{ytableau}}};
\node[inner sep=0]at($(5.5*\unit,0)+(vstep)$){$1$};
\node[inner sep=0]at($(5*\unit,-0.5*\unit)+(hstep)$){$2$};
\node[inner sep=0]at($(5*\unit,-1.5*\unit)+(hstep)$){$3$};
\node[inner sep=0]at($(4.5*\unit,-2*\unit)+(vstep)$){$4$};
\node[inner sep=0]at($(3.5*\unit,-2*\unit)+(vstep)$){$5$};
\node[inner sep=0]at($(3*\unit,-2.5*\unit)+(hstep)$){$6$};
\node[inner sep=0]at($(2.5*\unit,-3*\unit)+(vstep)$){$7$};
\node[inner sep=0]at($(2*\unit,-3.5*\unit)+(hstep)$){$8$};
\node[inner sep=0]at($(1.5*\unit,-4*\unit)+(vstep)$){$9$};
\node[inner sep=0]at($(0.5*\unit,-4*\unit)+(vstep)$){$10$};
\node[inner sep=0]at($(5.5*\unit,0)-(vstep)$){$1$};
\node[inner sep=0]at($(4.5*\unit,0)-(vstep)$){$4$};
\node[inner sep=0]at($(3.5*\unit,0)-(vstep)$){$5$};
\node[inner sep=0]at($(2.5*\unit,0)-(vstep)$){$7$};
\node[inner sep=0]at($(1.5*\unit,0)-(vstep)$){$9$};
\node[inner sep=0]at($(0.5*\unit,0)-(vstep)$){$10$};
\node[inner sep=0]at($(0,-0.5*\unit)-(hstep)$){$2$};
\node[inner sep=0]at($(0,-1.5*\unit)-(hstep)$){$3$};
\node[inner sep=0]at($(0,-2.5*\unit)-(hstep)$){$6$};
\node[inner sep=0]at($(0,-3.5*\unit)-(hstep)$){$8$};
\end{tikzpicture}
$$
\caption{A \Le -diagram with $\lambda=(5,5,3,2)$, $n=10$, and $k=4$, and its corresponding pipe dream with $\pi = (\underline{1},5,4,9,7,\overline{6},2,10,3,8)$.}
 \label{fig:Le}
\end{center}
\end{figure}

\begin{defn}
Fix $k$ and $n$.  Given a partition $\lambda$,
we 
let
$Y_{\lambda}$ denote the Young diagram associated to $\lambda$.  
A {\it \Le-diagram}
(or Le-diagram) $D$ of shape $\lambda$ and type $(k,n)$
is a Young diagram of shape $Y_\lambda$  
contained in a $k \times (n-k)$ rectangle,
whose boxes are filled with $0$'s and $+$'s in such a way that the
{\it \Le-property} is satisfied:
there is no $0$ which has a $+$ above it in the same column and a $+$ to its
left in the same row.
See \cref{fig:Le} for
an example of a \Le-diagram.
\end{defn}

\begin{lem}[{\cite[Section 20]{postnikov}}]\label{le_permutation_bijection}
The following algorithm is a bijection between \Le-diagrams $D$ of type 
$(k,n)$ and decorated permutations $\pi$ on $[n]$ with exactly $k$ anti-excedances.
\begin{enumerate}
\item Replace each $+$ in the \Le-diagram $D$
 with an elbow joint $\textelbow$, and each $0$ in $D$ with a cross
$\textcross$.
\item Note that the southeast
border of $Y_\lambda$ gives rise to a length-$n$
path from the northeast corner to the southwest corner of the
$k \times (n-k)$ rectangle.  Label the edges of this path with the
numbers $1$ through $n$.
\item Now label the edges of the north and west border of $Y_{\lambda}$
so that opposite horizontal edges and opposite vertical edges
have the same label.
\item View the resulting `pipe dream'
as a permutation $\pi = \pi(D)$ on $[n]$, by following the `pipes' from 
the southeastern border to the northwest border of the Young diagram.
If the pipe originating at label $i$ ends at the label $j$,
we define $\pi(i) := j$.
\item If $\pi(i)=i$ and $i$ labels two horizontal (respectively,
vertical) edges of $Y_{\lambda}$, then $\pi(i):=\underline{i}$
(respectively, $\pi(i):=\overline{i}$).
\end{enumerate}
\end{lem}
\cref{fig:Le} illustrates this procedure.

\begin{defn}
A {\it plabic graph\/}\footnote{``Plabic'' stands for {\itshape planar bi-colored}.}  is an undirected planar graph $G$ drawn inside a disk
(considered modulo homotopy)
with $n$ {\it boundary vertices\/} on the boundary of the disk,
labeled $1,\dots,n$ in clockwise order, as well as some
colored {\it internal vertices\/}.
These  internal vertices
are strictly inside the disk and are each colored either black or white. Moreover, 
each boundary vertex $i$ in $G$ is incident to a single edge.
If a boundary vertex is adjacent to a leaf (a vertex of degree $1$),
we refer to that leaf as a \emph{lollipop}.

A {\it perfect orientation\/} $\O$ of a plabic graph $G$ is a
choice of orientation of each of its edges such that each
black internal vertex $u$ is incident to exactly one edge
directed away from $u$, and each white internal vertex $v$ is incident
to exactly one edge directed towards $v$.
A plabic graph is called {\it perfectly orientable\/} if it admits a perfect orientation.
Let $G_\O$ denote the directed graph associated with a perfect orientation $\O$ of $G$. Since each boundary vertex is incident to a single edge, it is either a {\itshape source} (if it is incident to an outgoing edge) or a {\itshape sink} (if it is incident to an incoming edge) in $G_\O$. The {\it source set\/} $I_\O \subseteq [n]$ is the set of boundary vertices which are sources in $G_\O$.
\end{defn}

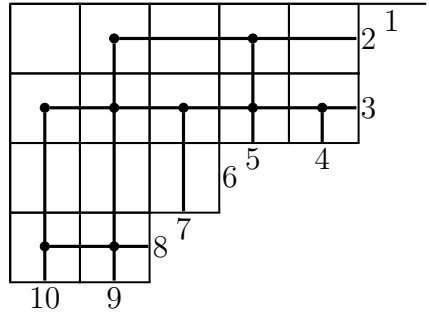
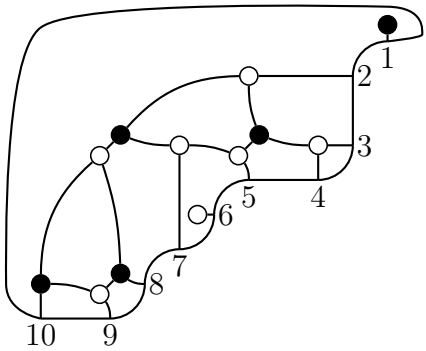
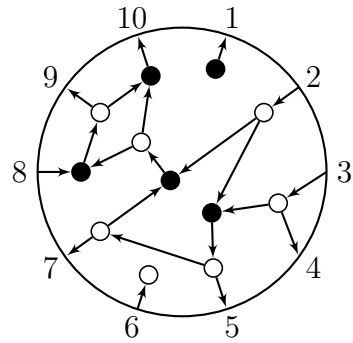
\begin{figure}[ht]
\centering
\subfloat[][The hook diagram $H(D)$.]{
\begin{tikzpicture}[baseline=(current bounding box.center)]
\tikzstyle{out1}=[inner sep=0,minimum size=1.2mm,circle,draw=black,fill=black]
\tikzstyle{in1}=[inner sep=0,minimum size=1.2mm,circle,draw=black,fill=white]
\pgfmathsetmacro{\unit}{0.922};
\useasboundingbox(-0.5*\unit,0.5*\unit)rectangle(6.5*\unit,-4.5*\unit);
\coordinate (vstep)at(0,-0.24*\unit);
\coordinate (hstep)at(0.17*\unit,0);
\coordinate (vepsilon)at(0,-0.02*\unit);
\coordinate (hepsilon)at(0.02*\unit,0);
\draw[thick](0,0)--(6*\unit,0) (0,0)--(0,-4*\unit);
\node[inner sep=0]at(0,0){\scalebox{1.6}{\begin{ytableau}
\none \\
\none \\
\none \\
\none \\
\none & \none & \none & \none & \none & & & & & \\
\none & \none & \none & \none & \none & & & & & \\
\none & \none & \none & \none & \none & & & \\
\none & \none & \none & \none & \none & &
\end{ytableau}}};
\node[inner sep=0]at($(5.5*\unit,0)+(vstep)$){$1$};
\node[inner sep=0]at($(5*\unit,-0.5*\unit)+(hstep)$){$2$};
\node[inner sep=0]at($(5*\unit,-1.5*\unit)+(hstep)$){$3$};
\node[inner sep=0]at($(4.5*\unit,-2*\unit)+(vstep)$){$4$};
\node[inner sep=0]at($(3.5*\unit,-2*\unit)+(vstep)$){$5$};
\node[inner sep=0]at($(3*\unit,-2.5*\unit)+(hstep)$){$6$};
\node[inner sep=0]at($(2.5*\unit,-3*\unit)+(vstep)$){$7$};
\node[inner sep=0]at($(2*\unit,-3.5*\unit)+(hstep)$){$8$};
\node[inner sep=0]at($(1.5*\unit,-4*\unit)+(vstep)$){$9$};
\node[inner sep=0]at($(0.5*\unit,-4*\unit)+(vstep)$){$10$};
\node[inner sep=0](b1)at($(5.5*\unit,0)+(vepsilon)$){};
\node[inner sep=0](b2)at($(5*\unit,-0.5*\unit)+(hepsilon)$){};
\node[inner sep=0](b3)at($(5*\unit,-1.5*\unit)+(hepsilon)$){};
\node[inner sep=0](b4)at($(4.5*\unit,-2*\unit)+(vepsilon)$){};
\node[inner sep=0](b5)at($(3.5*\unit,-2*\unit)+(vepsilon)$){};
\node[inner sep=0](b6)at($(3*\unit,-2.5*\unit)+(hepsilon)$){};
\node[inner sep=0](b7)at($(2.5*\unit,-3*\unit)+(vepsilon)$){};
\node[inner sep=0](b8)at($(2*\unit,-3.5*\unit)+(hepsilon)$){};
\node[inner sep=0](b9)at($(1.5*\unit,-4*\unit)+(vepsilon)$){};
\node[inner sep=0](b10)at($(0.5*\unit,-4*\unit)+(vepsilon)$){};
\node[out1](i12)at($(2*\unit,-1*\unit)+(-0.5*\unit,0.5*\unit)$){};
\node[out1](i14)at($(4*\unit,-1*\unit)+(-0.5*\unit,0.5*\unit)$){};
\node[out1](i21)at($(1*\unit,-2*\unit)+(-0.5*\unit,0.5*\unit)$){};
\node[out1](i22)at($(2*\unit,-2*\unit)+(-0.5*\unit,0.5*\unit)$){};
\node[out1](i23)at($(3*\unit,-2*\unit)+(-0.5*\unit,0.5*\unit)$){};
\node[out1](i24)at($(4*\unit,-2*\unit)+(-0.5*\unit,0.5*\unit)$){};
\node[out1](i25)at($(5*\unit,-2*\unit)+(-0.5*\unit,0.5*\unit)$){};
\node[out1](i41)at($(1*\unit,-4*\unit)+(-0.5*\unit,0.5*\unit)$){};
\node[out1](i42)at($(2*\unit,-4*\unit)+(-0.5*\unit,0.5*\unit)$){};
\path[very thick](b2)edge(i12) (b3)edge(i21) (b8)edge(i41) (i25)edge(b4) (i14)edge(b5) (i23)edge(b7) (i12)edge(b9) (i21)edge(b10);
\end{tikzpicture}\label{fig:hook}}\\
\subfloat[][The plabic graph $G(D)$.]{
\begin{tikzpicture}[baseline=(current bounding box.center)]
\tikzstyle{out1}=[inner sep=0,minimum size=2.4mm,circle,draw=black,fill=black,semithick]
\tikzstyle{in1}=[inner sep=0,minimum size=2.4mm,circle,draw=black,fill=white,semithick]
\pgfmathsetmacro{\unit}{0.922};
\useasboundingbox(-0.5*\unit,0.5*\unit)rectangle(6.5*\unit,-4.5*\unit);
\coordinate (vstep)at(0,-0.24*\unit);
\coordinate (hstep)at(0.17*\unit,0);
\coordinate (vloll)at(0,0.24*\unit);
\coordinate (hloll)at(-0.24*\unit,0);
\coordinate (dstep)at(0.15*\unit,0.15*\unit);
\node[inner sep=0]at($(5.5*\unit,0)+(vstep)$){$1$};
\node[inner sep=0]at($(5*\unit,-0.5*\unit)+(hstep)$){$2$};
\node[inner sep=0]at($(5*\unit,-1.5*\unit)+(hstep)$){$3$};
\node[inner sep=0]at($(4.5*\unit,-2*\unit)+(vstep)$){$4$};
\node[inner sep=0]at($(3.5*\unit,-2*\unit)+(vstep)$){$5$};
\node[inner sep=0]at($(3*\unit,-2.5*\unit)+(hstep)$){$6$};
\node[inner sep=0]at($(2.5*\unit,-3*\unit)+(vstep)$){$7$};
\node[inner sep=0]at($(2*\unit,-3.5*\unit)+(hstep)$){$8$};
\node[inner sep=0]at($(1.5*\unit,-4*\unit)+(vstep)$){$9$};
\node[inner sep=0]at($(0.5*\unit,-4*\unit)+(vstep)$){$10$};
\coordinate(v05h)at(5.5*\unit,0*\unit);
\coordinate(v0h5)at(5*\unit,-0.5*\unit);
\coordinate(v1h5)at(5*\unit,-1.5*\unit);
\coordinate(v24h)at(4.5*\unit,-2*\unit);
\coordinate(v23h)at(3.5*\unit,-2*\unit);
\coordinate(v2h3)at(3*\unit,-2.5*\unit);
\coordinate(v32h)at(2.5*\unit,-3*\unit);
\coordinate(v3h2)at(2*\unit,-3.5*\unit);
\coordinate(v41h)at(1.5*\unit,-4*\unit);
\coordinate(v40h)at(0.5*\unit,-4*\unit);
\node[out1](v0e5h)at($(5.5*\unit,0*\unit)+(vloll)$){};
\node[in1](v0h3h)at(3.5*\unit,-0.5*\unit){};
\node[inner sep=0](v0h2)at(2*\unit,-0.5*\unit){};
\node[inner sep=0](v11h)at(1.5*\unit,-1*\unit){};
\node[in1](v1h4h)at(4.5*\unit,-1.5*\unit){};
\node[out1](v1h3hout1)at($(3.5*\unit,-1.5*\unit)+(dstep)$){};
\node[in1](v1h3hin1)at($(3.5*\unit,-1.5*\unit)-(dstep)$){};
\node[in1](v1h2h)at(2.5*\unit,-1.5*\unit){};
\node[out1](v1h1hout1)at($(1.5*\unit,-1.5*\unit)+(dstep)$){};
\node[in1](v1h1hin1)at($(1.5*\unit,-1.5*\unit)-(dstep)$){};
\node[inner sep=0](v1h1)at(1*\unit,-1.5*\unit){};
\node[inner sep=0](v20h)at(0.5*\unit,-2*\unit){};
\node[in1](v2h3e)at($(3*\unit,-2.5*\unit)+(hloll)$){};
\node[out1](v3h1hout1)at($(1.5*\unit,-3.5*\unit)+(dstep)$){};
\node[in1](v3h1hin1)at($(1.5*\unit,-3.5*\unit)-(dstep)$){};
\node[out1](v3h0h)at(0.5*\unit,-3.5*\unit){};
\node[inner sep=0](border0)at(6*\unit,0.1*\unit){};
\node[inner sep=0](border1)at(5.5*\unit,0.5*\unit){};
\node[inner sep=0](border2)at(0.5*\unit,0.2*\unit){};
\node[inner sep=0](border3)at(0,-3.5*\unit){};
\path[thick](v05h)edge[bend right=45](v0h5) (v0h5)edge(v1h5) (v1h5)edge[bend left=45](v24h) (v24h)edge(v23h) (v23h)edge[bend right=45](v2h3) (v2h3)edge[bend left=45](v32h) (v32h)edge[bend right=45](v3h2) (v3h2)edge[bend left=45](v41h) (v41h)edge(v40h) (v05h)edge(v0e5h) (v0h5)edge(v0h3h) (v0h3h)edge[bend right=24](v1h1hout1) (v0h3h)edge[bend right=10](v1h3hout1) (v1h5)edge(v1h4h) (v1h4h)edge(v24h) (v1h4h)edge[bend left=10](v1h3hout1) (v1h3hout1)edge(v1h3hin1) (v1h3hin1)edge[bend left=16](v23h) (v1h3hin1)edge[bend right=10](v1h2h) (v1h2h)edge[bend left=10](v1h1hout1) (v1h2h)edge(v32h) (v1h1hout1)edge(v1h1hin1) (v1h1hin1)edge[bend right=24](v3h0h) (v1h1hin1)edge[bend left=8](v3h1hout1) (v2h3)edge(v2h3e) (v3h2)edge[bend left=16](v3h1hout1) (v3h1hout1)edge(v3h1hin1) (v3h1hin1)edge[bend right=10](v3h0h) (v3h1hin1)edge[bend left=16](v41h) (v3h0h)edge(v40h);
\draw[thick]plot[smooth,tension=0.36]coordinates{(v05h) (border0) (border1) (border2)  (border3) (v40h)};
\end{tikzpicture}\label{fig:plabic}}
\subfloat[][The plabic graph $G(D)$ redrawn and perfectly oriented.]{
\qquad\qquad\begin{tikzpicture}[baseline=(current bounding box.center)]
\tikzstyle{out1}=[inner sep=0,minimum size=2.4mm,circle,draw=black,fill=black,semithick]
\tikzstyle{in1}=[inner sep=0,minimum size=2.4mm,circle,draw=black,fill=white,semithick]
\pgfmathsetmacro{\radius}{1.92};
\draw[thick](0,0)circle[radius=\radius];
\foreach \x in {1,...,10}{
\node[inner sep=0](b\x)at(108-\x*36:\radius){};
\node at(108-\x*36:\radius+.24){$\x$};}
\node[out1](i1)at(72:0.75*\radius){};
\node[in1](i2)at(36:0.7*\radius){};
\node[in1](i3)at(-18:0.7*\radius){};
\node[out1](i4)at(-144:0.1*\radius){};
\node[in1](i5)at(-72:0.7*\radius){};
\node[in1](i6)at(-108:0.75*\radius){};
\node[in1](i7)at(-144:0.7*\radius){};
\node[out1](i8)at(180:0.7*\radius){};
\node[in1](i9)at(144:0.7*\radius){};
\node[out1](i10)at(108:0.7*\radius){};
\node[in1](i11)at(144:0.35*\radius){};
\node[out1](i12)at(-54:0.35*\radius){};
\path[-latex',thick](i1)edge(b1) (b2.center)edge(i2) (i2)edge(i12) (i2)edge(i4) (b3.center)edge(i3) (i3)edge(i12) (i3)edge(b4) (i12)edge(i5) (i5)edge(b5) (i5)edge(i7) (b6.center)edge(i6) (i7)edge(i4) (i7)edge(b7) (b8.center)edge(i8) (i8)edge(i9) (i9)edge(b9) (i9)edge(i10) (i10)edge(b10) (i4)edge(i11) (i11)edge(i8) (i11)edge(i10);
\end{tikzpicture}\qquad\qquad\label{fig:perfect_orientation}}
\caption{The hook diagram and plabic graph associated to the \Le -diagram $D$ from \cref{fig:Le}.}
\label{fig:orientation}
\end{figure}

\cref{fig:orientation} shows a plabic graph with a perfect orientation. In that example, $I_{\O} = \{2,3,6,8\}$.

All perfect orientations of a fixed plabic graph
$G$ have source sets of the same size $k$, where
$k-(n-k) = \sum \mathrm{color}(v)\cdot(\deg(v)-2)$. 
Here the sum is over all internal vertices $v$, where $\mathrm{color}(v) = 1$ if $v$ is black, and $\mathrm{color}(v) = -1$ if $v$ is white;
see~\cite[Lemma 9.4]{postnikov}.  In this case we say that $G$ is of {\it type\/} $(k,n)$.

The following construction of Postnikov \cite[Sections 6 and 20]{postnikov} associates a perfectly orientable plabic graph to any \Le -diagram.

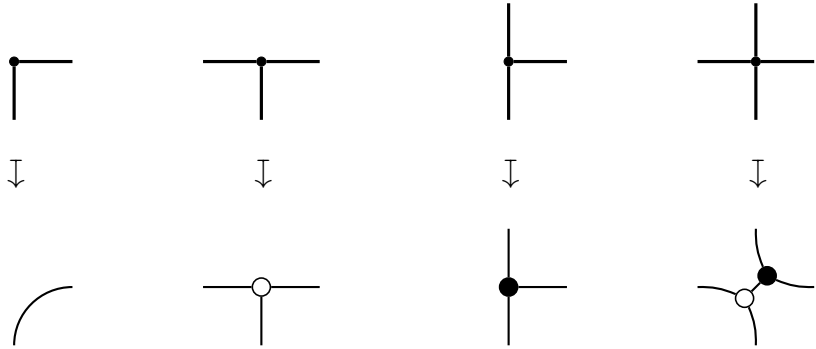
\begin{figure}[ht]
\begin{center}
\begin{tikzpicture}[baseline=(current bounding box.center)]
\tikzstyle{out1}=[inner sep=0,minimum size=2.4mm,circle,draw=black,fill=black,semithick]
\tikzstyle{in1}=[inner sep=0,minimum size=2.4mm,circle,draw=black,fill=white,semithick]
\tikzstyle{hookvertex}=[inner sep=0,minimum size=1.2mm,circle,draw=black,fill=black]
\pgfmathsetmacro{\unit}{1.5};
\pgfmathsetmacro{\side}{0.8};
\node[hookvertex](uc)at(0,\unit){};
\node[inner sep=0](un)at(0,\unit+\side){};
\node[inner sep=0](ue)at(\side,\unit){};
\node[inner sep=0](us)at(0,\unit-\side){};
\node[inner sep=0](uw)at(-\side,\unit){};
\node[inner sep=0](dc)at(0,-\unit){};
\node[inner sep=0](dn)at(0,-\unit+\side){};
\node[inner sep=0](de)at(\side,-\unit){};
\node[inner sep=0](ds)at(0,-\unit-\side){};
\node[inner sep=0](dw)at(-\side,-\unit){};
\node[inner sep=0]at(0,0){\rotatebox[origin=c]{-90}{$\mapsto$}};
\path[very thick](uc)edge(ue) (uc)edge(us);
\path[thick](de)edge[bend right=45](ds);
\end{tikzpicture}\qquad\qquad
\begin{tikzpicture}[baseline=(current bounding box.center)]
\tikzstyle{out1}=[inner sep=0,minimum size=2.4mm,circle,draw=black,fill=black,semithick]
\tikzstyle{in1}=[inner sep=0,minimum size=2.4mm,circle,draw=black,fill=white,semithick]
\tikzstyle{hookvertex}=[inner sep=0,minimum size=1.2mm,circle,draw=black,fill=black]
\pgfmathsetmacro{\unit}{1.5};
\pgfmathsetmacro{\side}{0.8};
\node[hookvertex](uc)at(0,\unit){};
\node[inner sep=0](un)at(0,\unit+\side){};
\node[inner sep=0](ue)at(\side,\unit){};
\node[inner sep=0](us)at(0,\unit-\side){};
\node[inner sep=0](uw)at(-\side,\unit){};
\node[in1](dc)at(0,-\unit){};
\node[inner sep=0](dn)at(0,-\unit+\side){};
\node[inner sep=0](de)at(\side,-\unit){};
\node[inner sep=0](ds)at(0,-\unit-\side){};
\node[inner sep=0](dw)at(-\side,-\unit){};
\node[inner sep=0]at(0,0){\rotatebox[origin=c]{-90}{$\mapsto$}};
\path[very thick](uc)edge(ue) (uc)edge(us) (uc)edge(uw);
\path[thick](dc)edge(de) (dc)edge(ds) (dc)edge(dw);
\end{tikzpicture}\qquad\qquad
\begin{tikzpicture}[baseline=(current bounding box.center)]
\tikzstyle{out1}=[inner sep=0,minimum size=2.4mm,circle,draw=black,fill=black,semithick]
\tikzstyle{in1}=[inner sep=0,minimum size=2.4mm,circle,draw=black,fill=white,semithick]
\tikzstyle{hookvertex}=[inner sep=0,minimum size=1.2mm,circle,draw=black,fill=black]
\pgfmathsetmacro{\unit}{1.5};
\pgfmathsetmacro{\side}{0.8};
\node[hookvertex](uc)at(0,\unit){};
\node[inner sep=0](un)at(0,\unit+\side){};
\node[inner sep=0](ue)at(\side,\unit){};
\node[inner sep=0](us)at(0,\unit-\side){};
\node[inner sep=0](uw)at(-\side,\unit){};
\node[out1](dc)at(0,-\unit){};
\node[inner sep=0](dn)at(0,-\unit+\side){};
\node[inner sep=0](de)at(\side,-\unit){};
\node[inner sep=0](ds)at(0,-\unit-\side){};
\node[inner sep=0](dw)at(-\side,-\unit){};
\node[inner sep=0]at(0,0){\rotatebox[origin=c]{-90}{$\mapsto$}};
\path[very thick](uc)edge(ue) (uc)edge(us) (uc)edge(un);
\path[thick](dc)edge(dn) (dc)edge(de) (dc)edge(ds);
\end{tikzpicture}\qquad\qquad
\begin{tikzpicture}[baseline=(current bounding box.center)]
\tikzstyle{out1}=[inner sep=0,minimum size=2.4mm,circle,draw=black,fill=black,semithick]
\tikzstyle{in1}=[inner sep=0,minimum size=2.4mm,circle,draw=black,fill=white,semithick]
\tikzstyle{hookvertex}=[inner sep=0,minimum size=1.2mm,circle,draw=black,fill=black]
\pgfmathsetmacro{\unit}{1.5};
\pgfmathsetmacro{\side}{0.8};
\pgfmathsetmacro{\shift}{0.15};
\node[hookvertex](uc)at(0,\unit){};
\node[inner sep=0](un)at(0,\unit+\side){};
\node[inner sep=0](ue)at(\side,\unit){};
\node[inner sep=0](us)at(0,\unit-\side){};
\node[inner sep=0](uw)at(-\side,\unit){};
\node[out1](dout1)at(\shift,-\unit+\shift){};
\node[in1](din1)at(-\shift,-\unit-\shift){};
\node[inner sep=0](dn)at(0,-\unit+\side){};
\node[inner sep=0](de)at(\side,-\unit){};
\node[inner sep=0](ds)at(0,-\unit-\side){};
\node[inner sep=0](dw)at(-\side,-\unit){};
\node[inner sep=0]at(0,0){\rotatebox[origin=c]{-90}{$\mapsto$}};
\path[very thick](uc)edge(ue) (uc)edge(us) (uc)edge(uw) (uc)edge(un);
\path[thick](dout1)edge[bend left=12](dn) (dout1)edge[bend right=12](de) (dout1)edge(din1) (din1)edge[bend left=12](ds) (din1)edge[bend right=12](dw);
\end{tikzpicture}
\caption{Local substitutions for getting the plabic graph $G(D)$ from the hook diagram $H(D)$.}
\label{fig:local}
\end{center}
\end{figure}

\begin{defn}\label{def:Le-plabic}
Let $D$ be a \Le-diagram and $\pi$ its decorated permutation. 
Delete the $0$'s of $D$, and replace each $+$ with a vertex.
From each vertex we construct a hook which goes east
and south, to the border of the Young diagram.
The resulting diagram is called the {\itshape hook diagram} $H(D)$. After
replacing the edges along the southeast border of the Young diagram with  boundary vertices
labeled by $1,\dots,n$, we obtain a planar graph in a disk, with $n$ boundary vertices and one internal vertex
for each $+$ of $D$.  Then
we replace the local region around each internal vertex
as in \cref{fig:local}, and add a black (respectively, white) lollipop for each black (respectively, white) fixed point of $\pi$. This gives rise to a plabic graph which we call $G(D)$. By orienting the edges of $G(D)$ down and to the left, we obtain a perfect orientation.
\end{defn}

\cref{fig:hook} depicts the hook diagram $H(D)$ corresponding to the 
\Le-diagram $D$ from \cref{fig:Le},
and \cref{fig:plabic} shows the corresponding plabic graph $G(D)$.

More generally, each \Le-diagram $D$ 
is associated with a family of \emph{reduced plabic graphs}
consisting of $G(D)$ together with other plabic graphs which can be obtained
from $G(D)$ by certain \emph{moves}; see \cite[Section 12]{postnikov}.

From the plabic graph constructed in \cref{def:Le-plabic} 
(and more generally from a reduced plabic graph $G$), 
one may read off the corresponding decorated permutation $\pi_G$ as follows.

\begin{defn}\label{def:rules}
Let $G$ be a reduced plabic graph as above
with boundary vertices $1,\dots, n$.  
The \emph{trip from $i$} is the path 
obtained by starting from $i$ and traveling along 
edges of $G$ according to the rule that each time we reach an internal black vertex we turn (maximally)
right, and each time we reach an internal white vertex we turn (maximally) left.
This trip ends at some boundary vertex ${\pi(i)}$. By \cite[Section 13]{postnikov}, the fact that $G$ is reduced implies that each fixed point of $\pi$ is attached to a lollipop; we color each fixed point by the color of its lollipop. 
In this way we obtain
the \emph{decorated permutation}
$\pi_G=(\pi(1),\dots,\pi(n))$
 of $G$.
\end{defn}

We invite the reader to verify that when we apply these rules to plabic graph $G$ of \cref{fig:plabic}, we obtain the decorated permutation $\pi_G = (\underline{1},5,4,9,7,\overline{6},2,10,3,8)$.

\subsection{Matroids and positroids}

A matroid is a combinatorial object which unifies several notions of independence. Among the many equivalent ways of defining a matroid we will adopt the point of view of bases, which is one of the most convenient for the study of positroids. We refer the reader to \cite{Oxley} for an in-depth introduction to matroid theory.

\begin{defn}\label{defn_matroid}
Let $E$ be a finite set. A \emph{matroid with ground set $E$} is a subset $M\subseteq 2^E$ satisfying the \emph{basis exchange axiom}:
\begin{center}
 if $B, B' \in M$ and $b\in B\setminus B'$, then there exists $b' \in B'\setminus B$ such that $(B \setminus \{b\}) \cup \{b'\} \in M$.
\end{center}
The elements of $M$ are called \emph{bases}. All bases of $M$ have the same size, called the \emph{rank} of $M$. We say that $i\in E$ is a \emph{loop} if $i$ is contained in no basis of $M$, and a \emph{coloop} if $i$ is contained in every basis of $M$.
\end{defn}

\begin{eg}\label{eg_matroid}
Let $A$ be a $k\times n$ matrix of rank $k$ 
with entries in a field $\mathbb{F}$. Then the subsets $B\in\binom{[n]}{k}$ such that the columns $B$ of $A$ are linearly independent form the bases of a matroid $M(A)$ with rank $k$ and ground set $[n]$. In terms of the Grassmannian, the rows of $A$ span an element of $\Gr_{k,n}(\mathbb{F})$, whose nonzero Pl\"{u}cker coordinates are indexed by $M(A)$. For example, the matrix
$$
A := \begin{bmatrix}1 & 0 & 2 & 0 & 0 \\ 0 & 1 & 3 & 0 & 1\end{bmatrix}
$$
over $\mathbb{F} := \mathbb{Q}$ gives rise to the matroid $M(A) = \{\{1,2\}, \{1,3\}, \{1,5\}, \{2,3\}, \{3,5\}\}$ with ground set $[5]$. In this example, $4$ is a loop of $M(A)$, and $M(A)$ has no coloops.

Matroids arising in this way are called \emph{representable} (over $\mathbb{F}$).
\end{eg}

\begin{eg}\label{eg_matroid2}
Given $k\in\mathbb{N}$ and a finite set $E$, all $k$-subsets of $E$ form a matroid $M$, called the {\itshape uniform matroid} (of rank $k$ with ground set $E$). Note that $M$ can be represented over any infinite field, by a generic $k\times n$ matrix.
\end{eg}

Recall the definition of a positroid from \cref{def:positroid}. In the language of matroid theory, a positroid is a matroid representable by an element of the totally nonnegative Grassmannian. Every perfectly orientable plabic graph gives rise to a positroid as follows.
\begin{defn}[{\cite[Proposition 11.7]{postnikov}}]\label{prop:perf}
Let $G$ be a plabic graph of type $(k,n)$.
Then we have a positroid 
$M_G$ on $[n]$ whose bases are precisely
$$
\{I_{\O} : \O \text{ is a perfect orientation of }G\},
$$
where $I_\O$ is the set of sources of $\O$.

If $D$ is a \Le -diagram contained in a $k\times (n-k)$ rectangle, we let $M(D)$ denote the positroid $M_{G(D)}$ of the plabic graph $G(D)$ from \cref{def:Le-plabic}.
\end{defn}
Postnikov \cite[Theorem 17.1]{postnikov} showed that every positroid can be realized as $M(D)$ for some \Le-diagram $D$. We observe that we can describe the loops and coloops of $M(D)$ in terms of $D$ as follows: $i$ is a loop if and only if $i$ labels a horizontal step whose column contains only $0$'s, and $i$ is a coloop if and only if $i$ labels a vertical step in the southeast border whose row contains only $0$'s.

We introduce some further notions from matroid theory which we will use later: duality, direct sum, connectedness, restriction, and a partial order.
\begin{defn}\label{defn_dual}
Let $M$ be a matroid with ground set $E$. Then $\{E\setminus B : B\in M\}$ is the set of bases of a matroid $M^*$ with ground set $E$, called the {\itshape dual} of $M$.
\end{defn}
See \cite[Section 2]{Oxley} for a proof that $M^*$ is indeed a matroid. We make the following observations about matroid duality:
\begin{itemize}
\item $(M^*)^* = M$;
\item the ranks of $M$ and $M^*$ sum to $|E|$;
\item $i\in E$ is a loop of $M$ if and only if $i$ is a coloop of $M^*$;
\item if $E = [n]$, then $M$ is a positroid if and only if $M^*$ is a positroid (see \cref{dual_translation}(ii)).
\end{itemize}
\begin{eg}\label{eg_matroid_dual}
Let $A := \begin{bmatrix}1 & 0 & 2 & 0 & 0 \\ 0 & 1 & 3 & 0 & 1\end{bmatrix}
$, as in \cref{eg_matroid}. Then
$$
M(A)^* = \{\{3,4,5\},\{2,4,5\},\{2,3,4\},\{1,4,5\},\{1,2,4\}\},
$$
and is represented by the matrix
$$
\begin{bmatrix}
0 & 0 & 0 & 1 & 0 \\
0 & 1 & 0 & 0 & -1 \\
-2 & 0 & 1 & 0 & -3
\end{bmatrix},
$$
whose rows are orthogonal to the rows of $A$.
\end{eg}

\begin{defn}\label{defn_direct_sum}
Let $M$ and $N$ be matroids with ground sets $E$ and $F$, respectively. The \emph{direct sum} $M\oplus N$ is the matroid with ground set $E\sqcup F$, and bases $\{B\sqcup C : B\in M, C\in N\}$. The rank of $M\oplus N$ is the sum of the ranks of $M$ and $N$.

A matroid is \emph{connected} if we cannot write it as the direct sum of two matroids whose ground sets are nonempty. Any matroid $M$ can be written uniquely (up to permuting the summands) as the direct sum of connected matroids, whose ground sets are the \emph{connected components} of $M$; see \cite[Corollary 4.2.9]{Oxley}.
\end{defn}

\begin{eg}\label{eg_direct_sum}
Consider the matrix $A := \begin{bmatrix}1 & 0 & 0 & 1 \\ 0 & 1 & 1 & 0\end{bmatrix}$ and its associated matroid $M(A) = \{\{1,2\},\{1,3\},\{2,4\},\{3,4\}\}$. We have $M(A) = M_1\oplus M_2$, where $M_1$ is the uniform matroid of rank $1$ with ground set $\{1,4\}$, and $M_2$ is the uniform matroid of rank $1$ with ground set $\{2,3\}$. Since $M_1$ and $M_2$ are connected, the connected components of $M(A)$ are $\{1,4\}$ and $\{2,3\}$. In particular, $M(A)$ is disconnected.
\end{eg}

\begin{defn}\label{defn_restriction}
Let $M$ be a matroid with ground set $E$. For a subset $F\subseteq E$, the {\itshape restriction $M|_F$} of $M$ to $F$ is the matroid with ground set $F$ whose bases are the inclusion-maximal sets of $\{B\cap F : B\in M\}$; see \cite[p.\ 20]{Oxley}.
\end{defn}
For example, if $M$ and $N$ are matroids with ground sets $E$ and $F$, respectively, then the restriction of $M\oplus N$ to $E$ is $M$, and the restriction of $M\oplus N$ to $F$ is $N$.

\begin{figure}[ht]
\begin{center}
\includegraphics[width=6.4in]{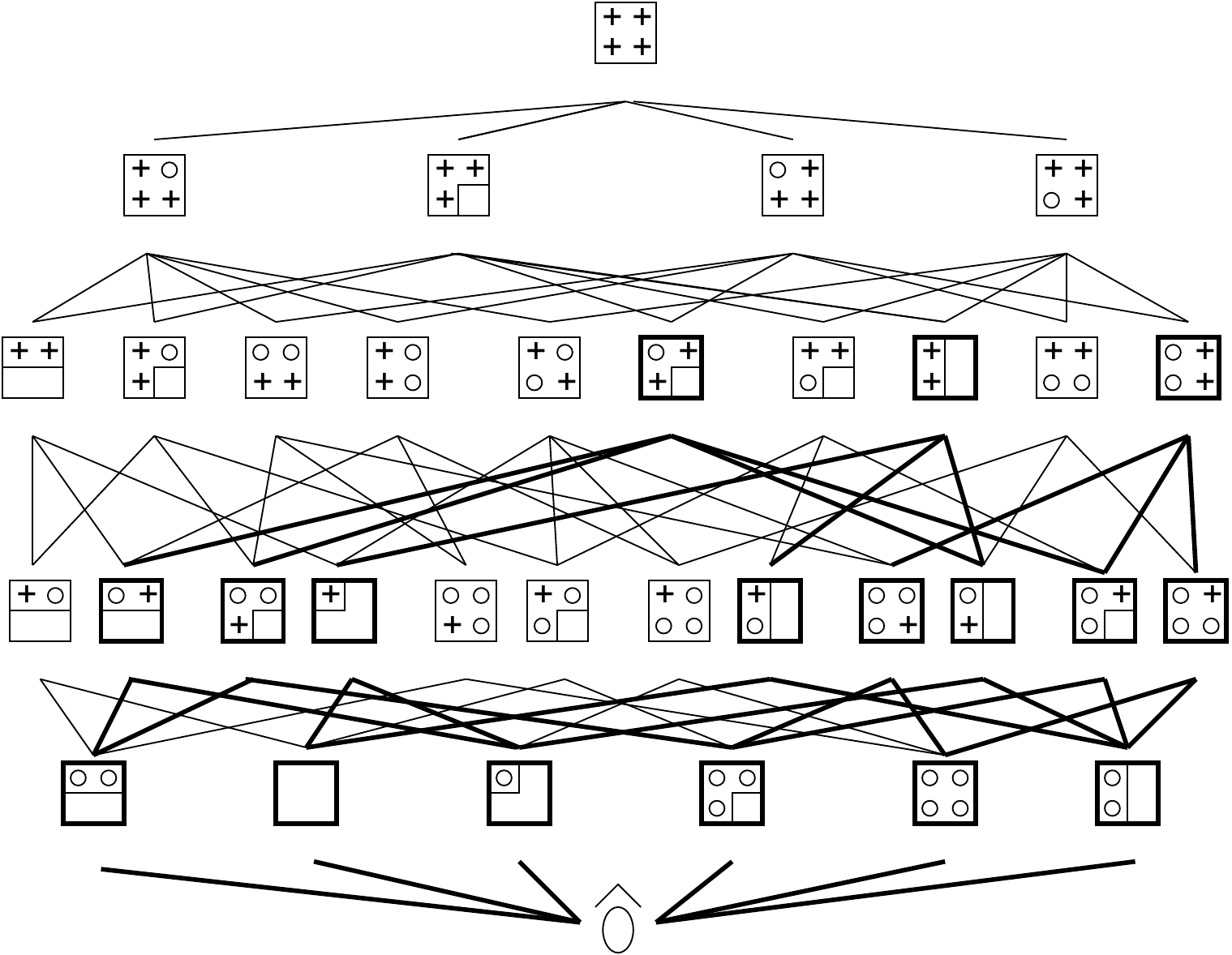}
\end{center}
\caption{The poset $Q_{2,4}$ of cells of $\Gr_{2,4}^{\ge 0}$, where each cell 
is identified with the corresponding \Le -diagram.  The bold
edges indicate  the subcomplex (an induced subposet) 
which gets identified with 
the amplituhedron $\mathcal{A}_{4,2,1}(Z)$.}
\label{G24}
\end{figure}

\begin{defn}\label{defn_weak_maps}
We define a partial order on matroids of rank $k$ with ground set $E$ as follows: $M'\le M$ if and only if every basis of $M'$ is a basis of $M$. (In the matroid theory literature, one says that the identity map on $E$ is a {\itshape weak map} from $M$ to $M'$.)
\end{defn}
Note that $M'\le M$ if and only if $M'^* \le M^*$. This partial order, restricted to positroids of rank $k$ with ground set $[n]$, recovers the poset 
$Q_{k,n}$ of cells of $\Gr_{k,n}^{\ge 0}$ coming from containment
of closures \cite[Section 17]{postnikov}. The poset $Q_{2,4}$ for $\Gr_{2,4}^{\ge 0}$ is shown in \cref{G24}.
\begin{rmk}
All bijections that we have defined in this section are compatible.  This gives us
a canonical way to label each positroid of rank $k$ with ground set $[n]$ by a set of bases,
a decorated permutation, 
a \Le -diagram, and an equivalence
class of reduced plabic graphs. The partial order on positroids (\cref{defn_weak_maps}) gives a partial order on these other objects (of type $(k,n)$).
\end{rmk}

\section{A complementary view of the amplituhedron \texorpdfstring{$\mathcal{A}_{n,k,m}$}{A(n,k,m)}}\label{sec_complementary_view}

\subsection{Background on sign variation}\label{sec_background}

\begin{defn}\label{defn_var}
Given $v\in\mathbb{R}^n$, let $\var(v)$ be the number of times $v$ changes sign, when viewed as a sequence of $n$ numbers and ignoring any zeros. We use the convention $\var(0) := -1$. We also define
$$
\overline{\var}(v) := \max\{\var(w) : \text{$w\in\mathbb{R}^n$ such that $w_i = v_i$ for all $i\in [n]$ with $v_i\neq 0$}\},
$$
i.e.\ $\overline{\var}(v)$ is the maximum number of times $v$ changes sign after we choose a sign for each zero component.
\end{defn}
For example, if $v := (4, -1, 0, -2)\in\mathbb{R}^4$, then $\var(v) = 1$ and $\overline{\var}(v) = 3$.

We now explain how $\var(\cdot)$ and $\overline{\var}(\cdot)$ are dual to each other.
\begin{defn}\label{defn_alt}
We define $\alt:\mathbb{R}^n\to\mathbb{R}^n$ by $\alt(v) := (v_1, -v_2, v_3, -v_4, \dots, (-1)^{n-1}v_n)$ for $v\in\mathbb{R}^n$. If $S\subseteq\mathbb{R}^n$, we let $\alt(S)$ denote $\{\alt(v) : v\in S\}$.
\end{defn}

\begin{lem}[Duality via $\alt$]\label{dual_translation}~\\
(i)\footnote{This result is stated without proof as \cite[Equation II.(67)]{gantmakher_krein_50}. See \cite[Equation (5.1)]{ando} for a proof.} \cite{gantmakher_krein_50} We have $\var(v) + \overline{\var}(\alt(v)) = n-1$ for all $v\in\mathbb{R}^n\setminus\{0\}$. \\
(ii)\footnote{The earliest reference we found for this result is \cite[Section 7]{hochster}, where it appears without proof. Hochster says this result ``was basically known to Hilbert.'' The idea is that if $[I_k | A]$ is a $k\times n$ matrix whose rows span $V\in\Gr_{k,n}$, where $A$ is a $k\times (n-k)$ matrix, then $V^\perp$ is the row span of the matrix $[A^T | -I_{n-k}]$. This appears implicitly in \cite[Equation (14)]{hilbert}, and more explicitly in \cite[Theorem 2.2.8]{Oxley} and \cite[Proposition 3.1(i)]{musiker_reiner}.} \cite{hilbert,hochster} Given $V\in\Gr_{k,n}$, let $V^\perp\in\Gr_{n-k,n}$ be the orthogonal complement of $V$. Then $V$ and $\alt(V^\perp)$ have the same Pl\"{u}cker coordinates, i.e.\ $\Delta_I(V) = \Delta_{[n]\setminus I}(\alt(V^\perp))$ for all $I\in\binom{[n]}{k}$.
\end{lem}
Note that part (ii) above implies that a subspace $V$ is totally nonnegative if and only if $\alt(V^\perp)$ is totally nonnegative, and totally positive if and only if $\alt(V^\perp)$ is totally positive.

The following result of Gantmakher and Krein, which characterizes totally nonnegative and totally positive subspaces in terms of sign variation, will be essential for us.
\begin{thm}[{\cite[Theorems V.3, V.7, V.1, V.6]{gantmakher_krein_50}}]\label{gantmakher_krein}
Let $V\in\Gr_{k,n}$. \\
(i) $V\in\Gr_{k,n}^{\ge 0}\iff\var(v)\le k-1\text{ for all }v\in V\iff\overline{\var}(w)\ge k\text{ for all }w\in V^\perp\setminus\{0\}$. \\
(ii) $V\in\Gr_{k,n}^{>0}\iff\overline{\var}(v)\le k-1\text{ for all }v\in V\setminus\{0\}\iff\var(w)\ge k\text{ for all }w\in V^\perp\setminus\{0\}$.
\end{thm}
\begin{cor}\label{cor:disjoint}
If $V\in \Gr_{k,n}^{\ge 0}$ and $W\in \Gr_{r,n}^{>0}$, where $r \geq k$,
then $V \cap W^\perp = \{0\}$.
\end{cor}

\begin{pf}
By \cref{gantmakher_krein}, $v\in V$ implies that 
$\var(v) \leq k-1$.  And $w\in W^\perp \setminus\{0\}$ implies that 
$\var(w) \geq r \geq k$.  Therefore $V \cap W^\perp = \{0\}$.
\end{pf}

We will also need to know which sign vectors appear in elements of $\Gr_{k,n}^{>0}$.
\begin{defn}
For $t\in\mathbb{R}$ we define
$$
\sign(t) := \begin{cases}
0, & \text{if $t = 0$}, \\
+, & \text{if $t > 0$}, \\
-, & \text{if $t < 0$}.
\end{cases}
$$
(We will sometimes use $1$ and $-1$ in place of $+$ and $-$.) Given $v\in\mathbb{R}^n$, define the {\itshape sign vector} $\sign(v)\in\{0,+,-\}^n$ of $v$ by $\sign(v)_i := \sign(v_i)$ for $i\in [n]$. For example, $\sign(5,0,-1,2) = (+, 0, -, +) = (1, 0, -1, 1)$. If $S\subseteq\mathbb{R}^n$, we let $\sign(S)$ denote $\{\sign(v) : v\in S\}$.
\end{defn}

\begin{lem}\label{positive_covectors}
Suppose that $V\in\Gr_{k,n}^{>0}$ with orthogonal complement $V^\perp$. \\
(i) $\sign(V) = \{\sigma\in\{0,+,-\}^n : \overline{\var}(\sigma) \le k-1\}\cup\{0\}$. \\
(ii) $\sign(V^\perp) = \{\sigma\in\{0,+,-\}^n : \var(\sigma)\ge k\}\cup\{0\}$.
\end{lem}
This essentially follows from \cref{gantmakher_krein}, \cite[Proposition 9.4.1]{bjorner_las_vergnas_sturmfels_white_ziegler}, and \cref{dual_translation}. For a more thorough explanation, see the claim in the proof of \cite[Lemma 4.1]{karp}.

\subsection{An orthogonally complementary view of the amplituhedron $\mathcal{A}_{n,k,m}(Z)$}\label{sec_B_amplituhedron}

\noindent The amplituhedron $\mathcal{A}_{n,k,m}(Z)$ is a subset of $\Gr_{k,k+m}$. Since we are considering the case $m=1$, it will be convenient for us to take orthogonal complements and work with subspaces of dimension $m$, rather than codimension $m$. To this end, we define an object $\mathcal{B}_{n,k,m}(W)\subseteq\Gr_m(W)$ for $W\in\Gr_{k+m,n}^{>0}$, which we show is homeomorphic to $\mathcal{A}_{n,k,m}(Z)$ (\cref{B_isomorphism}), where $W = \rowspan(Z)$. We remark that in the context of scattering amplitudes when $m=4$, $W$ is the subspace of $\mathbb{R}^n$ spanned by $4$ bosonic variables and $k$ fermionic variables.
\begin{defn}\label{defn_B_amplituhedron}
Given $W\in\Gr_{k+m,n}^{>0}$, let
\begin{align*}
\mathcal{B}_{n,k,m}(W) := \{V^\perp\cap W : V\in\Gr_{k,n}^{\ge 0}\}\subseteq\Gr_m(W),
\end{align*}
where $\Gr_m(W)$ denotes the subset of $\Gr_{m,n}$ of elements $X\in\Gr_{m,n}$ with $X\subseteq W$. Let us show that $\mathcal{B}_{n,k,m}(W)$ is well defined, i.e.\ $\dim(V^\perp\cap W) = m$ for all $V\in\Gr_{k,n}^{\ge 0}$.  By \cref{cor:disjoint} we have $V\cap W^\perp = \{0\}$, so the sum $V + W^\perp$ is direct. Hence $\dim(V+W^\perp) = \dim(V) + \dim(W^\perp) = n-m$. Since $(V^\perp\cap W)^\perp = V + W^\perp$, we get $\dim(V^\perp\cap W) = m$. (We remark that this is the same idea used in \cite[Section 4]{karp} to determine when, given an arbitrary linear map $Z:\mathbb{R}^n\to\mathbb{R}^{k+m}$, the image $\tilde{Z}(\Gr_{k,n}^{\ge 0})$ is well defined in $\Gr_{k,k+m}$.) 
\end{defn}

\begin{rmk}
While we were preparing this manuscript, we noticed that a 
similar construction appeared in Lam's definition of {\itshape universal amplituhedron varieties} \cite[Section 18]{lam}. There are two main differences between his construction and ours. First, Lam 
allows $Z$ to vary (hence the term ``universal''). 
Second, he 
works with complex varieties, and does not impose any positivity conditions on $V$ or $Z$ (rather, he restricts $V$ to lie in a closed complex positroid cell in $\Gr_{k,n}(\mathbb{C})$). 
Correspondingly he works with rational maps, 
while we will need our maps to be well defined everywhere.
\end{rmk}

We now show that $\mathcal{B}_{n,k,m}(W)$ is homeomorphic to $\mathcal{A}_{n,k,m}(Z)$, where $Z$ is any $(k+m)\times n$ matrix ($n \ge k+m$) with positive maximal minors and row span $W$. The idea is that we obtain $\mathcal{B}_{n,k,m}(W)$ from $\mathcal{A}_{n,k,m}(Z)\subseteq\Gr_{k,k+m}$ by taking orthogonal complements in $\mathbb{R}^{k+m}$, and then applying an isomorphism from $\mathbb{R}^{k+m}$ to $W$, so that our subspaces lie in $W$, not $\mathbb{R}^{k+m}$.
\begin{lem}\label{A_B_lemma}
Let $Z:\mathbb{R}^n\to\mathbb{R}^{k+m}$ be a surjective linear map, which we also regard as a $(k+m)\times n$ matrix, and let $W\in\Gr_{k+m,n}$ be the row span of $Z$. Then the map $f_Z : \Gr_m(W)\to\Gr_{k,k+m}$ given by
$$
f_Z(X) := Z(X^\perp) = \{Z(x) : x\in X^\perp\} \quad \text{ for all }X\in\Gr_m(W)
$$
is well defined and an isomorphism. (Here $X^\perp\in\Gr_{n-m,n}$ denotes the orthogonal complement of $X$ in $\mathbb{R}^n$; we use the notation $Z(X^\perp)$, and not $\tilde{Z}(X^\perp)$, because $\dim(X^\perp) \neq k$.) 

Moreover, for $X\in\Gr_m(W)$ with corresponding point $Y := f_Z(X)\in\Gr_{k,k+m}$, we can write the Pl\"{u}cker coordinates of $X$ (as an element of $\Gr_{m,n}$) in terms of $Y$ and $Z$, as follows. Let $z_1, \dots, z_n\in\mathbb{R}^{k+m}$ be the columns of the $(k+m)\times n$ matrix $Z$, and $y_1, \dots, y_k\in\mathbb{R}^{k+m}$ be a basis of $Y$. Then
for $1 \le i_1 < \cdots < i_m \le n$, we have
\begin{align}\label{plucker_translation}
\Delta_{\{i_1, \dots, i_m\}}(X) = \det([y_1 \,|\, \cdots \,|\, y_k \,|\, z_{i_1} \,|\, \cdots \,|\, z_{i_m}]).
\end{align}
\end{lem}

The formula \eqref{plucker_translation} is stated in \cite[Section 18]{lam}, 
though the proof is deferred to a forthcoming paper.
\begin{pf}
First let us show that $f_Z$ is well defined, i.e.\ for $X\in\Gr_m(W)$, we have $\dim(f_Z(X)) = k$. Since $X\subseteq W$, we have $W^\perp\subseteq X^\perp$, so we can write 
$X^\perp = V\oplus W^\perp$ 
for some $V\in\Gr_{k,n}$. 
Since $\ker(Z) = W^\perp$, we have
$Z(X^\perp) = Z(V+W^\perp) = Z(V)$, and then 
$V \cap \ker(Z) = V\cap W^\perp = \{0\}$ implies
$\dim(Z(V)) = k$.

To see that $f_Z$ is injective, suppose that we have $X,X'\in\Gr_m(W)$ with $Z(X^\perp) = Z(X'^\perp)$. Then $X^\perp + \ker(Z) = X'^\perp + \ker(Z)$. Taking orthogonal complements and using the fact that $\ker(Z) = W^\perp$, we obtain $X\cap W = X'\cap W$, so $X = X'$. Now we describe the inverse of $f_Z$. Given $Y\in\Gr_{k,k+m}$, consider 
the subspace
$Z^{-1}(Y) = \{v\in\mathbb{R}^n : Z(v)\in Y\}$.  Since 
$Z^{-1}(Y)$ contains 
$\ker(Z) = W^\perp$, which has dimension $n-k-m$,
and $\dim(Y) = k$, we have
$\dim(Z^{-1}(Y)) = n-m$.
Therefore we can write 
$Z^{-1}(Y)=X^\perp$ for some  $X\in \Gr_m(W)$, and then $Y = f_Z(X)$. 
It follows that $f_Z$ is invertible, and hence an isomorphism.

Now given $X\in\Gr_m(W)$ and $Y := f_Z(X)\in\Gr_{k,k+m}$, we prove \eqref{plucker_translation}. Fix column vectors $y_1, \dots, y_k\in\mathbb{R}^{k+m}$ which form a basis of $Y$ and $x_1, \dots, x_{n-m}\in\mathbb{R}^n$ which form a basis of $X^\perp$. For $i\in [n-m]$, we can write $Z(x_i) = \sum_{j=1}^kC_{i,j}y_j$ for some $C_{i,j}\in\mathbb{R}$. Let $D$ be the $(n-m)\times n$ matrix whose rows are $x_1^T, \dots, x_{n-m}^T$. Then the elements of $\Gr_{n-m,n+k}$ and $\Gr_{k+m,n+k}$ which are given by the row spans of 
$$
[-C \,|\, D] \qquad\text{ and }\qquad [y_1 \,|\, \cdots \,|\, y_k \,|\, Z],
$$
respectively, are orthogonally complementary. Given $I = \{i_1 < \cdots < i_m\}\subseteq [n]$ with complement $[n]\setminus I = \{j_1 < \cdots < j_{n-m}\}$, applying \cref{dual_translation}(ii) twice gives
\begin{align*}
\Delta_I(X) &= \Delta_{[n]\setminus I}(\alt(X^\perp)) = \Delta_{\{k+j_1, \dots, k+j_{n-m}\}}(\alt([-C \,|\, D])) = \\
&\qquad\quad\Delta_{\{1, \dots, k, k+i_1, \dots, k+i_m\}}([y_1 \,|\, \cdots \,|\, y_k \,|\, Z]) = \det([y_1 \,|\, \cdots \,|\, y_k \,|\, z_{i_1} \,|\, \cdots \,|\, z_{i_m}]).\qedhere
\end{align*}
\end{pf}

\begin{prop}\label{B_isomorphism}
Suppose that $Z$ is a $(k+m)\times n$ matrix ($n\ge k+m$) with positive maximal minors, and $W\in\Gr_{k+m,n}^{>0}$ is the row span of $Z$. Then the map $f_Z$ from \cref{A_B_lemma} restricts to a homeomorphism from $\mathcal{B}_{n,k,m}(W)$ onto $\mathcal{A}_{n,k,m}(Z)$, which sends $V^\perp\cap W$ to $\tilde{Z}(V)$ for all $V\in\Gr_{k,n}^{\ge 0}$. The Pl\"{u}cker coordinates of $V^\perp\cap W$ can be written in terms of $\tilde{Z}(V)$ and $Z$ by \eqref{plucker_translation}.
\end{prop}

\begin{eg}
Let $(n,k,m) := (4,2,1)$, and $Z:\mathbb{R}^4\to\mathbb{R}^3$ be given by the matrix
$$
\begin{bmatrix}
1 & 0 & 0 & 1\\
0 & 1 & 0 & -1 \\
0 & 0 & 1 & 1
\end{bmatrix} =: [z_1 \,|\, z_2 \,|\, z_3 \,|\, z_4],
$$
whose $3\times 3$ minors are all positive. Also let $V\in\Gr_{2,4}^{\ge 0}$ be the row span of the matrix $\begin{bmatrix}1 & a & 0 & 0 \\ 0 & 0 & 1 & b\end{bmatrix}$, where $a,b\ge 0$, and define $Y := \tilde{Z}(V)$. We can explicitly find a basis $y_1, y_2\in\mathbb{R}^3$ of $Y$ as follows:
$$
[y_1 \,|\, y_2] := \begin{bmatrix}
1 & 0 & 0 & 1\\
0 & 1 & 0 & -1 \\
0 & 0 & 1 & 1
\end{bmatrix}\begin{bmatrix}
1 & 0 \\
a & 0 \\
0 & 1 \\
0 & b
\end{bmatrix} = \begin{bmatrix}
1 & b \\
a & -b \\
0 & 1+b
\end{bmatrix}.
$$
Now let $X := V^\perp\cap\rowspan(Z)$, so that $X\in\mathcal{B}_{4,2,1}(\rowspan(Z))$ is mapped to $Y\in\mathcal{A}_{4,2,1}(Z)$  under the homeomorphism of \cref{B_isomorphism}. We can write $X$ as the line spanned by $(a(b+1), -(b+1), -b(a+1), a+1)$. We can check that we have
\begin{align*}
& \Delta_{\{1\}}(X) = \begin{vmatrix}
1 & b & 1 \\
a & -b & 0 \\
0 & 1+b & 0
\end{vmatrix}, && \Delta_{\{3\}}(X) = \begin{vmatrix}
1 & b & 0 \\
a & -b & 0 \\
0 & 1+b & 1
\end{vmatrix}, \\
& \Delta_{\{2\}}(X) = \begin{vmatrix}
1 & b & 0 \\
a & -b & 1 \\
0 & 1+b & 0
\end{vmatrix}, && \Delta_{\{4\}}(X) = \begin{vmatrix}
1 & b & 1 \\
a & -b & -1 \\
0 & 1+b & 1
\end{vmatrix},
\end{align*}
as asserted by \eqref{plucker_translation}. (Here $|M|$ denotes the determinant of $M$.)
\end{eg}

\begin{pf}[of \cref{B_isomorphism}]
Given $V\in\Gr_{k,n}^{\ge 0}$, since $\ker(Z) = W^\perp$ we have
$$
f_Z(V^\perp\cap W) = Z((V^\perp\cap W)^\perp) = Z(V+W^\perp) = Z(V) = \tilde{Z}(V).
$$
Thus the image $f_Z(\mathcal{B}_{n,k,m}(W))$ equals $\mathcal{A}_{n,k,m}(Z)$, and the result follows from \cref{A_B_lemma}.
\end{pf}

\subsection{A hypothetical intrinsic description of the amplituhedron}\label{sec_intrinsic_description}

We now give a description of the amplituhedron $\mathcal{B}_{n,k,1}(W)$ 
which does not mention $\Gr_{k,n}^{\ge 0}$.  This description will extend
to $\mathcal{B}_{n,k,m}(W)$ for $m>1$ if  
part (i) of the following problem has a positive answer.
\begin{prob}\label{lam_problem}
Let $V\in\Gr_{m,n}$, and $l\ge m$. ~\\
(i) If $\var(v)\le l-1$ for all $v\in V$, can we extend $V$ to an element of $\Gr_{l,n}^{\ge 0}$? \\
(ii) If $\overline{\var}(v)\le l-1$ for all $v\in V\setminus\{0\}$, can we extend $V$ to an element of $\Gr_{l,n}^{>0}$?
\end{prob}

\begin{rmk}
In March 2017, after reading a preliminary version of this article, Pavel Galashin \cite{galashin} explained to us an example (for $l=3$, $n=6$) showing that both parts of \cref{lam_problem} have a negative answer. (\cref{prob:equality} remains open.) Galashin's example leads to interesting questions in the study of {\itshape Grassmann polytopes}. See \cref{grassmann_polytope_remark} for further discussion, as well as some motivation for having posed \cref{lam_problem}.
\end{rmk}

\begin{lem}\label{intrinsic_description}
For $W\in\Gr_{k+m,n}^{>0}$, we have
$$
\mathcal{B}_{n,k,m}(W) \subseteq \{X\in\Gr_m(W) : k \le \overline{\var}(v) \le k+m-1 \text{ for all } v\in X\setminus\{0\}\}.
$$
If \cref{lam_problem}(i) has a positive answer for $l = n-k$, then equality holds.
\end{lem}
In calculating $\overline{{\var}}(v)$ for $v\in X$ (where $X\in\Gr_m(W)$), we regard $v$ as a vector in $\mathbb{R}^n$.
\begin{pf}
Given $X\in\mathcal{B}_{n,k,m}(W)$, we can write $X = V^\perp\cap W$ for some $V\in\Gr_{k,n}^{\ge 0}$. Then for any $v\in X\setminus\{0\}$, we have $\overline{\var}(v)\ge k$ by \cref{gantmakher_krein}(i) applied to $V$, and $\overline{\var}(v)\le k+m-1$ by \cref{gantmakher_krein}(ii) applied to $W$. This proves the containment. Now suppose that \cref{lam_problem}(i) has a positive answer for $l = n-k$. Let $X\in\Gr_m(W)$ with $\overline{\var}(v)\ge k$ for all $v\in X\setminus\{0\}$.  Then $\var(w)\le n-k-1$ for all $w\in\alt(X)\setminus\{0\}$ by \cref{dual_translation}(i). Hence we can extend $\alt(X)$ to an element of $\Gr_{n-k,n}^{\ge 0}$, which by \cref{dual_translation}(ii) we can write as $\alt(V^\perp)$ for some $V\in\Gr_{k,n}^{\ge 0}$. Since $\alt(X)\subseteq\alt(V^\perp)$, we have $X\subseteq V^\perp$, whence $X\subseteq V^\perp\cap W$. Since $\dim(X) = m = \dim(V^\perp\cap W)$, we have $X = V^\perp\cap W\in\mathcal{B}_{n,k,m}(W)$.
\end{pf}

\begin{prob}\label{prob:equality}
Do we have equality in 
\cref{intrinsic_description}?  In other words, is it true that 
for $W\in\Gr_{k+m,n}^{>0}$, we have
$$
\mathcal{B}_{n,k,m}(W)= \{X\in\Gr_m(W) : k \le \overline{\var}(v) \le k+m-1 \text{ for all } v\in X\setminus\{0\}\}? $$
\end{prob}

\begin{rmk}
We observe that \cref{prob:equality} has a positive 
answer in the extreme cases $k=0$ (when $\mathcal{B}_{n,k,m}(W) = \{W\}$), $m=0$ (when $\mathcal{B}_{n,k,m}(W) = \{\{0\}\})$, and $k+m=n$ (when $\mathcal{B}_{n,k,m}(W) = \{V^\perp : V\in\Gr_{k,n}^{\ge 0}\}$). Also, by \cite[Lemma 4.1]{karp}, both parts (i) and (ii) of \cref{lam_problem}
have a positive answer for $m=1$ (and all $l$ and $n$), and so 
\cref{prob:equality} has a positive answer for  $m=1$. 
\end{rmk}

This gives the following explicit description of $\mathcal{B}_{n,k,1}(W)$, 
which will be important for us in our study of the structure of $\mathcal{B}_{n,k,1}(W)$.
\begin{cor}\label{B_m=1}
For $W\in\Gr_{k+1,n}^{>0}$, we have
$$
\mathcal{B}_{n,k,1}(W) = \{w\in\mathbb{P}(W): \overline{\var}(w) = k\}\subseteq\mathbb{P}(W).
$$
\end{cor}

We now translate the hypothetical description of \cref{prob:equality} into one in terms of sign changes of certain sequences of Pl\"{u}cker coordinates. 
This is reminiscent of a 
description of the amplituhedron conjectured by 
Arkani-Hamed, Thomas, and Trnka \cite{ATT}. Namely, our \cref{plucker_description_translation} suggests that the interior of $\mathcal{A}_{n,k,m}(Z)$ is described by looking at certain sequences indexed by subsets $I\in\binom{[n]}{m-1}$. Arkani-Hamed, Thomas, and Trnka conjecture that for $m=4$, we need only consider subsets $I$ of the form $\{i, j, j+1\}$ \cite[(5.11)]{ATT}.\footnote{We note that the sequence in \cite{ATT} corresponding to $I = \{i , j, j+1\}$ differs from ours by a cyclic shift. However, in light of \cref{gale_remark}, requiring that our sequence has at least $k$ sign changes is equivalent to requiring that their sequence has exactly $k$ sign changes.} For arbitrary $m$, they conjecture that we need only consider the single subset $I = [m-1]$, if we also impose the constraint $\det([y_1 \,|\, \cdots \,|\, y_k \,|\, z_{i_1} \,|\, \cdots \,|\, z_{i_{m}}]) > 0$ for all $\{i_1, \dots, i_m\}$ satisfying {\itshape Gale's evenness condition} \cite[(5.14)]{ATT} (see \cref{gale_remark}). We thank Jara Trnka for telling us about their work ahead of its appearance on the arXiv. 

\begin{prop}\label{plucker_description}
Let $W\in\Gr_{k+m,n}^{>0}$, and define the open subset of $\Gr_m(W)$
$$
\mathcal{G} := \{X\in\Gr_m(W) : \var(v) \ge k \text{ for all } v\in X\setminus\{0\}\}.
$$
We have $\interior(\mathcal{B}_{n,k,m}(W))\subseteq\mathcal{G}$ and $\mathcal{B}_{n,k,m}(W)\subseteq\overline{\mathcal{G}}$, where $\interior(\cdot)$ and $\overline{\,\cdot\,}$ denote interior and closure. If equality holds in \cref{intrinsic_description}, then both containments above are equalities. Independently, we can describe $\mathcal{G}$  in terms of Pl\"{u}cker coordinates: 
$$
\mathcal{G} =  \left\{X\in\Gr_m(W) : \begin{aligned}
& \var(((-1)^{|I\cap [j]|}\Delta_{I\cup\{j\}}(X))_{j\in [n]\setminus I})\ge k\text{ for all } I\in\textstyle\binom{[n]}{m-1} \text{ such} \\
& \text{\,that the sequence $(\Delta_{I\cup\{j\}}(X))_{j\in [n]\setminus I}$ is not identically zero}
\end{aligned}\right\}.
$$
\end{prop}

We prove \cref{plucker_description} below. First we translate this description to $\mathcal{A}_{n,k,m}(Z)$ using \cref{A_B_lemma}, make a remark, and give an example.
\begin{cor}\label{plucker_description_translation}
Let $Z$ be a $(k+m)\times n$ matrix whose maximal minors are all positive, with row span $W\in\Gr_{k+m,n}^{>0}$. Let $\mathcal{F} := f_Z(\mathcal{G})$ be the image of the set $\mathcal{G}$ from \cref{plucker_description} under the map $f_Z$ from \cref{A_B_lemma}, and $z_1, \dots, z_n\in\mathbb{R}^{k+m}$ be the columns of $Z$. Then
\begin{align*}
\mathcal{F} = \left\{\langle y_1, \dots, y_k\rangle : \begin{aligned}
& \var((\det([y_1 \,|\, \cdots \,|\, y_k \,|\, z_{i_1} \,|\, \cdots \,|\, z_{i_{m-1}} \,|\, z_j]))_{j\in [n]\setminus I})\ge k \\
& \text{\,for all } I = \{i_1 < \cdots < i_{m-1}\}\in\textstyle\binom{[n]}{m-1} \text{ such that} \\
& \text{\,this sequence of minors is not identically zero}
\end{aligned}\right\}.
\end{align*}
We have $\interior(\mathcal{A}_{n,k,m}(Z))\subseteq\mathcal{F}$ and $\mathcal{A}_{n,k,m}(Z)\subseteq\overline{\mathcal{F}}$. If equality holds in \cref{intrinsic_description}, then both containments above are equalities.
\end{cor}
\begin{rmk}\label{gale_remark}
For all $X\in\Gr_m(W)$, we have $\overline{\var}(v)\le k+m-1$ for all $v\in X\setminus\{0\}$ by \cref{gantmakher_krein}(ii), whence $\overline{\var}((\Delta_{I\cup\{j\}}(X))_{j\in [n]\setminus I})\le k$ for all $I\in\binom{[n]}{m-1}$ such that the sequence $(\Delta_{I\cup\{i\}}(X))_{j\in [n]\setminus I}$ is not identically zero (see \cite[Theorem 3.1(ii)]{karp}). Note that the sequence $(\Delta_{I\cup\{j\}}(X))_{j\in [n]\setminus I}$ is the one in the description of $\mathcal{G}$, without the sign $(-1)^{|I\cap [j]|}$. For the sequence in the description of $\mathcal{F}$, this sign change corresponds to moving $z_j$ from the right end of the matrix into its proper relative position in the submatrix of $Z$. If $I\in\binom{[n]}{m-1}$ satisfies the condition that $(-1)^{|I\cap [j]|}$ is the same for all $j\in [n]\setminus I$ (called {\itshape Gale's evenness condition} \cite{gale}), then these two sequences are the same up to multiplication by $\pm 1$.
\end{rmk}

\begin{eg}\label{pentagon_example}
Let $(n,k,m) := (5,1,2)$, and $Z:\mathbb{R}^5\to\mathbb{R}^3$ be given by the matrix
$$
\begin{bmatrix}
1 & 0 & 0 & 1 & 3\\
0 & 1 & 0 & -1 & -2 \\
0 & 0 & 1 & 1 & 1
\end{bmatrix} =: [z_1 \,|\, z_2 \,|\, z_3 \,|\, z_4 \,|\, z_5],
$$
whose $3\times 3$ minors are all positive. In this case $\mathcal{A}_{5,1,2}(Z)$ is a convex pentagon in $\mathbb{P}^2 = \Gr_{1,3}$. Given $y\in\mathbb{R}^3$, \cref{plucker_description_translation} says that $\langle y\rangle\in\mathcal{F}$ if and only if each sequence $(\det([y \,|\, z_i \,|\, z_j]))_{j\in [5]\setminus \{i\}}$, for $i\in [5]$, is either identically zero or changes sign at least once. For example, if $i=3$, this sequence is
$$
\left(\begin{vmatrix}y_1 & 0 & 1 \\ y_2 & 0 & 0 \\ y_3 & 1 & 0\end{vmatrix},\begin{vmatrix}y_1 & 0 & 0 \\ y_2 & 0 & 1 \\ y_3 & 1 & 0\end{vmatrix},\begin{vmatrix}y_1 & 0 & 1 \\ y_2 & 0 & -1 \\ y_3 & 1 & 1\end{vmatrix},\begin{vmatrix}y_1 & 0 & 3 \\ y_2 & 0 & -2 \\ y_3 & 1 & 1\end{vmatrix}\right) = (y_2, -y_1, y_1+y_2, 2y_1+3y_2).
$$
Geometrically, the sequence corresponding to $i$ records where the point $\langle y\rangle\in\mathbb{P}^2$ lies in relation to each of the $4$ lines joining vertex $i$ to another vertex of the pentagon. If this sequence does not change sign, then $\langle y\rangle$ lies on the same side of all $4$ of these lines, i.e.\ the line segment between $\langle y\rangle$ and vertex $i$ does not intersect the interior of the pentagon. We see that $\mathcal{F}$ is the interior of the pentagon. In general, $\mathcal{F}$ is the interior of $\mathcal{A}_{n,k,m}(Z)$ if $k=1$, independently of \cref{prob:equality}.
\end{eg}

\begin{pf}[of \cref{plucker_description}]
First we prove $\interior(\mathcal{B}_{n,k,m}(W))\subseteq\mathcal{G}$, i.e.\ given $X\in\mathcal{B}_{n,k,m}(W)\setminus\mathcal{G}$, we have $X\notin\interior(\mathcal{B}_{n,k,m}(W))$. Since $X\notin\mathcal{G}$, there exists $v\in X\setminus\{0\}$ with $\var(v) < k$. Let $\sigma := \sign(v)$, and take $\tau\in\{+,-\}^n$ such that $\tau\ge\sigma$ and $\var(\tau) = \var(\sigma)$. Then by \cref{positive_covectors}, there exists $w\in W$ with $\sign(w) = \tau$. Now extend $v$ to a basis $v, v_2, \dots, v_m$ of $X$, and for $t > 0$ let $X_t := \spn(v + tw, v_2, \dots, v_m)$, so that $X_t\in\Gr_m(W)$ except for at most one value of $t$. Since $\overline{\var}(v+tw) = \overline{\var}(\tau) = \var(v) < k$, we have $X_t\notin\mathcal{B}_{n,k,m}(W)$ by \cref{intrinsic_description}. But $\{X_t\in\Gr_m(W) : t > 0\}$ intersects every neighborhood of $X$, so $X\notin\interior(\mathcal{B}_{n,k,m}(W))$.
The fact that $\mathcal{B}_{n,k,m}(W)\subseteq\overline{\mathcal{G}}$ follows from these two facts:
\begin{itemize}
\item $\Gr_{k,n}^{\ge 0} = \overline{\Gr_{k,n}^{>0}}$ \cite[Section 17]{postnikov};
\item if $V\in\Gr_{k,n}^{>0}$, then $V^\perp\cap W\in \mathcal{G}$ (by \cref{gantmakher_krein}(ii)).
\end{itemize}
Conversely, if equality holds in \cref{intrinsic_description}, then $\mathcal{G}\subseteq\mathcal{B}_{n,k,m}(W)$. Since $\mathcal{G}$ is open and $\mathcal{B}_{n,k,m}(W)$ is closed, we get the reverse containments $\mathcal{G}\subseteq\interior(\mathcal{B}_{n,k,m}(W))$ and $\overline{\mathcal{G}}\subseteq\mathcal{B}_{n,k,m}(W)$.

Now we describe $\mathcal{G}$ in terms of Pl\"{u}cker coordinates. By \cref{dual_translation}(i), we have
$$
\mathcal{G} =  \{X\in\Gr_m(W) : \overline{\var}(v) \le n-k-1 \text{ for all } v\in \alt(X)\setminus\{0\}\}.
$$
Theorem 3.1(ii) of \cite{karp} states that for $X'\in\Gr_{m,n}$, we have $\overline{\var}(v) \le n-k-1$ for all $v\in X'\setminus\{0\}$ if and only if $\overline{\var}((\Delta_{I\cup\{j\}}(X'))_{j\in [n]\setminus I})\le n-k-m$ for all $I\in\binom{[n]}{m-1}$ such that the sequence $(\Delta_{I\cup\{j\}}(X'))_{j\in [n]\setminus I}$ is not identically zero. Also by \cref{dual_translation}(i), we have $\overline{\var}(p) + \var(\alt(p)) = n-m$ for all nonzero $p\in\mathbb{R}^{[n]\setminus I}$ ($I\in\binom{[n]}{m-1}$), where $\alt$ acts on $\mathbb{R}^{[n]\setminus I}$ by changing the sign of every second component. We get
$$
\mathcal{G} =  \left\{X\in\Gr_m(W) : \begin{aligned}
& \var(\alt(\Delta_{I\cup\{j\}}(\alt(X)))_{j\in [n]\setminus I})\ge k\text{ for all } I\in\textstyle\binom{[n]}{m-1} \text{ such} \\
& \text{\,that the sequence $(\Delta_{I\cup\{j\}}(X))_{j\in [n]\setminus I}$ is not identically zero}
\end{aligned}\right\}.
$$
In order to obtain the desired description of $\mathcal{G}$, it remains to show that given $I\in\binom{[n]}{m-1}$ such that the sequence $(\Delta_{I\cup\{j\}}(X))_{j\in [n]\setminus I}$ is not identically zero, we have
\begin{align}\label{var_equality}
\var(\alt(\Delta_{I\cup\{j\}}(\alt(X)))_{j\in [n]\setminus I}) = \var(((-1)^{|I\cap [j]|}\Delta_{I\cup\{j\}}(X))_{j\in [n]\setminus I}).
\end{align}
To this end, write $I = \{i_1, \dots, i_{m-1}\}\in\binom{[n]}{m-1}$. Then for $j\in [n]\setminus I$, component $j$ of $\alt(\Delta_{I\cup\{j'\}}(\alt(X)))_{j'\in [n]\setminus I}\in\mathbb{R}^{[n]\setminus I}$ equals
$$
(-1)^{|([n]\setminus I)\cap [j]|-1}(-1)^{(i_1 - 1) + \cdots + (i_{m-1} - 1) + (j-1)}\Delta_{I\cup\{j\}}(X) = \epsilon_I(-1)^{|I\cap [j]|}\Delta_{I\cup\{j\}}(X),
$$
where $\epsilon_I := (-1)^{(i_1 - 1) + \cdots + (i_{m-1} - 1)} = \pm 1$ does not depend on $j$. This gives \eqref{var_equality}.
\end{pf}

\begin{pf}[of \cref{plucker_description_translation}]
Applying \eqref{plucker_translation} to the description of $\mathcal{G}$ in \cref{plucker_description} gives
$$
\mathcal{F} = \left\{\langle y_1, \dots, y_k\rangle : \begin{aligned}
& \var(((-1)^{|I\cap [j]|}\det([y_1 \,|\, \cdots \,|\, y_k \,|\, Z_{I\cup\{j\}}]))_{j\in [n]\setminus I})\ge k \\
& \text{\,for all } I = \{i_1 < \cdots < i_{m-1}\}\in\textstyle\binom{[n]}{m-1} \text{ such that} \\
& \text{\,this sequence of minors is not identically zero}
\end{aligned}\right\},
$$
where $Z_J$ denotes the submatrix of $Z$ with columns $J$, for $J\subseteq [n]$. Moving the column of $Z_{I\cup\{j\}}$ labeled by $j$ to the right end of 
the matrix introduces a sign $(-1)^{|I\cap (j,n]|}$ in the determinant. Since $(-1)^{|I\cap [j]|}(-1)^{|I\cap (j,n]|} = (-1)^{|I|}$, which does not depend on $j$, we obtain the stated description of $\mathcal{F}$. The rest follows from \cref{plucker_description}, using \cref{A_B_lemma}.
\end{pf}

\subsection{Removing $k$ from the definition of the amplituhedron}\label{sec_removing_k}

In the definition of $\mathcal{B}_{n,k,1}(W)$ (\cref{defn_B_amplituhedron}), if we let $W$ and $k$ vary, then we obtain the following object:
\begin{align}\label{master_amplituhedron_definition}
\widehat{\mathcal{B}}_{m,n} := \bigcup_{0\le k \le n-m}\{V^\perp\cap W : V\in\Gr_{k,n}^{\ge 0}, W\in\Gr_{k+m,n}^{>0}\}\subseteq\Gr_{m,n}.
\end{align}
By \cref{gantmakher_krein}, we have
\begin{align}\label{conjectured_master_amplituhedron}
\widehat{\mathcal{B}}_{m,n}\subseteq\bigcup_{0\le k \le n-m}\{X\in\Gr_{m,n} : k \le \overline{\var}(v) \le k+m-1 \text{ for all } v\in X\setminus\{0\}\},
\end{align}
and it is an interesting problem to determine if equality holds. If $m=1$, then \cref{lam_problem} has a positive answer, and an argument similar to the proof of \cref{intrinsic_description} implies $\widehat{\mathcal{B}}_{1,n} = \Gr_{1,n}$. 
However, $\widehat{\mathcal{B}}_{m,n}\neq\Gr_{m,n}$ for $m\ge 2$.

Our motivation for considering $\widehat{\mathcal{B}}_{m,n}$ is that in the BCFW recursion \cite{BCFW}, which conjecturally provides a ``triangulation'' of $\mathcal{A}_{n,k,4}(Z)$, the parameter $k$ is allowed to vary. This suggests that the diagrams appearing in the BCFW recursion corresponding to a given $n$ might
 label pieces of an object which somehow encompasses $\mathcal{A}_{n,k,4}(Z)$ for all  $k$. 
Could $\widehat{\mathcal{B}}_{4,n}$ be such an object?

As a first result about $\widehat{\mathcal{B}}_{m,n}$, we prove that the union in \eqref{master_amplituhedron_definition} is disjoint if $m\ge 1$. It suffices to show that the union in \eqref{conjectured_master_amplituhedron} is disjoint, which follows from the 
lemma below. It also follows from the lemma that the union in \eqref{conjectured_master_amplituhedron} consists of the elements of $\Gr_{m,n}$ whose range of $\overline{\var}(\cdot)$ (over nonzero vectors) is contained in as small an interval as possible. That is, for $X\in\Gr_{m,n}$ with $m\ge 1$, we have $\sup_{v,w\in X\setminus\{0\}}\overline{\var}(v) - \overline{\var}(w)\ge m-1$, and equality holds if and only if $X$ is contained in the union in \eqref{conjectured_master_amplituhedron}.
\begin{lem}
Suppose that $X\in\Gr_{m,n}$ with $m\ge 1$, and $k\ge 0$ such that $k \le \overline{\var}(v) \le k+m-1$ for all $v\in X\setminus\{0\}$. Then $\overline{\var}(v) = k+m-1$ for some $v\in X\setminus\{0\}$ and $\overline{\var}(w) = k$ for some $w\in X\setminus\{0\}$.
\end{lem}

\begin{pf}
Let $v^{(1)}, \dots, v^{(m)}$ be the rows of an $m\times n$ matrix whose rows span $X$, after it has been put into reduced row echelon form. That is, if $i_1 < \cdots < i_m$ index the pivot columns of this matrix, then we have $v^{(r)}_{i_s} = \delta_{r,s}$ for all $r,s\in [m]$, and $v^{(r)}_j = 0$ for all $r\in [m]$ and $j < i_r$. Let $v := v^{(m)}$, and note that $\overline{\var}(v) = \overline{\var}(v|_{[i_m,n]}) + i_m - 1$. Now let
$$
w := v + t\sum_{r=1}^{m-1} \epsilon_r v^{(r)},
$$
where $t > 0$ is sufficiently small that in positions $j\in [i_{m},n]$ where $v$ is nonzero, $w_j$ has the same sign as $v_j$, and $\epsilon_r\in\{1,-1\}$ is chosen to be $1$
precisely if $|\{i \in \mathbb{Z} : i_r < i < i_{r+1}\}|$ is even. 
Then $$k\le\overline{\var}(w)\le\overline{\var}(v|_{[i_m,n]})+i_m-m = 
\overline{\var}(v) -m+1 \leq (k+m-1)-m+1 = k,$$
and hence equality holds everywhere above.
\end{pf}

\section{A BCFW-like recursion for \texorpdfstring{$m=1$}{m=1}}\label{sec_bcfw}

\noindent In the case that $m=4$, the \emph{BCFW recursion} (named after
Britto, Cachazo, Feng, and Witten \cite{BCF, BCFW}) can be viewed as a 
procedure which outputs a subset of cells
of $\Gr_{k,n}^{\geq 0}$ whose images conjecturally ``triangulate'' the 
amplituhedron $\mathcal{A}_{n,k,4}(Z)$ 
\cite{arkani-hamed_trnka}.
The procedure is described in \cite[Section 2]{abcgpt} as an operation on plabic graphs.

In this section we give an $m=1$ analogue of the BCFW recursion,
which naturally leads us to a subset of cells 
of $\Gr_{k,n}^{\geq 0}$ that we call 
the \emph{$m=1$ BCFW cells of $\Gr_{k,n}^{\ge 0}$}. We remark that unlike the recursion for $m=4$ as described in \cite{abcgpt}, there is no `shift' to the decorated permutation involved.
These cells can be easily described in terms
of their \Le -diagrams, positroids, or decorated permutations.
As we will show in \cref{sec_induced_subcomplex}, the images
of these cells ``triangulate'' the $m=1$ amplituhedron
$\mathcal{A}_{n,k,1}(Z)$.  More specifically, they each map injectively into $\mathcal{A}_{n,k,1}(Z)$, and their images are disjoint
and together form a dense subset of $\mathcal{A}_{n,k,1}(Z)$. In fact, 
we will show that the BCFW cells in $\Gr_{k,n}^{\ge 0}$ 
plus the cells in their boundaries 
give rise to a cell decomposition of $\mathcal{A}_{n,k,1}(Z)$.
\begin{figure}[ht]
\includegraphics[height=.8in]{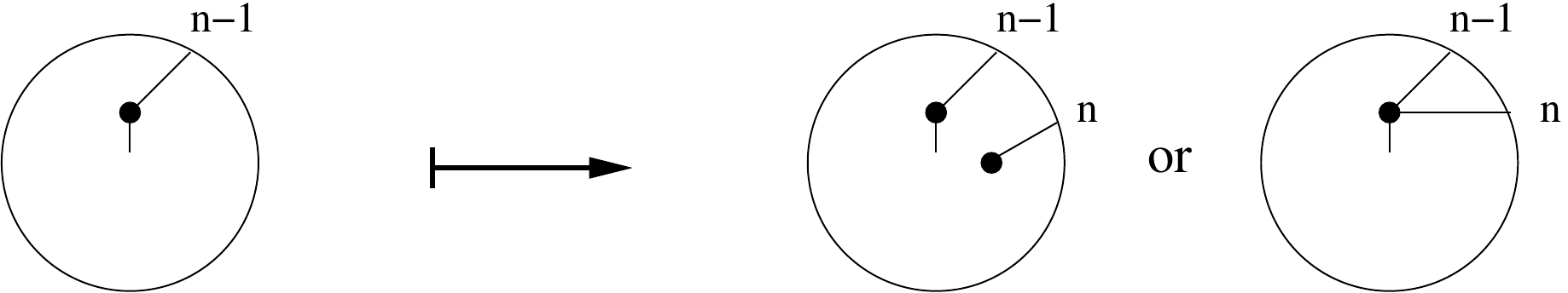}
\caption{The $m=1$ BCFW recursion: add a new 
vertex $n$, which is incident 
either to a black lollipop
or to the edge adjacent to 
vertex $n-1$.}
\label{BCFW}
\end{figure}

\begin{figure}[ht]
\begin{center}
\includegraphics[width=6.4in]{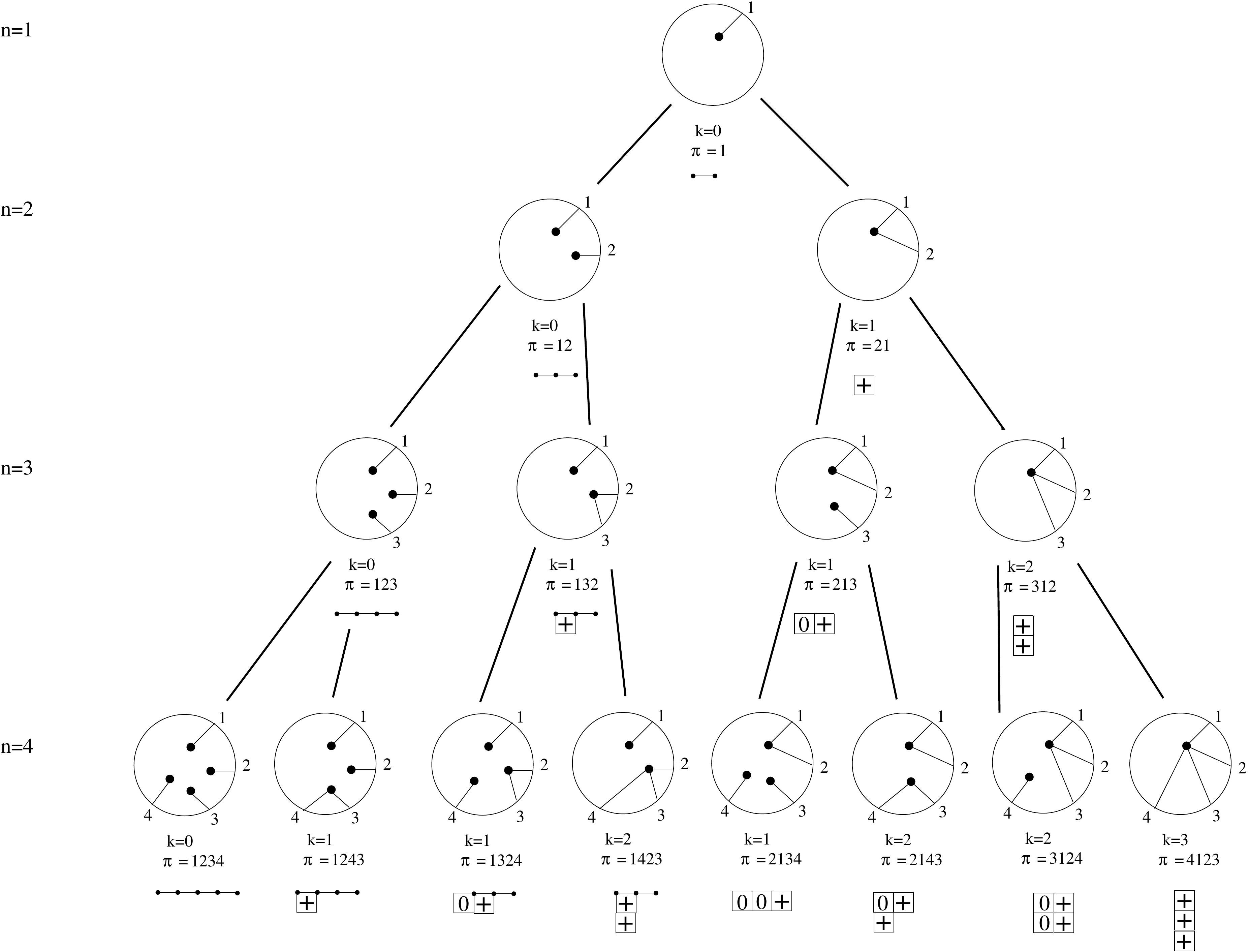}
\end{center}
\caption{The $m=1$ BCFW-style recursion for $n=1,2,3,4$.}
\label{BCFW2}
\end{figure}

\begin{defn}
The \emph{BCFW recursion for $m=1$} is defined as follows.  
\begin{itemize}
\item We start from the 
plabic graph with one boundary vertex, incident to a black lollipop.
This is our graph for $n=1$ (with $k=0$).  
\item Given a plabic graph with $n-1$ boundary vertices 
produced by our recursion (where $n\ge 2$), we perform one of 
the following two operations:
we add a new boundary vertex $n$ which is incident either to a black 
lollipop, or to the edge adjacent to boundary vertex $n-1$.
See \cref{BCFW}.  The first operation preserves the $k$
statistic, while the second operation increases it by $1$.
\end{itemize}
We refer to the set of all plabic graphs with fixed $n$ and $k$
statistics produced in this way as the 
\emph{$m=1$ BCFW cells of $\Gr_{k,n}^{\ge 0}$}. See \cref{BCFW2} for all $m=1$ BCFW cells of $\Gr_{k,n}^{\ge 0}$ with $n \leq 4$.
\end{defn}

The following lemma is easy to verify by inspection, using the bijections between
plabic graphs, \Le-diagrams, and decorated permutations which we gave in 
\cref{sec_TNN_Grassmannian}.  

\begin{lem}\label{BCFW_diagrams_permutations}
The $m=1$ BCFW cells of $\Gr_{k,n}^{\ge 0}$ are indexed by the 
\Le-diagrams of type $(k,n)$ such that each of the $k$ rows contains a unique $+$, which 
is at the far right of the row. The decorated permutation of such a \Le-diagram $D$ can be written in cycle notation as
$$
\pi(D) = (i_1, i_1-1, i_1-2, \dots , 1) (i_2 , i_2 - 1 ,\dots, i_1+1)
\cdots (n ,n-1, \dots , i_{n-k-1}+1),$$
where $i_1 < i_2 < \dots < i_{n-k-1} < i_{n-k} = n$ label the horizontal steps of the southeast border of $D$ (read northeast to southwest), and all fixed points are colored black. In particular, the number of $m=1$ BCFW cells of $\Gr_{k,n}^{\ge 0}$ equals $\binom{n-1}{k}$.
\end{lem}

\section{\texorpdfstring{$\mathcal{A}_{n,k,1}$}{A(n,k,1)} as a subcomplex of the totally nonnegative Grassmannian}\label{sec_induced_subcomplex}

\noindent In this section we show that the amplituhedron $\mathcal{B}_{n,k,1}(W)$ (for $W\in\Gr_{k+1,n}^{>0}$) is isomorphic to a subcomplex of $\Gr_{k,n}^{\ge 0}$. We begin by defining a stratification of $\mathcal{B}_{n,k,1}(W)$,
whose strata are indexed by sign vectors $\overline{\PSign_{n,k,1}}$ and have a natural poset structure.  We will also define a poset of certain
\Le-diagrams $\overline{\mathcal{D}_{n,k,1}}$ and a poset of 
positroids $\overline{\mathcal{M}_{n,k,1}}$ and show that 
all three posets are isomorphic.  Finally we will give an isomorphism 
between 
the amplituhedron $\mathcal{B}_{n,k,1}(W)$ and 
the subcomplex of $\Gr_{k,n}^{\ge 0}$ indexed by 
the cells associated to $\overline{\mathcal{D}_{n,k,1}}$, which 
induces an isomorphism on the posets of closures of strata, ordered
by containment.

Recall that by \cref{B_m=1}, we have
$$
\mathcal{B}_{n,k,1}(W) = \{w\in W\setminus\{0\}: \overline{\var}(w) = k\}\subseteq\mathbb{P}(W),
$$
where we identify a nonzero vector in $W$ with the line it spans in $\mathbb{P}(W)$. We define a stratification of $\mathcal{B}_{n,k,1}(W)$ using 
sign vectors.

\begin{defn}\label{defn_Sign}
Let $\overline{\Sign_{n,k,1}}$ denote the set of nonzero sign vectors $\sigma\in\{0,+,-\}^n$ with $\overline{\var}(\sigma) = k$, such that if $i\in [n]$ indexes the first nonzero component of $\sigma$, then $\sigma_i = (-1)^{i-1}$. (Equivalently, the first nonzero component of $\alt(\sigma)$ equals $+$.) 
Also let $\overline{\PSign_{n,k,1}}$ denote the set of nonzero sign vectors $\sigma\in\{0,+,-\}^n$ with $\overline{\var}(\sigma) = k$, modulo multiplication by $\pm 1$. 
We let $\Sign_{n,k,1}$ 
and ${\PSign_{n,k,1}}$ be the subsets of 
$\overline{\Sign_{n,k,1}}$ and
$\overline{\PSign_{n,k,1}}$ consisting of vectors with no zero components.
\end{defn}

\begin{defn}\label{defn_B_stratification}
We stratify the amplituhedron $\mathcal{B}_{n,k,1}(W)$ by $\overline{\PSign_{n,k,1}}$,
i.e.\ its strata are $\mathcal{B}_{\sigma}(W):=\{w\in W\setminus\{0\} : \sign(w) = \pm \sigma\}$ for $\sigma\in \overline{\PSign_{n,k,1}}$. All strata are nonempty by \cref{positive_covectors}(i). The strata are partially ordered 
by containment of closures, i.e.\ $\mathcal{B}_\sigma(W)\le\mathcal{B}_\tau(W)$ if and only if $\mathcal{B}_\sigma(W)\subseteq\overline{\mathcal{B}_\tau(W)}$.
\end{defn}
This stratification for $\mathcal{B}_{5,2,1}(W)$ is shown in \cref{hyperplane_arrangement_sign_vectors_5_2_1}. We will show in \cref{triangulation_m=1} and \cref{cyclic_arrangement} that the partial order on strata corresponds to 
a very natural partial order on $\overline{\PSign_{n,k,1}}$, which we now 
describe.

\begin{defn}\label{defn_sign_vector_partial_order}
We define a partial order on the set of sign vectors 
$\{0,+,-\}^n$ as follows:
$\sigma\le\tau$ if and only if $\sigma_i = \tau_i$ for all $i\in [n]$ such that $\sigma_i\neq 0$. Equivalently, $\sigma\le\tau$ if and only if we can obtain $\sigma$ by setting some components of $\tau$ to $0$.
This gives a partial order on 
$\overline{\Sign_{n,k,1}}$  by restriction.
And for nonzero sign vectors $\sigma,\tau$ representing elements in
$\overline{\PSign_{n,k,1}}$, 
we say that $\sigma\le\tau$  if and only if $\sigma\le\tau$ or $\sigma\le -\tau$ in $\overline{\Sign_{n,k,1}}$.
\end{defn}

For example, $(+,0,+,0,+) \le (+,-,+,+,+)$, but $(+,0,+,0,+)\nleq (+,-,0,+,+)$. 
\cref{hyperplane_arrangement_sign_vectors_5_2_1} shows
 $\overline{\Sign_{5,2,1}}$ 
as labels of the bounded faces of a hyperplane arrangement.

We will now show that $\overline{\Sign_{n,k,1}}$ and $\overline{\PSign_{n,k,1}}$ are isomorphic as posets. Our reason for using both posets
 is as follows. Since $\mathcal{B}_{n,k,1}(W)$ is a subset of the 
projective space $\mathbb{P}(W)$, the sign vectors used to index strata 
should be considered modulo multiplication by $\pm 1$, which leads us 
naturally to $\overline{\PSign_{n,k,1}}$. However, in \cref{cyclic_arrangement} we will show that $\mathcal{B}_{n,k,1}(W)$ is isomorphic to the bounded complex of a hyperplane arrangement, and the bounded faces of this arrangement are labeled by sign vectors {\itshape not} considered modulo multiplication by $\pm 1$. 
To prove this result, we will need to work with $\overline{\Sign_{n,k,1}}$, which requires a more careful analysis.
\begin{lem}\label{Sign_PSign}~\\
(i) The map $\overline{\Sign_{n,k,1}}\to\overline{\PSign_{n,k,1}}, \sigma\mapsto\sigma$ is an isomorphism of posets. \\
(ii) Conversely, suppose that $P$ is an induced subposet of $\{\sigma\in\{0,+,-\}^n\setminus\{0\} : \overline{\var}(\sigma) = k\}$ such that the map $P\to\overline{\PSign_{n,k,1}}, \sigma\mapsto\sigma$ is an isomorphism of posets. 
Then $P$ equals $\overline{\Sign_{n,k,1}}$ or $-\overline{\Sign_{n,k,1}}$.
\end{lem}
This says that $\pm\overline{\Sign_{n,k,1}}$ are the unique liftings of $\overline{\PSign_{n,k,1}}$ to $\{0,+,-\}^n$ which preserve its poset structure. For example, if $n=2, k=1$, then the lifting $P := \{(+,-), (+,0), (0,+)\}$ does not have the same poset structure as $\overline{\Sign_{2,1,1}} = \{(+,-), (+,0), (0,-)\}$, since $P$ has two maximal elements, but $\overline{\Sign_{2,1,1}}$ has the unique maximum $(+,-)$.
\begin{pf}
(i) The map $\overline{\Sign_{n,k,1}}\to\overline{\PSign_{n,k,1}}, \sigma\mapsto\sigma$ is a bijection and a poset homomorphism. To show that it is a poset isomorphism, we must show that there do not exist $\sigma,\tau\in\overline{\Sign_{n,k,1}}$ with $\sigma\le -\tau$. Suppose otherwise that such $\sigma$ and $\tau$ exist. Let $\sigma' := \alt(\sigma)$ and $\tau' := \alt(\tau)$. By \cref{dual_translation}(i), we have $\var(\sigma') = \var(\tau') = n-k-1$. Let $i,j\in [n]$ index the first nonzero components of $\sigma',\tau'$, respectively, so that $\sigma'_i = \tau'_j = +$. But also by our assumption, $\sigma'_i = -\tau'_i$.  Since  $j < i$, $-\tau'$ changes sign at least once from $j$ to $i$. This implies $\var(-\tau') > \var(\sigma')$, a contradiction.

(ii) Let $G$ be the graph with vertex set $\PSign_{n,k,1}$, where distinct $\sigma,\tau\in\PSign_{n,k,1}$ are adjacent if and only if $\sigma$ differs in a single component from either $\tau$ or $-\tau$.
\begin{claim}
$G$ is connected.
\end{claim}
\begin{claimpf}
Let us show that every vertex $\sigma$ of $G$ is connected to
$$
\widehat{\sigma} := (+,-,+,-,\dots, (-1)^{k-1}, (-1)^k, (-1)^k, \dots),
$$
i.e.\ $\widehat{\sigma}$ is the sign vector which alternates in sign on $[k+1]$, and is constant thereafter. Take $i\in [n]$ maximum
such that $\sigma$ alternates in sign on $[i]$. If $i = k+1$, then $\sigma = \widehat{\sigma}$. Otherwise, take $j > i$ minimum with $\sigma_j\neq\sigma_i$. That is, on the interval $[i,j]$, $\sigma$ equals $(+,+,\ldots,+,-)$ up to sign. By performing sign flips at components $j-1, j-2, \dots, i+1$, we obtain a sign vector which equals $(+,-,-,\dots,-)$ on $[i,j]$, and hence alternates in sign on $[i+1]$. Repeating this procedure for $i+ 1, i+2, \dots, k$, we obtain $\widehat{\sigma}$. For example, if $\sigma = (+,+,+,-,-,+,-)$ (where $n=7$, $k=3$), then we obtain $\widehat{\sigma}$ as follows:
$$
\sigma\xmapsto{i=1} (+,-,-,-,-,+,-) \xmapsto{i=2} (+,-,+,+,+,+,-) \xmapsto{i=3} (+,-,+,-,-,-,-) = \widehat{\sigma}.
$$
This proves the claim.
\end{claimpf}

Let $H$ be the corresponding graph for $P$, i.e.\ $H$ has vertex set $Q := P\cap\{+,-\}^n$, and two distinct sign vectors are adjacent in $H$ if and only if they differ in a single component. Note that $G$ (respectively $H$) depends only on the poset $\overline{\PSign_{n,k,1}}$ (respectively $P$): two distinct sign vectors are adjacent if and only if they cover a common sign vector in the poset. Since $P\cong\overline{\PSign_{n,k,1}}$ and $G$ is connected by sign flips, we get that $H$ is connected by sign flips. Note that we can never flip the first component, which would change sign variation. Hence all elements of $Q$ have the same first component. After replacing $P$ with $-P$ if necessary, we may assume that $\sigma_1 = +$ for all $\sigma\in Q$.

Let $\sigma\in P$. We will show that $\sigma_i = (-1)^{i-1}$, where $i$ indexes the first nonzero component of $\sigma$. Since $\overline{\var}(\sigma) = k$, we can extend $\sigma$ to $\tau\in\{+,-\}^n$ (i.e.\ $\tau\ge\sigma)$ with $\var(\tau) = k$. Then $\tau$ alternates in sign on $[i]$, i.e.\ it equals $((-1)^{i-1}\sigma_i, (-1)^{i-2}\sigma_i, \dots, \sigma_i)$ on $[i]$. Now 
since $P\to\overline{\PSign_{n,k,1}}, \sigma\mapsto\sigma$ is a poset isomorphism, 
exactly one of $\tau, -\tau$ is in $P$. If $-\tau\in P$, then $\sigma \nleq -\tau$ in $P$ but $\sigma\leq -\tau$ in $\overline{\PSign_{n,k,1}}$, a contradiction. Hence $\tau$ is in $P$ (and hence in $Q$), whence $\tau_1 = +$, i.e.\ $\sigma_i = (-1)^{i-1}$.
\end{pf}

We now define a subcomplex of $\Gr_{k,n}^{\ge 0}$ which will turn out to be isomorphic to $\mathcal{B}_{n,k,1}(W)$.

\begin{defn}\label{Le_m=1}
Let $\mathcal{D}_{n,k,1}$ 
(respectively,  $\overline{\mathcal{D}_{n,k,1}}$)
be the set of \Le-diagrams
contained in a $k \times (n-k)$ rectangle whose rows each have 
precisely one $+$ (respectively, at most one $+$), 
and each $+$ appears at the right end of its row.
For 
$D\in \overline{\mathcal{D}_{n,k,1}}$, we let
$\dim(D) := \dim(S_D)$ be the number of $+$'s in $D$.
\end{defn}
Note that $\mathcal{D}_{n,k,1}$ indexes the $m=1$ BCFW cells of $\Gr_{k,n}^{\ge 0}$ by \cref{BCFW_diagrams_permutations}.

Since \Le -diagrams index the cells of $\Gr_{k,n}^{\ge 0}$, 
$\overline{\mathcal{D}_{n,k,1}}$ 
has a poset structure as a subposet of $Q_{k,n}$. 
However, it is more convenient for us to define our own partial order on $\overline{\mathcal{D}_{n,k,1}}$, as follows; we then show in \cref{le_positroid_isomorphism} that our partial order agrees with the one coming from $Q_{k,n}$, and 
that our poset on 
$\overline{\mathcal{D}_{n,k,1}}$ 
is in fact an order ideal (a downset)  of $Q_{k,n}$.
\begin{defn}
We define a partial order on $\overline{\mathcal{D}_{n,k,1}}$, with the following cover relations $\lessdot$:
\begin{itemize}
\item (Type 1) Let $D\in \overline{\mathcal{D}_{n,k,1}}$, 
where $\dim(D) \geq 1$, and  choose some $+$ in $D$ which has no $+$'s below it
in the same column.  We  obtain
$D'$ from $D$ by deleting the box containing that $+$ and every box below it 
and in the same column.  
Then $D'\lessdot D$.
\item (Type 2) 
Let $D\in \overline{\mathcal{D}_{n,k,1}}$, 
where $\dim(D) \geq 1$, and  choose some $+$ in $D$.  We  obtain
$D'$ from $D$ by replacing the $+$ with a $0$.
Then $D'\lessdot D$.
\end{itemize}
\end{defn}

We now give a poset isomorphism 
$\DS: \overline{\mathcal{D}_{n,k,1}} \to \overline{\Sign_{n,k,1}}$.
\begin{defn}\label{def:Lesign}
Let $D\in {\mathcal{D}_{n,k,1}}$ with Young diagram $Y_{\lambda}$, 
where the steps of the southeast border of $Y_{\lambda}$ are
labeled from $1$
to $n$.
Then we define $\sigma(D) \in {\Sign_{n,k,1}}$ recursively by 
setting:
\begin{itemize}
\item $\sigma_1 := +$;
\item $\sigma_{i+1} = \sigma_i$ if and only if 
$i$ is the label of a horizontal step of $Y_{\lambda}$, $i = 1, \dots, n-1$.
\end{itemize}
Now let $D\in \overline{\mathcal{D}_{n,k,1}}$ with Young diagram $Y_{\lambda}$.
We obtain $\widehat{D}$ from $D$ by putting a $+$ at the far right of every row which has no $+$.
To define $\sigma(D)$, we first compute $\sigma(\widehat{D})$ as above, but then for every vertical 
step $i$ corresponding to an all-zero row of $D$, we set $\sigma(D)_i = 0$.
This gives a map $\DS: \overline{\mathcal{D}_{n,k,1}} \to \overline{\Sign_{n,k,1}}$ defined by $\DS(D) = \sigma(D)$.
\end{defn}
For examples of this bijection, compare \cref{hyperplane_arrangement_le_diagrams_5_2_1} with \cref{hyperplane_arrangement_sign_vectors_5_2_1}.
\begin{lem}\label{DS_bijection}
The map $\DS$ from \cref{def:Lesign} is an isomorphism of posets.
\end{lem}
\begin{pf}
First we show that $\DS$ is well defined, i.e.\ given $D\in\overline{\mathcal{D}_{n,k,1}}$ we have $\sigma(D)_i = (-1)^{i-1}$, where $i\in [n]$ indexes the first nonzero component of $\sigma(D)$. Indeed, we have that $1, \dots, i-1$ label vertical steps of the southeast border of $D$ whose rows contain only $0$'s, and $i$ labels either a vertical step whose row contains a $+$, or a horizontal step. Hence when we lift $D$ to $\widehat{D}$ in \cref{def:Lesign}, $\sigma(\widehat{D})$ alternates in sign on $[1,i]$, whence $\sigma(D)_i = \sigma(\widehat{D})_i = (-1)^{i-1}$.

To show that $\DS$ is a bijection, we describe its inverse. Given $\sigma \in \overline{\Sign_{n,k,1}}$, we construct 
$D\in \overline{\mathcal{D}_{n,k,1}}$ with $\sigma(D) = \sigma$. We first lift $\sigma$ to $\widehat{\sigma}\in\Sign_{n,k,1}$ by reading $\sigma$ from right to left:
\begin{itemize}
\item Set $\widehat{\sigma}_n := (-1)^k$.
\item For $i = n-1, \dots 1$, if $\sigma_i = 0$, we set $\widehat{\sigma}_i := -\widehat{\sigma}_{i+1}$, and otherwise we set $\widehat{\sigma}_i := \sigma_{i}$.
\end{itemize}
The reason we set $\widehat{\sigma}_n := (-1)^k$ is that if $\tau\in\{+,-\}^n$ with $\tau\ge\sigma$ and $\var(\tau) = k$, then $\tau_1 = +$ since $\tau\in\Sign_{n,k,1}$, which implies that $\tau_n = (-1)^k$. The same reasoning implies that $\widehat{\sigma}\in\Sign_{n,k,1}$. Now we construct $D$ as follows. The vertical steps of the southeast border of $D$ are in bijection with 
positions $i$ of $\widehat{\sigma}$ such that $\widehat{\sigma}_i = - \widehat{\sigma}_{i+1}$.
In each row with vertical step labeled $i$, where $\sigma_i\neq 0$, we place a $+$ in the far 
right.  We then fill the remaining boxes with $0$'s. Note that $\sigma(D) = \sigma$. Moreover, we know that $0$'s of $\sigma$ must correspond to vertical steps in the southeast border of $D$, which uniquely determines $D$. Therefore $\DS$ is a bijection.

We can check that $\DS$ is a poset homomorphism. It remains to check that its inverse is a poset homomorphism. To this end, let $\sigma'\lessdot\sigma$ be a cover relation in $\overline{\Sign_{n,k,1}}$, so that $\sigma'$ is obtained from $\sigma$ by setting some $\sigma_i$  to $0$. Construct $D\in\overline{\mathcal{D}_{n,k,1}}$ as above with $\sigma(D) = \sigma$. Since $\sigma_i\neq 0$, either $i$ labels 
a horizontal step in the southeast border of $D$ and $i-1$ labels a vertical step (otherwise we would have $\overline{\var}(\sigma') > \overline{\var}(\sigma)$), or 
$i$ labels a vertical step and there is a $+$ in that row. In the first case, let
$D' \lessdot D$ be the Type 1 cover relation obtained by 
deleting the box containing the lowest $+$ in column of $D$ labeled by $i$ (and all boxes below it).  
In the second case, let
$D' \lessdot D$ be the Type 2  cover relation obtained by 
replacing the unique $+$ in the row labeled by $i$ with a $0$. In either case we have $\sigma(D') = \sigma'$.
\end{pf}

We now introduce the positroids corresponding to $\overline{\mathcal{D}_{n,k,1}}$.
\begin{defn}\label{positroids_m=1}
Let $\mathcal{M}_{n,k,1}$ be the set of matroids of rank $k$ with ground set $[n]$ which are direct sums $M_1\oplus\cdots\oplus M_{n-k}$, where $E_1\sqcup\cdots\sqcup E_{n-k}$ is a partition of $[n]$ into nonempty intervals, and $M_j$ is the uniform matroid of rank $|E_j|-1$ with ground set $E_j$, for $j\in [n-k]$. Let $\overline{\mathcal{M}_{n,k,1}}$ be the order ideal of matroids of $\mathcal{M}_{n,k,1}$, i.e.\ the set of matroids $M'$ of rank $k$ with ground set $[n]$ such that $M' \le M$ for some $M\in\mathcal{M}_{n,k,1}$. (The partial order was defined in \cref{defn_weak_maps}.)

We also define $\mathcal{M}_{n,k,1}^* := \{M^* : M\in\mathcal{M}_{n,k,1}\}$ and $\overline{\mathcal{M}_{n,k,1}^*} := \{M^* : M\in\overline{\mathcal{M}_{n,k,1}}\}$, so that $\overline{\mathcal{M}_{n,k,1}^*}$ is the set of matroids $M'$ of rank $n-k$ with ground set $[n]$ such that $M'\le M$ for some $M\in\mathcal{M}_{n,k,1}^*$. Since taking duals commutes with direct sum \cite[Proposition 4.2.21]{Oxley}, $\mathcal{M}_{n,k,1}^*$ is the set of matroids of rank $k$ with ground set $[n]$ which are direct sums $M_1\oplus\cdots\oplus M_{n-k}$, where $E_1\sqcup\cdots\sqcup E_{n-k}$ is a partition of $[n]$ into nonempty intervals, and $M_j$ is the uniform matroid of rank $1$ with ground set $E_j$, for $j\in [n-k]$.
\end{defn}
The matroids in $\overline{\mathcal{M}_{n,k,1}}$ 
and $\overline{\mathcal{M}_{n,k,1}^*}$ are in fact positroids; see \cref{le-diagram_m=1}. (Alternatively, using \cref{positroid_downset_m=1}, we can easily construct a matrix representing any matroid in $\overline{\mathcal{M}_{n,k,1}^*}$.) This implies that $\overline{\mathcal{M}_{n,k,1}}$ indexes an order ideal 
in $Q_{k,n}$, the poset of cells of $\Gr_{k,n}^{\ge 0}$.
\begin{lem}\label{positroid_downset_m=1}
The matroids in  $\overline{\mathcal{M}_{n,k,1}}$ 
and $\overline{\mathcal{M}_{n,k,1}^*}$ 
correspond to pairs 
$(E_1\sqcup\cdots\sqcup E_{n-k}, C)$, where
\begin{itemize}
\item $E_1\sqcup\cdots\sqcup E_{n-k}$ is a partition of $[n]$ into nonempty intervals and $C \subseteq [n]$;
\item for all $j\in [n-k]$, $E_j\setminus C$ is nonempty; and
\item $\max(E_j)\notin C$ whenever $\max(E_j)\neq n$.
\end{itemize}

The matroid $M \in 
\overline{\mathcal{M}_{n,k,1}}$ associated to 
$(E_1\sqcup\cdots\sqcup E_{n-k}, C)$ is 
the direct sum $M_1\oplus\cdots\oplus M_{n-k}$, where 
$M_j$ ($j\in [n-k]$) is the matroid with ground set $E_j$ such that 
$E_j\cap C$ are coloops and the restriction of $M_j$ to $E_j\setminus C$ is a uniform matroid of rank $|E_j\setminus C| - 1$.

And the matroid $M^* \in 
\overline{\mathcal{M}_{n,k,1}^*}$ associated to 
$(E_1\sqcup\cdots\sqcup E_{n-k}, C)$ is 
the direct sum $M_1^*\oplus\cdots\oplus M^*_{n-k}$, where 
$M^*_j$ ($j\in [n-k]$) is the matroid with ground set $E_j$ 
such that $E_j\cap C$ are loops and the restriction of $M_j$ to $E_j\setminus C$ is a uniform matroid of rank $1$.
\end{lem}

\begin{pf}
The description of $\overline{\mathcal{M}_{n,k,1}}$ follows from the description of 
$\overline{\mathcal{M}_{n,k,1}^*}$, so we prove the latter. In considering $\overline{\mathcal{M}_{n,k,1}^*}$, we will use the following claim.
\begin{claim}
Let $M$ be a matroid with ground set $F$ and connected components $F_1, \dots, F_l$, and let $M'\le M$. Then each $F_j$ ($j\in [l]$) is a union of connected components of $M'$, and $M'|_{F_j}\le M|_{F_j}$.
\end{claim}

\begin{claimpf}
By \cite[Proposition 4.1.2]{Oxley} (see also \cite[Proposition 7.2]{ardila_rincon_williams}), elements $i,j$ of the ground set are in the same connected component of a matroid $N$ if and only if there exist bases $B, B'$ of $N$ with $B' = (B\setminus\{i\})\cup\{j\}$. It follows that if $i,j\in F$ are in the same connected component of $M'$, then they are in the same connected component of $M$. The fact that $M'|_{F_j}\le M|_{F_j}$ for $j\in [l]$ follows from \cref{defn_restriction}, as long as $M'|_{F_j}$ and $M|_{F_j}$ have the same rank. But this is true since $\sum_{j=1}^l\rank(M'|_{F_j}) = \rank(M') = \rank(M) = \sum_{j=1}^l\rank(M|_{F_j})$.
\end{claimpf}

Now note that if $M$ is a uniform matroid of rank $1$, then the matroids $M'$ satisfying $M'\le M$ are precisely all matroids of rank $1$ with the same ground set as $M$. Hence by the claim, the elements of $\overline{\mathcal{M}_{n,k,1}^*}$ are obtained precisely by taking a matroid $M_1\oplus\cdots\oplus M_{n-k}\in\mathcal{M}_{n,k,1}$ with ground set $E_1\sqcup\cdots\sqcup E_{n-k}$, and choosing some subset $C\subseteq [n]$ of the ground set, satisfying $E_j\setminus C\neq\emptyset$ for all $j\in [n-k]$, to turn into loops (i.e.\ we delete all bases with a nonempty intersection with $C$). The condition that $\max(E_j)\notin C$ if $\max(E_j)\neq n$ comes from our convention that a loop which appears `between' two intervals is associated to the interval on its right.
\end{pf}

\begin{rmk}
Using \cref{positroid_downset_m=1}, one can write down the generating function for the stratification of $\mathcal{B}_{n,k,1}(W)$ with respect to dimension; see \cref{generating_function_m=1}. We will 
give a different proof of \cref{generating_function_m=1}, using the fact that $\mathcal{B}_{n,k,1}(W)$ is isomorphic to the bounded complex of a generic hyperplane arrangement, whose rank generating function is known.
\end{rmk}

We now give a poset isomorphism $\DM: \overline{\mathcal{D}_{n,k,1}} \to \overline{\mathcal{M}_{n,k,1}}$.

\begin{figure}[ht]
\begin{center}
 \includegraphics[width=6.4in]{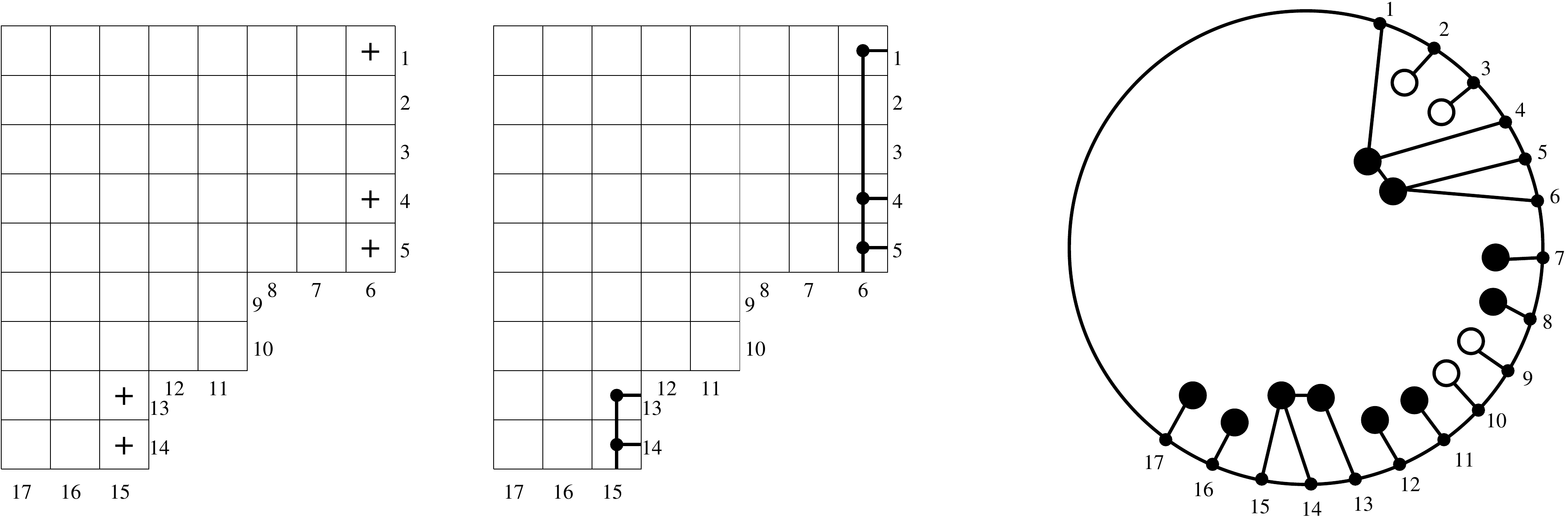}
 \caption{Going from a \Le -diagram in $\overline{\mathcal{D}_{n,k,1}}$ to 
the corresponding plabic graph, which in turn determines the positroid in $\overline{\mathcal{M}_{n,k,1}}$.}
 \label{fig:LeToPlabic}
\end{center}
\end{figure}

\begin{defn}\label{le-diagram_m=1}
Given $D\in\overline{\mathcal{D}_{n,k,1}}$, 
we label the southeast border of $D$
by the numbers $1$ through $n$.  Let the labels of 
the horizontal steps be denoted by 
$h_1,\dots,h_{n-k}$.  Then set $E_1 := \{1, 2, \dots, h_1\}$, 
$E_j := \{h_{j-1} + 1, h_{j-1} + 2, \dots, h_j\}$ for $1 <  j < n-k$,
and $E_{n-k}:= \{h_{n-k-1}+1, h_{n-k-1}+2,\dots,n\}$.   Let $C$ be the set of labels of all vertical steps
indexing rows of $D$ with no $+$'s.
Then $(E_1 \sqcup \dots \sqcup E_{n-k},C)$ determines a positroid in $\overline{\mathcal{M}_{n,k,1}}$ as in 
\cref{positroid_downset_m=1}, which we denote by $M(D)$.  By inspection, we see that the map
$\DM: D \mapsto M(D)$ is a bijection, and we will denote the inverse by 
$\DM^{-1}:M \mapsto D(M)$.
\end{defn}
We observe that $M(D)$ is precisely the positroid of $D$ defined in \cref{prop:perf}. In particular, the elements of $\overline{\mathcal{M}_{n,k,1}}$ and $\overline{\mathcal{M}_{n,k,1}^*}$ are all positroids. See \cref{fig:LeToPlabic}.
The white lollipops correspond to coloops and the black lollipops
correspond to loops.  By considering perfect
orientations, it is easy to see that every component of the graph 
which is not a white lollipop gives rise to a uniform matroid of corank $1$.

\begin{lem}\label{le_positroid_isomorphism}
The map $\DM$ from \cref{le-diagram_m=1} is an isomorphism of posets. In particular, $\overline{\mathcal{D}_{n,k,1}}$ can be identified with an order ideal of $Q_{k,n}$, the poset of cells of $\Gr_{k,n}^{\ge 0}$.
\end{lem}

\begin{pf}
By inspection, $\DM$ is a bijection and a poset homomorphism. To see that its inverse is a poset homomorphism, let $M'\lessdot M$ be a cover relation in $\overline{\mathcal{M}_{n,k,1}}$, and take $D\in\overline{\mathcal{D}_{n,k,1}}$ with $M = M(D)$. Let $E_1\sqcup\cdots\sqcup E_{n-k}$ correspond to $M$ as in \cref{positroid_downset_m=1}. Then we obtain $M'$ from $M$ by taking some interval $E_j$ which contains at least $2$ elements which are not coloops, and turning some $i\in E_j$ into a coloop. If $i$ is the greatest element of $E_j$ not already a coloop, then $i$ labels a horizontal step in the southeast border of $D$. In this case, let $D' \lessdot D$ be the Type 1 cover relation given by deleting the box containing the lowest $+$ in column of $D$ labeled by $i$ (and all boxes below it). Otherwise $i$ labels a vertical step in the southeast border of $D$ whose row contains a $+$; let $D' \lessdot D$ be the Type 2 cover relation given by replacing this $+$ with a $0$. Then $M(D') = M'$.
\end{pf}

By composing the bijections from Definitions \ref{def:Lesign} and \ref{le-diagram_m=1}, we obtain the following result.

\begin{cor}\label{cor:posetiso}
The map $\DS \circ {\DM}^{-1}: \overline{\mathcal{M}_{n,k,1}}\to\overline{\mathcal{\Sign}_{n,k,1}}, M \mapsto \sigma(M)$ is an isomorphism of posets, 
where we define $\sigma(M) := \sigma(D)$ for $D\in\overline{\mathcal{D}_{n,k,1}}$ with $M = M(D)$. We can compute $\sigma(M)$ from $M$ as follows.
If $M\in {\mathcal{M}_{n,k,1}}$ is the direct sum
$M_1 \oplus \dots \oplus M_{n-k}$ with ground set 
$E_1 \sqcup \dots \sqcup E_{n-k}$, then $\sigma(M)$ is uniquely determined by the following two properties:
\begin{itemize}
\item $\sigma(M)_1 = +$;
\item $\sigma(M)_i = \sigma(M)_{i+1}$ if and only if $i$ and 
$i+1$ are in different blocks of
$E_1 \sqcup \dots \sqcup E_{n-k}$.
\end{itemize}
And if $M\in \overline{\mathcal{M}_{n,k,1}}$ is the direct sum
$M_1 \oplus \dots \oplus M_{n-k}$ with ground set
$E_1 \sqcup \dots \sqcup E_{n-k}$ and coloops at $C$,
then $\sigma(M)$ is given exactly as above, except that $\sigma(M)_i=0$ for each $i\in C$.
\end{cor}
For example, if $M \in \mathcal{M}_{8,5,1}$ is associated to 
$E_1\sqcup E_2 \sqcup E_3$ with 
$E_1 = \{1,2,3\}, E_2 = \{4,5,6\}, E_3 = \{7,8\}$, then $\sigma(M)$ equals $(+,-,+,+,-,+,+,-)$.

\begin{defn}\label{defphi}
Let $\mathcal{S} := \bigsqcup_{M\in\overline{\mathcal{M}_{n,k,1}}}S_M = \overline{\bigsqcup_{M\in\mathcal{M}_{n,k,1}}S_M}$ be the subcomplex of $\Gr_{k,n}^{\ge 0}$ corresponding to $\mathcal{M}_{n,k,1}$. We define a map
(\emph{cf. \cref{defn_B_amplituhedron}})
$$
\phi_{W}:\mathcal{S}\to\mathcal{B}_{n,k,1}(W), \quad V\mapsto V^\perp\cap W.
$$
\end{defn}

\begin{prop}\label{prop:induce}
If $V\in S_M$ for $M\in \overline{\mathcal{M}_{n,k,1}}$, then 
$V^\perp \cap W \in \mathcal{B}_{\sigma(M)}(W)$.
In other words, the map $\phi_{W}$ from \cref{defphi} induces
the map $M \mapsto \sigma(M)$ on strata.
\end{prop}
Technically $\sigma(M)$ is an element of $\overline{\Sign_{n,k,1}}$, while
 $\mathcal{B}_{n,k,1}(W)$ is stratified by $\overline{\PSign_{n,k,1}}$. By \cref{Sign_PSign}(i) the map $\overline{\Sign_{n,k,1}}\to\overline{\PSign_{n,k,1}}, \sigma\mapsto\sigma$ is a poset isomorphism, so we need not concern ourselves with this distinction.
\begin{pf}
We first consider the case that $M\in\mathcal{M}_{n,k,1}$. 
Let us describe the sign vectors of $V^\perp$ for $V\in S_M$. 
Write $M = M_1\oplus\cdots\oplus M_{n-k}$, where
$E_1\sqcup\cdots\sqcup E_{n-k}$
 is a partition of $[n]$ into nonempty intervals,  and $M_j$ is the uniform matroid of rank $|E_j| -1$ with ground set $E_j$, for $j\in [n-k]$.  
By \cref{positive_covectors}(ii), if $V_j\in S_{M_j}$ then $\sign(V_j^\perp) = \{\sigma\in\{0,+,-\}^{E_j} : \var(\sigma)\ge |E_j| - 1\}\cup\{0\}$, i.e.\
$$
\sign(V_j^\perp) = \{0, (+,-,+,-,\dots), (-,+,-,+,\dots)\}.
$$
Hence for $V\in S_M$, we have $\sigma\in\sign(V^\perp)$ if and only if $\sigma|_{E_j}$ equals $0$ or strictly alternates in sign, for all $j\in [n-k]$. 

Recall from \cref{B_m=1} that 
$$\mathcal{B}_{n,k,1}(W) = \{w\in W\setminus\{0\}: \overline{\var}(w) = k\}\subseteq\mathbb{P}(W).$$
Note that there is a unique nonzero $\sigma\in\sign(V^\perp)$ (modulo multiplication by $\pm 1$) with $\overline{\var}(\sigma) = k$: $\sigma$ has no zero components, and $\sigma_i = \sigma_{i+1}$ if and only if $i$ and $i+1$ are in different 
blocks  of $E_1\sqcup\cdots\sqcup E_{n-k}$, for all $i\in [n-1]$. 
This is precisely the definition of $\sigma(M)$.
Therefore we must have $V^\perp \cap W \in \mathcal{B}_{\sigma(M)}(W)$.

Now we consider the general case of sign vectors of $V'^\perp$, where $M'\in\overline{\mathcal{M}_{n,k,1}}$ and $V'\in S_{M'}$. We have $M'\le M$ for some $M\in\mathcal{M}_{n,k,1}$, 
and $M'$ is obtained from $M$ by making some subset $C\subseteq [n]$ of the ground set coloops.
Thus the sign vectors of $V'^\perp$ are precisely obtained from those of $V^\perp$ by setting the components indexed by
$C$ to $0$. In particular, we again have a unique nonzero $\sigma'\in\sign(V'^\perp)$ (modulo multiplication by $\pm 1$) with $\overline{\var}(\sigma') = k$, which we obtain from $\sigma(M)$ by setting the components $C$ to $0$. Then $\sigma' = \sigma(M')$,
and therefore ${V'}^\perp \cap W \in \mathcal{B}_{\sigma(M')}(W)$.
\end{pf}

\begin{thm}\label{triangulation_m=1}
The map $\phi_{W}$ from \cref{defphi} is a homeomorphism
which induces 
a poset isomorphism on the stratifications of $\mathcal{S}$ and $\mathcal{B}_{n,k,1}(W)$. 
\end{thm}

\begin{pf}
We know from \cref{prop:induce} that $\phi_{W}$ induces a poset isomorphism
on the strata. To show that $\phi_{W}$ is a bijection, we 
construct the inverse map $\phi_{W}^{-1}:\mathcal{B}_{n,k,1}\to\mathcal{S}$ as follows. Given an element of $\mathcal{B}_{n,k,1}(W)$ spanned by $w\in W\setminus\{0\}$ (so $\overline{\var}(w) = k$), let $\sigma$ be either $\sign(w)$ or $-\sign(w)$, whichever is in $\overline{\Sign_{n,k,1}}$. Also let $M\in\overline{\mathcal{M}_{n,k,1}}$ be the positroid such that $\sigma(M)=\sigma$, corresponding to $(E_1\sqcup\cdots\sqcup E_{n-k}, C)$ in \cref{positroid_downset_m=1}. If $V\in\mathcal{S}$ with $\phi_{W}(V) = w$, then $V\in S_M$ because $\phi_{W}$ induces a poset isomorphism
on the strata. Since $E_j\setminus C\neq\emptyset$ for all $j\in [n-k]$, $\sigma$ is nonzero when restricted to any interval $E_j$. Hence the unique $V\in\mathcal{S}$ with $\phi_{W}(V) = w$ is determined by the conditions
$$
V^\perp|_{E_j} = \spn(w|_{E_j}) \qquad \text{ for all }j\in [n-k].
$$
Explicitly, $V^\perp$ has the basis $w^{(1)}, \dots, w^{(n-k)}$, where $w^{(i)}|_{E_j} = \delta_{i,j}w|_{E_j}$ for $i,j\in [n-k]$. 
Thus $\phi_{W}$ is invertible, with an inverse which is piecewise polynomial
(each stratum is a domain of polynomiality).
The map $\phi_{W}$ is continuous, and therefore a homeomorphism.
\end{pf}

\section{\texorpdfstring{$\mathcal{A}_{n,k,1}$}{A(n,k,1)} as the bounded complex of a cyclic hyperplane arrangement}\label{sec_cyclic_arrangement}

\noindent We show that the $m=1$ amplituhedron $\mathcal{B}_{n,k,1}(W)$ (or $\mathcal{A}_{n,k,1}(Z)$) is homeomorphic to the complex of bounded faces of a {\itshape cyclic hyperplane arrangement} of $n$ hyperplanes in $\mathbb{R}^k$. It then follows from a result of Dong \cite{Dong} that it is homeomorphic to a ball. This story is somewhat analogous to that of $k=1$ amplituhedra
$\mathcal{A}_{n,1,m}$, which are {\itshape cyclic polytopes} with $n$ vertices in $\mathbb{P}^m$.\footnote{In fact, the 
$k=1$ amplituhedra $\mathcal{A}_{n,1,m}$ 
are precisely the {\itshape alternating polytopes} of dimension $m$
with $n$ vertices in $\mathbb{P}^m$, as follows from work of  Sturmfels \cite{Sturmfels}. Alternating polytopes are  cyclic polytopes which have the additional property that every induced subpolytope is also cyclic. See \cite[pp.\ 396-397]{bjorner_las_vergnas_sturmfels_white_ziegler} for an example of a cyclic polytope which is not alternating.} (We do not know whether this is a coincidence, or a specific instance of some form of duality for amplituhedra.)

Cyclic hyperplane arrangements have been studied by Shannon \cite{shannon}, Ziegler \cite{ziegler}, Ram\'{i}rez Alfons\'{i}n \cite{ramirez_alfonsin}, and Forge and Ram\'{i}rez Alfons\'{i}n \cite{forge_ramirez_alfonsin}. 
For an introduction to hyperplane arrangements, see \cite{stanley_hyperplanes}.

\begin{rmk}
In the literature, a {\itshape cyclic hyperplane arrangement} of $n$ hyperplanes in $\mathbb{R}^k$ is usually defined to be an arrangement with hyperplanes
\begin{align}\label{cyclic_arrangement_standard}
\{x\in\mathbb{R}^k : t_ix_1 + t_i^2x_2 + \cdots + t_i^kx_k + 1 = 0\} \qquad (i\in [n]),
\end{align}
where $0 < t_1 < \cdots < t_n$. We will need to consider more general hyperplane arrangements, whose hyperplanes are of the form
\begin{align}\label{cyclic_arrangement_general}
\{x\in\mathbb{R}^k : A_{i,1}x_1 + \cdots + A_{i,k}x_k + A_{i,0} = 0\} \qquad (i\in [n]),
\end{align}
such that
$$
\colspan((A_{i,j})_{1 \le i \le n, 1 \le j \le k})\in\Gr_{k,n}^{>0} \qquad\text{ and }\qquad \colspan((A_{i,j})_{1 \le i \le n, 0 \le j \le k})\in\Gr_{k+1,n}^{>0}.
$$
(We will require that this latter subspace is $W$.) The hyperplane arrangement \eqref{cyclic_arrangement_standard} is of this form by Vandermonde's identity. \cref{cyclic_arrangement} implies that all hyperplane arrangements of the form \eqref{cyclic_arrangement_general} are isomorphic, so we will not be concerned with the distinction between \eqref{cyclic_arrangement_standard} and \eqref{cyclic_arrangement_general}. This is analogous to the situation for cyclic polytopes (see \cite{Sturmfels}).
\end{rmk}

\begin{defn}\label{defn_generic_arrangement}
An arrangement $\{H_1, \dots, H_n\}$ of hyperplanes in $\mathbb{R}^k$ is called 
{\itshape generic} if for all $I\subseteq [n]$, we have $\dim(\bigcap_{i\in I}H_i) = n-|I|$ if $|I| \le k$, and $\bigcap_{i\in I}H_i = \emptyset$ if $|I| > k$.
\end{defn}

\begin{rmk}
Cyclic polytopes have many faces of each dimension, in the sense of the upper bound theorem of McMullen \cite{mcmullen} and Stanley \cite{stanley_75}. An analogous property of cyclic hyperplane arrangements is that they have few simplicial
faces of each dimension, in the sense of Shannon \cite{shannon}. Note that it does not make sense to look at the total number of faces, because the number of faces of a given dimension of a generic hyperplane arrangement depends only on its dimension and the number of hyperplanes \cite{buck}.
\end{rmk}

In order to define our hyperplane arrangement, we will use the following convention.
\begin{defn}\label{w0_convention}
Suppose that $V\in\Gr_{k,n}^{>0}$ and $w\in\mathbb{R}^n\setminus V$ such that $V+w\in\Gr_{k+1,n}^{>0}$. Let $u\in V^\perp$ be the orthogonal projection of $w$ to $V^\perp$, i.e.\ $w - u\in V$. Since $u\neq 0$, \cref{gantmakher_krein}(ii)
implies that $\var(u)\geq k$.  But since $u$ also lies in $V+w\in \Gr_{k+1,n}^{>0}$,
\cref{gantmakher_krein}(ii)
implies that $\overline{\var}(u) \leq k$.  Therefore
$\var(u) = \overline{\var}(u) = k$. But $\var(u)=\overline{\var}(u)$ implies that 
the first component $u_1$ of $u$ is nonzero. 
We call $w$ {\itshape positively oriented with respect to $V$} if $u_1 > 0$, and 
{\itshape negatively oriented} if $u_1 < 0$.
\end{defn}

\begin{eg}
Let $V\in\Gr_{1,3}^{>0}$ be the span of $(1,1,1)$ (so $n=3$, $k=1$), and $w := (1,2,4)$, so that $V+w\in\Gr_{2,3}^{>0}$. Then the projection of $w$ to $V^\perp$ is $u := (-\frac{4}{3}, -\frac{1}{3}, \frac{5}{3})$, since $w - u = \frac{7}{3}(1,1,1)\in V$. Since $u_1 = -\frac{4}{3} < 0$, $w$ is negatively oriented with respect to $V$.
\end{eg}

Let us show that the cyclic hyperplane arrangements defined in \eqref{cyclic_arrangement_standard} satisfy this positive orientation property. 
Although we will not need \cref{Vandermonde_positively_oriented} in what follows,
it will mean that our future characterization of face labels 
of $\mathcal{H}^W$ in \cref{cyclic_arrangement_faces} applies to 
``classical'' cyclic hyperplane arrangements \eqref{cyclic_arrangement_standard}. (This is due to the assumption in \cref{defn_cyclic_arrangement} that $w^{(0)}$ is positively oriented with respect to $V$.)
\begin{prop}\label{Vandermonde_positively_oriented}
Let $0 < t_1 < \cdots < t_n$, and $V\in\Gr_{k,n}^{>0}$ be the span of the vectors $(t_1^j, \dots, t_n^j)$ for $1 \le j \le k$. Then $w=(1, \dots, 1)$ is positively oriented with respect to $V$.
\end{prop}

\begin{pf}
We define the $n\times (k+1)$ matrix $A$ with entries
$A_{i,j} := t_i^j$ ($1 \le i \le n, 0 \le j \le k$). 
Then $A$ is a {\itshape totally positive matrix}, i.e.\ all its minors are positive. Indeed, for $I = \{i_1 < \cdots < i_l\}\subseteq [n]$ and $J = \{j_1 < \cdots < j_l\}\subseteq\{0,1,\dots, k\}$, the classical definition of Schur functions implies that
$$
\det(A_{I,J}) = 
s_{(j_l - l + 1, j_{l-1} - l + 2, \dots, j_1)}(t_{i_1}, \dots, t_{i_l})
\prod_{r,s\in [l], \; r < s}(t_{i_s} - t_{i_r}).
$$
Since Schur functions are monomial-positive,
$\det(A_{I,J}) >0$.

Because $V+w$ is the column span of $A$, and $V$ is the span of the last
$k$ columns of $A$, we get that 
$V\in \Gr_{k,n}^{>0}$
and $V + w \in \Gr_{k+1,n}^{>0}$. As in \cref{w0_convention}, we let $u\in V^\perp$ be the projection of $w$ to $V^\perp$. That is, $w - u\in V$ and $\var(u) = \overline{\var}(u) = k$. 
We must show that $u_1 > 0$. 
We will use the following properties of totally positive
$n \times (k+1)$ matrices $A$:
\begin{itemize}
\item $\var(Ax)\le\var(x)$ for all $x\in\mathbb{R}^{k+1}$ \cite{schoenberg};
\item if $\var(Ax) = \var(x)$, then the first nonzero components of $Ax$ and $x$ have the same sign \cite[Theorem V.5]{gantmakher_krein_50}.
\end{itemize}
Since $w - u\in V$, we can write $u = Ax$ for some $x\in\mathbb{R}^{k+1}$ with $x_1 = 1$. And because $\var(u) = k$, we have $\var(x) = k$, whence the first nonzero component of $u$ is positive.
\end{pf}

We will also need the following result of 
Rietsch.\footnote{Rietsch
\cite{rietsch_private} in fact proved that the totally nonnegative part $\Fl_n^{\ge 0}$ of the complete flag variety 
(as defined by Lusztig \cite[Section 8]{lusztig}) projects surjectively onto $\Gr_{k,n}^{\ge 0}$, and that the Lusztig-Rietsch 
stratification of $\Fl_n^{\ge 0}$ projects onto Postnikov's stratification of $\Gr_{k,n}^{\ge 0}$. In particular, given $V\in\Gr_{k,n}^{>0}$, there exists a complete flag $V_0 \subset V_1 \subset \dots \subset V_n$ 
in the totally positive part $\Fl_n^{>0}$ of $\Fl_n^{\ge 0}$ with $V_k = V$. This immediately implies \cref{TP_restriction}, because if $V_0 \subset V_1 \subset \dots \subset V_n\in\Fl_n^{>0}$ then $V_j\in\Gr_{j,n}^{>0}$ 
for $0 \le j \le n$. (See \cite[Corollary 7.2]{tsukerman_williams} for a related result.)
An alternative proof of Rietsch's result was given by Talaska and Williams \cite[Theorem 6.6]{talaska_williams}, by relating Postnikov's parameterizations of cells of $\Gr_{k,n}^{\ge 0}$ \cite[Theorem 6.5]{postnikov} to Marsh and Rietsch's parameterizations of the Lusztig-Rietsch cells \cite{marsh_rietsch}. 
A direct proof of \cref{TP_restriction} using similar tools was given in \cite[Lemma 15.6]{lam}.}

\begin{lem}[\cite{rietsch_private}]\label{TP_restriction}
If $W\in\Gr_{k+m,n}^{>0}$, where $m\ge 0$, then $W$ contains a subspace in $\Gr_{k,n}^{>0}$.
\end{lem}

We now define our hyperplane arrangement. A hyperplane arrangement $\mathcal{H}$ partitions its ambient space into {\itshape faces}; maximal faces (equivalently, connected components of the complement of $\mathcal{H}$) are called {\itshape regions}.
\begin{defn}\label{defn_cyclic_arrangement}
Given $W\in\Gr_{k+1,n}^{>0}$, 
we use \cref{TP_restriction} 
to choose vectors $w^{(1)}, \dots, w^{(k)}\in W$ 
such that $V:=\spn(w^{(1)}, \dots, w^{(k)})\in\Gr_{k,n}^{>0}$. We extend $\{w^{(1)}, \dots, w^{(k)}\}$ to a basis $\{w^{(0)}, w^{(1)}, \dots, w^{(k)}\}$ of $W$; after replacing $w^{(0)}$ with $-w^{(0)}$ if necessary, we assume that $w^{(0)}$ is positively oriented with respect to $V$ (see \cref{w0_convention}).

We let $\mathcal{H}^W$ be the hyperplane arrangement in $\mathbb{R}^k$ with hyperplanes
$$
H_i := \{x\in\mathbb{R}^k : w^{(1)}_ix_1 + \cdots + w^{(k)}_ix_k + w^{(0)}_i = 0\} \text{ for }i\in [n].
$$
Note that $\mathcal{H}^W$ is generic by the first three sentences of the proof of \cite[Proposition 5.13]{stanley_hyperplanes}.\footnote{We thank Richard Stanley for pointing out this argument to us.} Also note that $\mathcal{H}^W$ depends not only on $W$ but also on our choice of basis of $W$.

Given $w\in W$, we let $\langle w \rangle \in \mathbb{P}(W)$ 
denote the line spanned by $w$.
We define the maps
\begin{align*}
& \Psi_{\mathcal{H}^W} : \mathbb{R}^k\to\mathbb{P}(W),\hspace*{-48pt} && x\mapsto 
\langle x_1w^{(1)} + \cdots + x_kw^{(k)} + w^{(0)} \rangle, \\
& \psi_{\mathcal{H}^W} : \mathbb{R}^k\to\{0,+,-\}^n,\hspace*{-48pt} && x\mapsto\sign(x_1w^{(1)} + \cdots + x_kw^{(k)} + w^{(0)}).
\end{align*}
Note that the faces of ${\mathcal{H}^W}$ are precisely the nonempty fibers of $\psi_{\mathcal{H}^W}$. If $\sigma\in\{0,+,-\}^n$ has a nonempty preimage under $\psi_{\mathcal{H}^W}$, we call this fiber the face of ${\mathcal{H}^W}$ {\itshape labeled by $\sigma$}. 
When we identify faces with labels in this way, the face poset of 
${\mathcal{H}^W}$ 
is an induced subposet of the sign vectors $\{0,+,-\}^n$ (see \cref{defn_sign_vector_partial_order}).

Finally, we let $B({\mathcal{H}^W})$ be the subcomplex of bounded faces of ${\mathcal{H}^W}$. We denote the set of sign vectors which label the bounded faces of ${\mathcal{H}^W}$ by $\mathcal{V}_{B({\mathcal{H}^W})}$.
\end{defn}

\begin{rmk}
By \cref{Vandermonde_positively_oriented},
the classical cyclic hyperplane arrangement \eqref{cyclic_arrangement_standard} is an example of such an ${\mathcal{H}^W}$ from 
\cref{defn_cyclic_arrangement}.
\end{rmk}

We will show that $\Psi_{\mathcal{H}^W}$ gives a homeomorphism from $B({\mathcal{H}^W})$ to $\mathcal{B}_{n,k,1}(W)$ (\cref{cyclic_arrangement}). The key to the proof is establishing that $\mathcal{V}_{B({\mathcal{H}^W})} = \overline{\Sign_{n,k,1}}$, which we do in \cref{cyclic_arrangement_faces}. 
(Recall that $\overline{\Sign_{n,k,1}}$ is the set of nonzero sign vectors $\sigma\in\{0,+,-\}^n$ such that $\overline{\var}(\sigma) = k$, and if $i\in [n]$ indexes the first nonzero component of $\sigma$, then $\sigma_i = (-1)^{i-1}$. See \cref{defn_Sign}.) 

In what follows, we fix $W\in\Gr_{k+1,n}^{>0}$, as well as a basis $\{w^{(0)}, w^{(1)}, \dots, w^{(k)}\}$ of $W$ and a corresponding hyperplane arrangement ${\mathcal{H}^W}$, as in \cref{defn_cyclic_arrangement}.
\begin{lem}\label{unbounded_faces_lemma} ~\\
(i) If $\sigma\in\{+,-\}^n$ satisfies $\var(\sigma) \le k-1$, then $\sigma$ labels an unbounded region of ${\mathcal{H}^W}$. \\
(ii) If $\sigma\in\{0,+,-\}^n$ labels an unbounded face of ${\mathcal{H}^W}$, then $\overline{\var}(\sigma) \le k-1$.
\end{lem}

\begin{pf}
(i) Suppose that $\sigma\in\{+,-\}^n$ with $\var(\sigma) \le k-1$.
Since $\spn(w^{(1)}, \dots, w^{(k)})\in\Gr_{k,n}^{>0}$,  
\cref{positive_covectors}(i) implies that we 
can write $\sigma = \sign(w)$ for some $w\in\spn(w^{(1)}, \dots, w^{(k)})$. 
Then for all $t > 0$ sufficiently large we have $\sign(tw + w^{(0)}) = \sigma$, whence $\psi_{\mathcal{H}^W}^{-1}(\sigma)$ is unbounded. 

(ii) Given $\sigma\in\{0,+,-\}^n$ such that $\psi_{\mathcal{H}^W}^{-1}(\sigma)$ is an unbounded face of ${\mathcal{H}^W}$, take a sequence $(x^{(j)})_{j\in\mathbb{N}}$ in $\psi_{\mathcal{H}^W}^{-1}(\sigma)$ 
with $\lim_{j\to\infty}\|x^{(j)}\| = \infty$.

\begin{claim}
$\lim_{j\to\infty}\|x^{(j)}_1w^{(1)} + \cdots + x^{(j)}_kw^{(k)} + w^{(0)}\| = \infty$.
\end{claim}

\begin{claimpf}
Let $c\in\mathbb{R}$ be the minimum value of $\|x_1w^{(1)} + \cdots + x_kw^{(k)}\|$ on the compact set $\{x\in\mathbb{R}^k : \|x\| = 1\}$. Since $w^{(1)}, \dots, w^{(k)}$ are linearly independent, we have $c > 0$. Hence by the triangle inequality,
$$
\|x^{(j)}_1w^{(1)} + \cdots + x^{(j)}_kw^{(k)} + w^{(0)}\| \ge \|x^{(j)}_1w^{(1)} + \cdots + x^{(j)}_kw^{(k)}\| - \|w^{(0)}\| \ge c\|x^{(j)}\| - \|w^{(0)}\|
$$
for $j\in\mathbb{N}$.
The claim follows when we send $j$ to $\infty$.
\end{claimpf}

By the claim, the following sequence is well defined for $j$ sufficiently large:
$$
\left(\frac{x^{(j)}_1w^{(1)} + \cdots + x^{(j)}_kw^{(k)} + w^{(0)}}{\|x^{(j)}_1w^{(1)} + \cdots + x^{(j)}_kw^{(k)} + w^{(0)}\|}\right)_{j\in\mathbb{N}}.
$$
Since $\{x\in\mathbb{R}^n : \|x\| = 1\}$ is compact, this sequence has a convergent subsequence;
we let $v\in\mathbb{R}^n$ denote one of its limit points.
By the claim, we have $v\in\spn(w^{(1)}, \dots, w^{(k)})\setminus\{0\}$. We get $\overline{\var}(\sigma)\le\overline{\var}(v)\le k-1$, where the first inequality holds since $\sign(v)$ is obtained from $\sign(\sigma)$ by possibly 
setting some nonzero components to zero, and the second inequality follows from \cref{gantmakher_krein}(ii).
\end{pf}

\begin{lem}\label{bounded_faces_lemma}
The map $\mathcal{V}_{B({\mathcal{H}^W})}\to\overline{\PSign_{n,k,1}}, \sigma\mapsto\sigma$ is a poset isomorphism.
In other words, the face poset of the bounded faces of ${\mathcal{H}^W}$
is isomorphic to the face poset of $\mathcal{B}_{n,k,1}(W)$.
\end{lem}

\begin{pf}
First let us show that the map $\mathcal{V}_{B({\mathcal{H}^W})}\to\overline{\PSign_{n,k,1}}$ is well defined, i.e.\ given $\sigma\in\mathcal{V}_{B({\mathcal{H}^W})}$, we have $\overline{\var}(\sigma) = k$. Since $B({\mathcal{H}^W})$ equals the closure of the union of the bounded regions of ${\mathcal{H}^W}$ 
\cite{Dong} (see also \cite[Chapter 1, Exercise 7(d)]{stanley_hyperplanes}), the face labeled by $\sigma$ is contained in the closure of some bounded region of ${\mathcal{H}^W}$, labeled by, say, $\tau\in\{+,-\}^n$. Therefore $\tau\ge\sigma$, where the partial order is the one on $\{0,+,-\}^n$ from 
\cref{defn_sign_vector_partial_order}. 
By \cref{unbounded_faces_lemma}(i) we have $\var(\tau)\ge k$, and so $\overline{\var}(\sigma)\ge\var(\tau)\ge k$. 
But by
\cref{gantmakher_krein}(ii) we 
also have $\overline{\var}(\sigma)\le k$, so 
$\overline{\var}(\sigma)= k$, as desired.

Therefore the map $\mathcal{V}_{B({\mathcal{H}^W})}\to\overline{\PSign_{n,k,1}}$ is a poset homomorphism. To see that it is surjective, note that if $\sigma\in\{0,+,-\}^n\setminus\{0\}$ satisfies $\overline{\var}(\sigma)=k$, then $\sigma = \sign(w)$ for some $w\in W$ by \cref{positive_covectors}(i). When we write $w$ in terms of the basis $w^{(0)}, w^{(1)}, \dots, w^{(k)}$, the coefficient of $w^{(0)}$ is nonzero. (Otherwise $w\in\spn(w^{(1)}, \dots, w^{(k)})\in\Gr_{k,n}^{>0}$, implying $\overline{\var}(w) \leq k-1$ by \cref{gantmakher_krein}(ii).)  Rescaling $w$ by a positive real number so that this coefficient is $\pm 1$, we see 
that we can write 
$$w = x_1 w^{(1)} + \dots + x_k w^{(k)} \pm w^{(0)}$$ for some $x\in\mathbb{R}^k$.
Therefore
either $\sigma$ or $-\sigma$ labels a face of ${\mathcal{H}^W}$, and such a face is bounded by \cref{unbounded_faces_lemma}(ii).

It remains to show that the map $\mathcal{V}_{B({\mathcal{H}^W})}\to\overline{\PSign_{n,k,1}}$ is injective and that its inverse is a poset homomorphism. It suffices to prove that there do not exist $\sigma,\tau\in\{0,+,-\}^n$ labeling bounded faces of ${\mathcal{H}^W}$ such that $\sigma\le -\tau$. Suppose otherwise that there exist such $\sigma,\tau$. Take $x,y\in\mathbb{R}^k$ with $\psi_{\mathcal{H}^W}(x) = \sigma$ and $\psi_{\mathcal{H}^W}(y) = \tau$. Subtracting,
we get $\sign((x_1-y_1)w^{(1)} + \cdots + (x_k-y_k)w^{(k)}) = -\tau$.
Since 
$\spn(w^{(1)}, \dots, w^{(k)})\in\Gr_{k,n}^{>0}$, 
 \cref{gantmakher_krein}(ii) 
implies that 
$\overline{\var}(\tau)\le k-1$. But we showed in the first paragraph
 that $\overline{\var}(\tau) = k$.
\end{pf}

\begin{prop}[The face labels of ${\mathcal{H}^W}$]\label{cyclic_arrangement_faces} ~\\
(i) The labels of the bounded faces of ${\mathcal{H}^W}$ are precisely $\overline{\Sign_{n,k,1}}$, i.e.\ $\mathcal{V}_{B({\mathcal{H}^W})} = \overline{\Sign_{n,k,1}}$. \\
(ii) The labels of the unbounded faces of ${\mathcal{H}^W}$ are precisely $\sigma\in\{0,+,-\}^n$ with $\overline{\var}(\sigma)\le k-1$.
\end{prop}

\begin{eg}\label{hyperplane_arrangement_example_n=5}
Let $n := 5, k := 2, m := 1$, and $W\in\Gr_{3,5}^{>0}$ have the basis
$$
w^{(0)} := (-1, -1, -1, -1, -1), \quad w^{(1)} := (0, 1, 2, 3 ,4), \quad w^{(2)} := (10, 6, 3, 1, 0).
$$
Note that $V:=\spn(w^{(1)}, w^{(2)})\in\Gr_{2,5}^{>0}$. 
The projection $u$ of $w^{(0)}$ to $V^\perp$ is
$$
u := w^{(0)} + \textstyle\frac{1}{831}(232w^{(1)} + 90w^{(2)}) = \textstyle\frac{1}{831}(69, -59, -97, -45, 97).
$$
Since $u_1 > 0$, by \cref{w0_convention} $w^{(0)}$ is positively oriented with respect to $V$.

The hyperplane arrangement ${\mathcal{H}^W}$ from \cref{defn_cyclic_arrangement} consists of $5$ lines in $\mathbb{R}^2$:
$$\text{$\ell_1$: }  10y = 1, \qquad
\text{$\ell_2$: }  x + 6y = 1, \qquad
\text{$\ell_3$: }  2x + 3y = 1, \qquad
\text{$\ell_4$: }  3x + y = 1, \qquad
\text{$\ell_5$: }  4x  = 1.$$
See \cref{hyperplane_arrangement_sign_vectors_5_2_1} and \cref{hyperplane_arrangement_le_diagrams_5_2_1} for $B({\mathcal{H}^W})$ labeled by sign vectors and \Le -diagrams, respectively. (The positive side of each line is above and to the right of it.) By \cref{cyclic_arrangement_faces}, the bounded faces of ${\mathcal{H}^W}$ are labeled by sign vectors in $\overline{\Sign_{5,2,1}}$, and the bounded regions by sign vectors in $\Sign_{5,2,1}$. The unbounded regions are labeled by the sign vectors $\sigma\in\{0,+,-\}^5$ satisfying $\overline{\var}(\sigma) \le 1$. By \cref{cyclic_arrangement}, we have $\mathcal{B}_{5,2,1}(W)\cong B({\mathcal{H}^W})$.
\end{eg}

\afterpage{\begin{figure}[p!]
$$
\begin{tikzpicture}[baseline=(current bounding box.center)]
\pgfmathsetmacro{\unit}{93.6};
\pgfmathsetmacro{\epsilon}{0.04};
\pgfmathsetmacro{\SWx}{0.240};
\pgfmathsetmacro{\SWy}{0.095};
\pgfmathsetmacro{\NEx}{0.410};
\pgfmathsetmacro{\NEy}{0.265};
\node[inner sep=0,left]at(\SWx*\unit,\unit/10){$\ell_1$};
\node[inner sep=0,left]at(\SWx*\unit,{(1-\SWx)*\unit/6}){$\ell_2$};
\node[inner sep=0,left]at(\SWx*\unit,{(1-2*\SWx)*\unit/3}){$\ell_3$};
\node[inner sep=0,left]at({(1-\NEy)*\unit/3},\NEy*\unit){$\ell_4$};
\node[inner sep=0,above]at(\unit/4,\NEy*\unit){$\ell_5$};
\node[inner sep=0]at(0.264*\unit,0.175*\unit){\scriptsize{${+}{+}{+}{-}{+}$}};
\node[inner sep=0]at(0.268*\unit,0.134*\unit){\scriptsize{${+}{+}{-}{-}{+}$}};
\node[inner sep=0]at(0.272*\unit,0.110*\unit){\scriptsize{${+}{-}{-}{-}{+}$}};
\node[inner sep=0]at(0.306*\unit,0.121*\unit){\scriptsize{${+}{+}{-}{+}{+}$}};
\node[inner sep=0]at(0.32*\unit,0.106*\unit){\scriptsize{${+}{-}{-}{+}{+}$}};
\node[inner sep=0]at(0.36*\unit,0.1027*\unit){\scriptsize{${+}{-}{+}{+}{+}$}};
\node[inner sep=0]at(\unit/4+\epsilon,5*\unit/24){\scriptsize{${+}{+}{+}{-}{0}$}};
\node[inner sep=0]at(\unit/4+\epsilon,7*\unit/48){\scriptsize{${+}{+}{-}{-}{0}$}};
\node[inner sep=0]at(\unit/4+\epsilon,9*\unit/80){\scriptsize{${+}{-}{-}{-}{0}$}};
\node[inner sep=0]at(15*\unit/56,11*\unit/56){\scriptsize{${+}{+}{+}{0}{+}$}};
\node[inner sep=0]at(15*\unit/56,13*\unit/84){\scriptsize{${+}{+}{0}{-}{+}$}};
\node[inner sep=0]at(37*\unit/136,33*\unit/272){\scriptsize{${+}{0}{-}{-}{+}$}};
\node[inner sep=0]at(11*\unit/40,\unit/10-\epsilon){\scriptsize{${0}{-}{-}{-}{+}$}};
\node[inner sep=0]at(69*\unit/238,31*\unit/238){\scriptsize{${+}{+}{-}{0}{+}$}};
\node[inner sep=0]at(101*\unit/340,37*\unit/340){\scriptsize{${+}{-}{-}{0}{+}$}};
\node[inner sep=0]at(13*\unit/42,8*\unit/63){\scriptsize{${+}{+}{0}{+}{+}$}};
\node[inner sep=0]at(16*\unit/51,35*\unit/306){\scriptsize{${+}{0}{-}{+}{+}$}};
\node[inner sep=0]at(13*\unit/40,\unit/10-\epsilon){\scriptsize{${0}{-}{-}{+}{+}$}};
\node[inner sep=0]at(41*\unit/120,19*\unit/180){\scriptsize{${+}{-}{0}{+}{+}$}};
\node[inner sep=0]at(11*\unit/30,19*\unit/180){\scriptsize{${+}{0}{+}{+}{+}$}};
\node[inner sep=0]at(3*\unit/8,\unit/10-\epsilon){\scriptsize{${0}{-}{+}{+}{+}$}};
\node[inner sep=0]at(\unit/4+\epsilon,\unit/4){\scriptsize{${+}{+}{+}{0}{0}$}};
\node[inner sep=0]at(\unit/4+\epsilon,\unit/6){\scriptsize{${+}{+}{0}{-}{0}$}};
\node[inner sep=0]at(\unit/4+\epsilon,\unit/8){\scriptsize{${+}{0}{-}{-}{0}$}};
\node[inner sep=0]at(\unit/4+\epsilon,\unit/10-\epsilon){\scriptsize{${0}{-}{-}{-}{0}$}};
\node[inner sep=0]at(2*\unit/7,\unit/7){\scriptsize{${+}{+}{0}{0}{+}$}};
\node[inner sep=0]at(5*\unit/17,2*\unit/17){\scriptsize{${+}{0}{-}{0}{+}$}};
\node[inner sep=0]at(3*\unit/10,\unit/10-\epsilon){\scriptsize{${0}{-}{-}{0}{+}$}};
\node[inner sep=0]at(\unit/3,\unit/9){\scriptsize{${+}{0}{0}{+}{+}$}};
\node[inner sep=0]at(7*\unit/20,\unit/10-\epsilon){\scriptsize{${0}{-}{0}{+}{+}$}};
\node[inner sep=0]at(2*\unit/5,\unit/10-\epsilon){\scriptsize{${0}{0}{+}{+}{+}$}};
\clip(\SWx*\unit,\SWy*\unit)rectangle(\NEx*\unit,\NEy*\unit);
\draw[domain=0:1,smooth,variable=\x]plot(\x*\unit,\unit/10);
\draw[domain=0:1,smooth,variable=\x]plot(\x*\unit,{(1-\x)*\unit/6});
\draw[domain=0:1,smooth,variable=\x]plot(\x*\unit,{(1-2*\x)*\unit/3});
\draw[domain=0:1,smooth,variable=\x]plot(\x*\unit,{(1-3*\x)*\unit});
\draw[domain=0:1,smooth,variable=\y]plot(\unit/4,\y*\unit);
\end{tikzpicture}
$$
\caption{The hyperplane arrangement ${\mathcal{H}^W}$ from \cref{hyperplane_arrangement_example_n=5}, with $\mathcal{B}_{5,2,1}(W)\cong B({\mathcal{H}^W})$. Its bounded faces are labeled by sign vectors.}\label{hyperplane_arrangement_sign_vectors_5_2_1}
\end{figure}\clearpage}
\afterpage{\begin{figure}[p!]
$$
\begin{tikzpicture}[baseline=(current bounding box.center)]
\pgfmathsetmacro{\unit}{93.6};
\pgfmathsetmacro{\epsilon}{0.04};
\pgfmathsetmacro{\SWx}{0.240};
\pgfmathsetmacro{\SWy}{0.095};
\pgfmathsetmacro{\NEx}{0.410};
\pgfmathsetmacro{\NEy}{0.265};
\node[inner sep=0,left]at(\SWx*\unit,\unit/10){$\ell_1$};
\node[inner sep=0,left]at(\SWx*\unit,{(1-\SWx)*\unit/6}){$\ell_2$};
\node[inner sep=0,left]at(\SWx*\unit,{(1-2*\SWx)*\unit/3}){$\ell_3$};
\node[inner sep=0,left]at({(1-\NEy)*\unit/3},\NEy*\unit){$\ell_4$};
\node[inner sep=0,above]at(\unit/4,\NEy*\unit){$\ell_5$};
\clip(\SWx*\unit,\SWy*\unit)rectangle(\NEx*\unit,\NEy*\unit);
\draw[domain=0:1,smooth,variable=\x]plot(\x*\unit,\unit/10);
\draw[domain=0:1,smooth,variable=\x]plot(\x*\unit,{(1-\x)*\unit/6});
\draw[domain=0:1,smooth,variable=\x]plot(\x*\unit,{(1-2*\x)*\unit/3});
\draw[domain=0:1,smooth,variable=\x]plot(\x*\unit,{(1-3*\x)*\unit});
\draw[domain=0:1,smooth,variable=\y]plot(\unit/4,\y*\unit);
\node[inner sep=0,fill=white]at(0.264*\unit,0.175*\unit){\scalebox{0.4}{$\begin{ytableau}+ \\ +\end{ytableau}$}};
\node[inner sep=0,fill=white]at(0.268*\unit,0.134*\unit){\scalebox{0.4}{$\begin{ytableau}0 & + \\ +\end{ytableau}$}};
\node[inner sep=0,fill=white]at(0.272*\unit,0.110*\unit){\scalebox{0.4}{$\begin{ytableau}0 & 0 & + \\ +\end{ytableau}$}};
\node[inner sep=0,fill=white]at(0.306*\unit,0.1204*\unit){\scalebox{0.4}{$\begin{ytableau}0 & + \\ 0 & +\end{ytableau}$}};
\node[inner sep=0,fill=white]at(0.32*\unit,0.107*\unit){\scalebox{0.4}{$\begin{ytableau}0 & 0 & + \\ 0 & +\end{ytableau}$}};
\node[inner sep=0,fill=white]at(0.3586*\unit,0.1031*\unit){\scalebox{0.4}{$\begin{ytableau}0 & 0 & + \\ 0 & 0 & +\end{ytableau}$}};
\node[inner sep=0,fill=white]at(\unit/4,5*\unit/24){\scalebox{0.4}{$\begin{ytableau}+\end{ytableau}$}};
\node[inner sep=0,fill=white]at(\unit/4,7*\unit/48){\scalebox{0.4}{$\begin{ytableau}0 & +\end{ytableau}$}};
\node[inner sep=0,fill=white]at(\unit/4,9*\unit/80){\scalebox{0.4}{$\begin{ytableau}0 & 0 & +\end{ytableau}$}};
\node[inner sep=0,fill=white]at(15*\unit/56,11*\unit/56){\scalebox{0.4}{$\begin{ytableau}+ \\ 0\end{ytableau}$}};
\node[inner sep=0,fill=white]at(15*\unit/56,13*\unit/84){\scalebox{0.4}{$\begin{ytableau}0 \\ +\end{ytableau}$}};
\node[inner sep=0,fill=white]at(37*\unit/136,33*\unit/272){\scalebox{0.4}{$\begin{ytableau}0 & 0 \\ +\end{ytableau}$}};
\node[inner sep=0,fill=white]at(11*\unit/40,\unit/10-\epsilon){\scalebox{0.4}{$\begin{ytableau}0 & 0 & 0 \\ +\end{ytableau}$}};
\node[inner sep=0,fill=white]at(69*\unit/238,31*\unit/238){\scalebox{0.4}{$\begin{ytableau}0 & + \\ 0\end{ytableau}$}};
\node[inner sep=0,fill=white]at(101*\unit/340,37*\unit/340){\scalebox{0.4}{$\begin{ytableau}0 & 0 & + \\ 0\end{ytableau}$}};
\node[inner sep=0,fill=white]at(13*\unit/42,8*\unit/63){\scalebox{0.4}{$\begin{ytableau}0 & + \\ 0 & 0\end{ytableau}$}};
\node[inner sep=0,fill=white]at(16*\unit/51,35*\unit/306){\scalebox{0.4}{$\begin{ytableau}0 & 0 \\ 0 & +\end{ytableau}$}};
\node[inner sep=0,fill=white]at(13*\unit/40,\unit/10-\epsilon){\scalebox{0.4}{$\begin{ytableau}0 & 0 & 0 \\ 0 & +\end{ytableau}$}};
\node[inner sep=0,fill=white]at(41*\unit/120,19*\unit/180){\scalebox{0.4}{$\begin{ytableau}0 & 0 & + \\ 0 & 0\end{ytableau}$}};
\node[inner sep=0,fill=white]at(11*\unit/30+0.001*\unit,19*\unit/180){\scalebox{0.4}{$\begin{ytableau}0 & 0 & + \\ 0 & 0 & 0\end{ytableau}$}};
\node[inner sep=0,fill=white]at(3*\unit/8,\unit/10-\epsilon){\scalebox{0.4}{$\begin{ytableau}0 & 0 & 0 \\ 0 & 0 & +\end{ytableau}$}};
\node[inner sep=0,fill=white]at(\unit/4,\unit/4){\scriptsize{$\emptyset$}};
\node[inner sep=0,fill=white]at(\unit/4,\unit/6){\scalebox{0.4}{$\begin{ytableau}0\end{ytableau}$}};
\node[inner sep=0,fill=white]at(\unit/4,\unit/8){\scalebox{0.4}{$\begin{ytableau}0 & 0\end{ytableau}$}};
\node[inner sep=0,fill=white]at(\unit/4,\unit/10-\epsilon){\scalebox{0.4}{$\begin{ytableau}0 & 0 & 0\end{ytableau}$}};
\node[inner sep=0,fill=white]at(2*\unit/7,\unit/7){\scalebox{0.4}{$\begin{ytableau}0 \\ 0\end{ytableau}$}};
\node[inner sep=0,fill=white]at(5*\unit/17,2*\unit/17){\scalebox{0.4}{$\begin{ytableau}0 & 0 \\ 0\end{ytableau}$}};
\node[inner sep=0,fill=white]at(3*\unit/10,\unit/10-\epsilon){\scalebox{0.4}{$\begin{ytableau}0 & 0 & 0 \\ 0\end{ytableau}$}};
\node[inner sep=0,fill=white]at(\unit/3,\unit/9){\scalebox{0.4}{$\begin{ytableau}0 & 0 \\ 0 & 0\end{ytableau}$}};
\node[inner sep=0,fill=white]at(7*\unit/20,\unit/10-\epsilon){\scalebox{0.4}{$\begin{ytableau}0 & 0 & 0 \\ 0 & 0\end{ytableau}$}};
\node[inner sep=0,fill=white]at(2*\unit/5,\unit/10-\epsilon){\scalebox{0.4}{$\begin{ytableau}0 & 0 & 0 \\ 0 & 0 & 0\end{ytableau}$}};
\end{tikzpicture}
$$
\caption{The hyperplane arrangement ${\mathcal{H}^W}$ from \cref{hyperplane_arrangement_example_n=5}, with $\mathcal{B}_{5,2,1}(W)\cong B({\mathcal{H}^W})$. Its bounded faces are labeled by \Le -diagrams.}\label{hyperplane_arrangement_le_diagrams_5_2_1}
\end{figure}
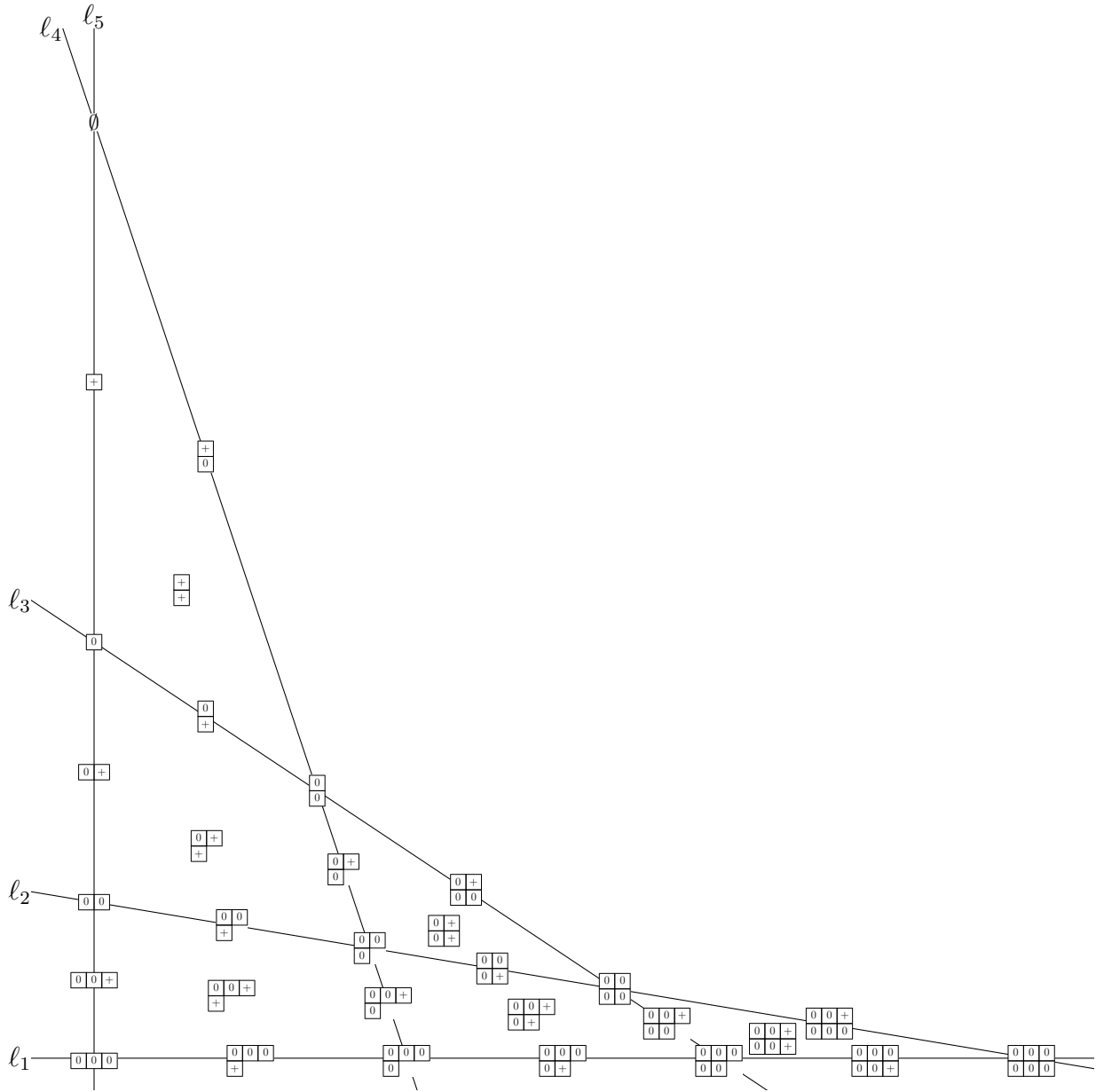\clearpage}

\begin{pf}[of \cref{cyclic_arrangement_faces}]
(i) By \cref{bounded_faces_lemma} and \cref{Sign_PSign}(ii), $\mathcal{V}_{B({\mathcal{H}^W})}$ equals either $\overline{\Sign_{n,k,1}}$ or $-\overline{\Sign_{n,k,1}}$. We must rule out the latter possibility. Recall that by construction, $w^{(0)}$ is positively oriented with respect to $V:=\spn(w^{(1)}, \dots, w^{(k)})$. According to \cref{w0_convention}, there exist $x_1, \dots, x_k\in\mathbb{R}$ such that $u := x_1w^{(1)} + \cdots + x_kw^{(k)} + w^{(0)}$ satisfies $\var(u) = \overline{\var}(u) = k$ and $u_1 > 0$. Let $\sigma := \sign(u)$. Then $\sigma$ labels a face of ${\mathcal{H}^W}$, which is a bounded face by \cref{unbounded_faces_lemma}(ii). We get $\sigma\in\mathcal{V}_{B({\mathcal{H}^W})}$, and $\sigma_1 = +$ implies $\mathcal{V}_{B({\mathcal{H}^W})}\neq -\overline{\Sign_{n,k,1}}$.

(ii) One direction follows from \cref{unbounded_faces_lemma}(ii). For the other direction, we will use the following fact about generic hyperplane arrangements $\mathcal{H}$ in $\mathbb{R}^k$: if $\tau$ labels a face of $\mathcal{H}$, then $\sigma$ also labels a face of $\mathcal{H}$ for all $\sigma\ge\tau$. 
(This follows from the fact that
the normal vectors of any $k$ or fewer hyperplanes are 
linearly independent.)

Given $\sigma\in\{0,+,-\}^n$ with $\overline{\var}(\sigma)\le k-1$, we must show that $\sigma$ labels a face of ${\mathcal{H}^W}$ (whence this face is unbounded by part (i)). Since ${\mathcal{H}^W}$ is generic, it suffices to construct $\tau\le\sigma$ which labels a face of ${\mathcal{H}^W}$. Our strategy will be to modify the sign vector $\sigma$ until we get a sign vector $\tau\le\sigma$ with $\tau\in\overline{\Sign_{n,k,1}}$, which then implies that 
$\tau$ labels a bounded face of ${\mathcal{H}^W}$ by part (i).

Set $\sigma' := \alt(\sigma)$. By \cref{dual_translation}(i), we have $\var(\sigma') \ge n-k$. In particular, the fact that $k < n$ implies that $\sigma'$ has a positive component. Take $i\in [n]$ minimum with $\sigma'_i = +$, and let $\sigma''$ be obtained from $\sigma'$ by setting to zero all components $j$ with $j < i$. Note that $\var(\sigma'') \ge \var(\sigma') - 1 \ge n-k-1$. Now we repeatedly set the last nonzero component of $\sigma''$ to zero, until we obtain a sign vector $\tau'$ with $\var(\tau') = n-k-1$. Letting $\tau := \alt(\tau')$, we have $\tau\le\sigma$, and $\overline{\var}(\tau) = k$ by \cref{dual_translation}(i). Since the first nonzero component of $\tau'$ equals $+$, we have $\tau\in\overline{\Sign_{n,k,1}}$, as desired.
\end{pf}
We are now ready to show that $\mathcal{B}_{n,k,1}(W)\cong B({\mathcal{H}^W})$.
\begin{thm}\label{cyclic_arrangement}
The restriction of $\Psi_{\mathcal{H}^W}$ to $B({\mathcal{H}^W})$ is a homeomorphism from 
$B({\mathcal{H}^W})$ to $\mathcal{B}_{n,k,1}(W)$, and induces an isomorphism of posets on the strata of $B({\mathcal{H}^W})$ and $\mathcal{B}_{n,k,1}(W)$. Explicitly, $\Psi_{\mathcal{H}^W}$ sends the stratum $\psi_{\mathcal{H}^W}^{-1}(\sigma)$ of $B({\mathcal{H}^W})$ to the stratum $\mathcal{B}_\sigma (W)$ of $\mathcal{B}_{n,k,1}(W)$, for all $\sigma\in\overline{\Sign_{n,k,1}}$.
\end{thm}

\begin{pf}
To see that $\Psi_{\mathcal{H}^W}$ is injective, note that 
if $\Psi_{{\mathcal{H}^W}}(x)=\Psi_{{\mathcal{H}^W}}(y)$ for some $x,y\in\mathbb{R}^k$, then there exists $t\in\mathbb{R}\setminus\{0\}$
such that 
$x_1w^{(1)} + \cdots + x_kw^{(k)} + w^{(0)} = 
t(y_1w^{(1)} + \cdots + y_kw^{(k)} + w^{(0)}).$  The linear 
independence of the vectors $w^{(0)}, w^{(1)}, \dots, w^{(k)}$ implies that 
$t = 1$ and $x_i = y_i$ for all $i\in [k]$, so $x = y$. 

Recall that by \cref{B_m=1}, 
$\mathcal{B}_{n,k,1}(W) 
= \{w\in\mathbb{P}(W) : \overline{\var}(w) = k\}$. 
Let us show that  $\Psi_{\mathcal{H}^W}(\R^k)$ contains 
$\mathcal{B}_{n,k,1}(W)$.
It suffices to observe that when we express an element of $\mathcal{B}_{n,k,1}(W)$ as a linear combination of $w^{(0)}, w^{(1)}, \dots, w^{(k)}$, the coefficient of $w^{(0)}$ is nonzero. This follows from \cref{gantmakher_krein}(ii),
and the fact that $\spn(w^{(1)}, \dots, w^{(k)}) \in \Gr_{k,n}^{>0}$.
 
Now we let $\mathcal{Q} := 
\Psi_{\mathcal{H}^W}^{-1}(
\mathcal{B}_{n,k,1}(W))$. 
We have shown that $\Psi_{\mathcal{H}^W}$ is a homeomorphism from 
$\mathcal{Q}$ to $\mathcal{B}_{n,k,1}$. Recall that $\mathcal{B}_{n,k,1}(W)$ is stratified by the sign vectors in $\overline{\PSign_{n,k,1}(W)}$. Therefore $\mathcal{Q}$ is the union of the faces of ${\mathcal{H}^W}$ labeled by $\sigma\in\{0,+,-\}^n$ satisfying $\overline{\var}(\sigma) = k$. By \cref{cyclic_arrangement_faces}, this is precisely $B({\mathcal{H}^W})$. The fact that $\Psi_{\mathcal{H}^W}$ induces a poset isomorphism on strata follows from \cref{bounded_faces_lemma}.
\end{pf}

\begin{rmk}\label{rmk:labeling}
It follows from \cref{cyclic_arrangement} 
that the amplituhedron $\mathcal{B}_{n,k,1}(W)$ is a regular cell complex, and in particular
its strata 
$\mathcal{B}_{\sigma}(W)$ are homeomorphic to open balls.  Using the results of 
\cref{sec_induced_subcomplex}, we can also index the cells of the amplituhedron by \Le-diagrams and matroids, 
in which case we will use the notation $\mathcal{B}_D(W)$ and $\mathcal{B}_M(W)$, respectively.
\end{rmk}

\begin{cor}\label{homeomorphic_to_ball}
The amplituhedron $\mathcal{B}_{n,k,1}(W)$ (and also $\mathcal{A}_{n,k,1}(Z)$)
is homeomorphic to a ball of dimension $k$.
\end{cor}

\begin{pf}
This follows from 
\cref{cyclic_arrangement} together with Dong's result
\cite[Theorem 3.1]{Dong} that the bounded complex of 
a uniform affine oriented matroid (of which the bounded complex of a generic hyperplane arrangement is a special case) is a piecewise linear ball.
\end{pf}

Recall from \cref{triangulation_m=1} that the amplituhedron
$\mathcal{B}_{n,k,1}(W)$ is homeomorphic to a subcomplex of 
$\Gr_{k,n}^{\geq 0}$; indeed, $\mathcal{B}_{n,k,1}(W)$
 inherits a cell decomposition
whose face poset is an induced subposet of the face poset of 
$\Gr_{k,n}^{\geq 0}$.
Another way to prove \cref{homeomorphic_to_ball} would be to show that 
the face poset of $\mathcal{B}_{n,k,1}(W)$, with a new top element
$\hat{1}$ adjoined, is \emph{shellable}.  Then since the 
cell decomposition of $\mathcal{B}_{n,k,1}(W)$ is regular, and its face poset is {\itshape pure} and {\itshape subthin} (this follows from \cref{cyclic_arrangement}), a 
result of Bj\"{o}rner \cite[Proposition 4.3(c)]{Bjorner} would imply that 
$\mathcal{B}_{n,k,1}(W)$ is homeomorphic to a ball.
\begin{prob}\label{EL}
Show that the face poset of $\mathcal{B}_{n,k,1}(W)$ 
with a top element $\hat{1}$ adjoined is shellable,
e.g.\ by finding an EL-labeling.
\end{prob}

\begin{rmk}
Note that in earlier work the second author proved that the 
face poset of $\Gr_{k,n}^{\geq 0}$ is thin and shellable
\cite{Williams}, which shows that the same is true for any 
induced subposet.  However, this does not solve 
\cref{EL}, because after adjoining $\hat{1}$, the face
poset of $\mathcal{B}_{n,k,1}(W)$ is no longer an induced subposet
of the face poset of $\Gr_{k,n}^{\geq 0}$.
\end{rmk}

As a further corollary 
of \cref{cyclic_arrangement},
we obtain the generating function for the stratification of $\mathcal{B}_{n,k,1}(W)$ with respect to dimension, since Buck found the corresponding generating function of $B(\mathcal{H})$ for a generic hyperplane arrangement $\mathcal{H}$ (which only depends on its dimension and the number of hyperplanes).
\begin{cor}[\cite{buck}]\label{generating_function_m=1}
Let $f_{n,k,1}(q) := \sum_{\text{strata $S$ of $\mathcal{B}_{n,k,1}(W)$}}q^{\dim(S)}\in\mathbb{N}[q]$ be the generating function for the stratification of $\mathcal{B}_{n,k,1}(W)$, with respect to dimension. Then
$$
f_{n,k,1}(q) = \sum_{i=0}^k\binom{n-k-1+i}{i}\binom{n}{k-i}q^i = \sum_{j=0}^k\binom{n-k-1+j}{j}(1+q)^j.
$$
\end{cor}
For example, we have $f_{5,3,1}(q) = 4q^3 + 15q^2 + 20q + 10$, which we invite the reader to verify from \cref{fig:hyperplanes}.
\begin{pf}
By the corollary to Theorem 3 of \cite{buck}, for $i\in\mathbb{N}$ the coefficient of $q^i$ in $f_{n,k,1}(q)$ equals $\frac{k+1}{n-k+i}\binom{k}{i}\binom{n}{k+1}$, which we can rewrite as $\binom{n-k-1+i}{i}\binom{n}{k-i}$. This gives the first sum above. For the last equality above, note that the coefficient of $q^i$ in $\sum_{j=0}^k\binom{n-k-1+j}{j}(1+q)^j$ equals
\begin{multline*}
\sum_{j=i}^k\binom{n-k-1+j}{j}\binom{j}{i} = \binom{n-k-1+i}{i}\sum_{j=i}^k\binom{n-k-1+j}{n-k-1+i} \\
= \binom{n-k-1+i}{i}\binom{n}{n-k+i}
\end{multline*}
by the hockey-stick identity.
\end{pf}

\begin{rmk}
By substituting $q=-1$ into 
the last expression in \cref{generating_function_m=1}, it is easy to check that the Euler characteristic of $\mathcal{B}_{n,k,1}(W)$ equals $1$.
\end{rmk}

\section{How cells of \texorpdfstring{$\mathcal{A}_{n,k,1}$}{A(n,k,1)} fit together}\label{sec_glue}

\noindent In this section we will address how cells of the $m=1$ amplituhedron fit 
together.  In particular, we will explicitly work out when
two maximal cells are adjacent, and which cells lie in the boundary 
of $\mathcal{B}_{n,k,1}(W)$, in terms of \Le -diagrams, sign vectors, positroids, and decorated permutations. See \cref{hyperplane_arrangement_sign_vectors_5_2_1} and \cref{hyperplane_arrangement_le_diagrams_5_2_1} for examples 
when  $n=5$
and  $k=2$.

\begin{prop}[Adjacency of maximal cells in the $m=1$ amplituhedron]\label{m=1_adjacency} ~\\
Given $D_1, D_2\in\mathcal{D}_{n,k,1}$, the following are equivalent:
\begin{enumerate}
\item[(i)] the cells $\mathcal{B}_{D_1}(W)$ and $\mathcal{B}_{D_2}(W)$ in 
$\mathcal{B}_{n,k,1}(W)$ are adjacent, i.e.\ their closures 
intersect in a cell $\mathcal{B}_{D'}(W)$ of codimension $1$, where  $D'\in\overline{\mathcal{D}_{n,k,1}}$ is
 necessarily  unique;
\item[(ii)] the Young diagrams of $D_1$ and $D_2$ differ by a single box, in which case we obtain $D'$ from either $D_1$ or $D_2$ by including this box with a $0$ inside it;
\item[(iii)] there exists $2 \le i \le n-1$ such that the sign vectors $\sigma(D_1)$ and $\sigma(D_2)$ differ precisely in component $i$, in which case we obtain $\sigma(D')$ from either $\sigma(D_1)$ or $\sigma(D_2)$ by setting component $i$ to $0$;
\item[(iv)] there exists $2 \le i \le n-1$ such that we can obtain the partition of $[n]$ for $M(D_1)$ (in the sense of \cref{positroids_m=1}) from the partition of $[n]$ for $M(D_2)$ by moving $i$ from one interval to another, in which case we obtain $M(D')$ from either $M(D_1)$ or $M(D_2)$ by turning $i$ into a coloop;
\item[(v)] there exists $2 \le i \le n-1$ such that 
$\pi(D_2) = 
s_{i-1}\pi(D_1)s_i,$
 in which case $\pi(D')$ equals either $\pi(D_1)s_i = s_{i-1}\pi(D_2)$ or $s_{i-1}\pi(D_1) = \pi(D_2)s_i$, whichever has exactly $k-1$ inversions. 
(Here $s_j$ denotes the simple transposition exchanging $j$ and $j+1$, and all fixed points are colored black.)
\end{enumerate}
\end{prop}

\begin{pf}
The uniqueness of $D'$ and the equivalence (i) $\Leftrightarrow$ (iii) follows from the fact that $\mathcal{B}_{n,k,1}(W)\cong B({\mathcal{H}^W})$ (\cref{cyclic_arrangement}); note that since $\sigma(D_1),
\sigma(D_2)\in \Sign_{n,k,1}$, the sign vectors cannot differ in their first or last component. The equivalence (ii) $\Leftrightarrow$ (iii) follows from \cref{DS_bijection}, and (ii) $\Leftrightarrow$ (iv) follows from the bijection in \cref{le-diagram_m=1}. The equivalence (ii) $\Leftrightarrow$ (v) follows from the bijection in \cref{BCFW_diagrams_permutations}, using the following explicit description of $\pi' := s_{i-1}\pi s_i$ for any $m=1$ BCFW permutation $\pi$:
\begin{itemize}
\item if $i$ is the minimum value in its cycle and not a fixed point, then we obtain $\pi'$ from $\pi$ in cycle notation by moving $i$ to the cycle with $i-1$, to the left of $i-1$;
\item if $i$ is the maximum value in its cycle and not a fixed point, then we obtain $\pi'$ from $\pi$ in cycle notation by moving $i$ to the cycle with $i+1$, to the right of $i+1$;
\item if $\pi(i) = i$, then $\pi'$ has  $k+1$ anti-excedances (and hence does not index a cell of $\Gr_{k,n}^{\ge 0}$, by \cref{le_permutation_bijection});
\item if $i-1$, $i$, and $i+1$ are all in the same cycle, then $\pi' = \pi$. \qedhere
\end{itemize}
\end{pf}

Since our stratification of the amplituhedron $\mathcal{B}_{n,k,1}(W)$ is a regular cell decomposition of 
a ball (see \cref{rmk:labeling} and
\cref{homeomorphic_to_ball}), it is interesting to characterize which cells
(necessarily of codimension at least $1$) comprise its boundary.
Note that by our identification of cell complexes
 $\mathcal{B}_{n,k,1}(W)\cong B({\mathcal{H}^W})$, every cell of $\mathcal{B}_{n,k,1}(W)$ lies in 
either the interior or the boundary of $\mathcal{B}_{n,k,1}(W)$.

\begin{prop}[Boundary of the $m=1$ amplituhedron]\label{m=1_boundary} ~\\
Given $D\in\overline{\mathcal{D}_{n,k,1}}$, the following are equivalent:
\begin{enumerate}
\item[(i)] 
the cell $\mathcal{B}_D(W)$ of $\mathcal{B}_{n,k,1}(W)$ is contained in the interior;
\item[(ii)] $D$ has $k$ (nonempty) rows, and for all $r\in [k]$ such that row $r$ of $D$ has no $+$'s, row $r-1$ of $D$ is longer than row $r$ (where row $0$ has length $n-k$);
\item[(iii)] $\var(\sigma(D)) = k$; 
\item[(iv)] $C\subseteq\{\min(E_1), \dots, \min(E_{n-k})\}\setminus\{1\}$, where $(E_1\sqcup\cdots\sqcup E_{n-k}, C)$ corresponds to $M(D)$ as in \cref{positroid_downset_m=1};
\item[(v)] if $i\in [n]$ is a white fixed point of $\pi$, then $2 \le i \le n-1$ and $i-1$ is not an anti-excedance of $\pi$.
\end{enumerate}
\end{prop}
We remark that the result applies even when $D$ is in $\mathcal{D}_{n,k,1}$ (i.e.\ the corresponding cell has full dimension), in which case all of the above properties hold.
\begin{pf}
We fix a hyperplane arrangement ${\mathcal{H}^W}$ as in \cref{defn_cyclic_arrangement}, so that $\mathcal{B}_{n,k,1}(W)\cong B({\mathcal{H}^W})$ by \cref{cyclic_arrangement}, and we let $\sigma$ denote $\sigma(D)$.

(i) $\Rightarrow$ (iii): Suppose that $\var(\sigma)\neq k$. Then $\var(\sigma) < k$, and we can construct $\tau\in\{+,-\}^n$ with $\tau\ge \sigma$ and $\var(\tau) = \var(\sigma)$. For example, do the following repeatedly: take $i\in [n]$ such that component $i$ is zero but either component $i-1$ or $i+1$ is nonzero, and make component $i$ nonzero and equal to either component $i-1$ or $i+1$. Then $\tau$ labels an unbounded face of ${\mathcal{H}^W}$ by \cref{cyclic_arrangement_faces}(ii), whose closure contains the face labeled by $\sigma$.

(iii) $\Rightarrow$ (i): The faces of ${\mathcal{H}^W}$ whose closure contains the face labeled by $\sigma$ are labeled by $\tau\in\{0,+,-\}^n$ with $\tau\ge \sigma$. If $\var(\sigma) = k$ then $\overline{\var}(\tau)\ge \var(\tau)\ge\var(\sigma) = k$, whence the face labeled by $\tau$ is bounded by \cref{cyclic_arrangement_faces}(ii).

(ii) $\Leftrightarrow$ (iii): Observe that $\var(\sigma) = \overline{\var}(\sigma)$ if and only if for all $i\in [n]$ such that $\sigma_i = 0$, we have that $i\neq 1,n$, and that $\sigma_{i-1},\sigma_{i+1}$ are nonzero and of opposite sign. This condition is equivalent to (ii) by \cref{DS_bijection}.

(ii) $\Leftrightarrow$ (iv): This follows from the bijection in 
\cref{le-diagram_m=1}.

(ii) $\Leftrightarrow$ (v): This follows from the bijection in 
\cref{le_permutation_bijection}.
\end{pf}

\section{The image in \texorpdfstring{$\mathcal{A}_{n,k,1}$}{A(n,k,1)} of an arbitrary cell of \texorpdfstring{$\Gr_{k,n}^{\ge 0}$}{the totally nonnegative Grassmannian}}\label{sec_injectivity}

\noindent In this section we study the image in the $m=1$ amplituhedron of 
an arbitrary cell of the totally nonnegative Grassmannian. In particular, 
we describe the image of an arbitrary cell  in terms of 
strata of $\mathcal{B}_{n,k,1}(W)$ (\cref{image_sign_vectors}), 
we compute the dimension of the image of an arbitrary cell (\cref{image_dimension}), and we characterize 
the cells which map injectively to the $m=1$ amplituhedron (\cref{image_as_union}).
 Since we have a regular cell decomposition of the amplituhedron 
$\mathcal{B}_{n,k,1}(W)$ using the $m=1$
BCFW cells and their closures (which can be indexed by 
the \Le-diagrams in $\overline{\mathcal{D}_{n,k,1}}$),
it is also natural to ask how to describe the image of an arbitrary cell
of $\Gr_{k,n}^{\geq 0}$ in terms of $\overline{\mathcal{D}_{n,k,1}}$.  We answer
this question (\cref{image_as_union}) for cells which map injectively
into the amplituhedron.

Let us fix a subspace $W\in\Gr_{k+1,n}^{>0}$ for the remainder of the section. Given a \Le -diagram $D$ inside a $k\times (n-k)$ rectangle, let $S_D$ denote its corresponding positroid cell, i.e.\ the cell $S_{M(D)}$ from \cref{def:positroid}, where $M(D)$ is the positroid corresponding to $D$ from \cref{sec_TNN_Grassmannian}. Recall from \cref{rmk:labeling} that
$$
\mathcal{B}_D(W) := \{V^\perp\cap W : V\in S_D\}
$$
is the image of $S_D$ in $\mathcal{B}_{n,k,1}(W)$. (It is equivalent to study the image $\tilde{Z}(S_D)$ in $\mathcal{A}_{n,k,1}(Z)$ by \cref{B_isomorphism}, where $Z$ is any $(k+1)\times n$ matrix whose rows span $W$, but we will find it more convenient to work in $\mathcal{B}_{n,k,1}(W)$.)
\begin{defn}\label{defn_sign(D)}
Let $D$ be a \Le -diagram of type $(k,n)$. Fix $V\in S_D$ and define
$$\mathcal{V}(D) := \{\sign(v) : v\in V^\perp\}\subseteq\{0,+,-\}^n.$$ 
In terms of oriented matroids \cite{bjorner_las_vergnas_sturmfels_white_ziegler}, $\mathcal{V}(D)$ is the set of {\itshape vectors} of the {\itshape positive orientation} of $M(D)$, and so does not depend on our choice of $V\in S_D$.
\end{defn}

A basic observation is that $\mathcal{B}_D(W)$ depends precisely on the sign vectors in $\mathcal{V}(D)$ which minimize $\overline{\var}(\cdot)$.

\begin{lem}\label{image_sign_vectors}
For $D$ a \Le-diagram of type $(k,n)$, 
we have 
$$\mathcal{B}_D(W) = \bigcup\{\mathcal{B}_\sigma(W) : \sigma\in\mathcal{V}(D)\text{ with }\overline{\var}(\sigma) = k\}.$$
\end{lem}
Recall that $\mathcal{B}_\sigma(W)$ is the $\sigma$-stratum of $\mathcal{B}_{n,k,1}(W)$ from \cref{defn_B_stratification}. While \cref{image_sign_vectors}
is not very explicit (it requires being able to compute the sign
vectors in $\mathcal{V}(D)$), we will give a more concrete description 
of the images of certain cells in
\cref{image_as_union}.
\begin{pf}
By \cref{B_m=1}, the left-hand side is contained in the right-hand side. Conversely, given $\sigma\in\mathcal{V}(D)$ with $\overline{\var}(\sigma) = k$ and an element $\spn(w)\in\mathcal{B}_\sigma(W)$ (where $w\in W\setminus\{0\}$), let us show that $\spn(w)\in\mathcal{B}_D(W)$. Take any $V\in S_D$. Since $\sigma\in\mathcal{V}(D)$, there exists $v\in V^\perp$ with $\sign(v) = \sigma$. 
Since $\sign(v)=\sign(w)$, there exist $c_1, \dots, c_n > 0$ such that $(c_1v_1, \dots, c_nv_n)=w$. 
We use the positive torus action (see \cref{positive_torus_remark})
to define 
$V' := \{(\frac{x_1}{c_1}, \dots, \frac{x_n}{c_n}) : x\in V\}\in S_D$, so that $(c_1v_1, \dots, c_nv_n)\in V'^\perp\cap W$. Since $\dim(V'^\perp\cap W) = 1$, we get $\spn(w) = V'^\perp\cap W\in\mathcal{B}_D(W)$. 
\end{pf}

\begin{rmk}
If $M$ is the positroid  corresponding to a \Le -diagram $D$, 
with dual positroid $M^*$, then by \cref{dual_translation}(ii) we have $\mathcal{V}(D) = 
\{\alt(\sign(u)) : u\in U\}$ for any $U\in S_{M^*}$. 
Therefore by \cref{dual_translation}(i), determining which sign vectors in $\mathcal{V}(D)$  minimize $\overline{\var}(\cdot)$ is equivalent to determining which sign vectors in $\sign(U)$ (for any $U\in S_{M^*}$) maximize $\var(\cdot)$.
\end{rmk}

Recall that 
$\overline{\mathcal{D}_{n,k,1}}$ 
is the set of \Le-diagrams with at most one $+$ per row, and each $+$ appears at the right end of its row. 
We showed in \cref{prop:induce} that 
for $D\in \overline{\mathcal{D}_{n,k,1}}$, 
$\mathcal{V}(D)$ contains a unique sign vector $\sigma$ (up to multiplication by $\pm 1$) with $\overline{\var}(\sigma) = k$, which we denoted by $\sigma(D)$ (\cref{def:Lesign}). In this case, we have $\mathcal{B}_D(W) = \mathcal{B}_{\sigma(D)}(W)$, as verified by \cref{image_sign_vectors}. 
Also note that by \cref{triangulation_m=1}, the dimension of $\mathcal{B}_D(W)$ is the number of $+$'s in $D$ (see \cref{hyperplane_arrangement_le_diagrams_5_2_1}). 
We now give a formula for $\dim(\mathcal{B}_D(W))$ for any \Le-diagram $D$.
\begin{prop}[Dimension of the image of an arbitrary cell]\label{image_dimension}~\\
Let $D$ be a \Le -diagram of type $(k,n)$. Then the dimension of $\mathcal{B}_D(W)$ is the number of rows of $D$ which contain a $+$.
\end{prop}
This implies that a cell of $\Gr_{k,n}^{\ge 0}$ has its dimension preserved when mapped by $\tilde{Z}$ to the $m=1$ amplituhedron if and only if its \Le -diagram has at most one $+$ in each row. Lam \cite[Theorem 4.2]{lam_amplituhedron_cells} gave an alternative criterion for general $m$ in terms of the affine Stanley symmetric function associated to the decorated permutation of the cell. This is related to the notion of {\itshape kinematical support} (see \cite[Chapter 10]{abcgpt}, \cite[Definition 4.3]{lam_amplituhedron_cells}).
\begin{pf}
Label the steps of the southeast border of $D$ by $1, \dots, n$ from northeast to southwest, and denote by $I\subseteq [n]$ the set of $i\in [n]$ such that $i$ labels a vertical step whose row contains no $+$'s. We will show that the codimension of $\mathcal{B}_D(W)$ equals $|I|$.

It follows from \cref{prop:perf} that $I$ is the set of coloops of $M(D)$, i.e.\ $I = \{i\in [n]: e^{(i)}\in V\}$ for any $V\in S_D$, where $e^{(i)}$ denotes the $i$th unit vector. Hence $I = \{i\in [n] : \sigma_i = 0\text{ for all }\sigma\in\mathcal{V}(D)\}$. Now let ${\mathcal{H}^W}$ be a hyperplane arrangement from \cref{defn_cyclic_arrangement}, so that $B({\mathcal{H}^W})$ is homeomorphic to $\mathcal{B}_{n,k,1}(W)$ by \cref{cyclic_arrangement}. By \cref{image_sign_vectors} we have
$$
\mathcal{B}_D(W)\subseteq\bigcup\{\mathcal{B}_\sigma(W) : \sigma\in\overline{\PSign_{n,k,1}} \text{ with } \sigma_i = 0\text{ for all }i\in I\},
$$
so the image of $\mathcal{B}_D(W)$ in $B({\mathcal{H}^W})$ under the homeomorphism $\mathcal{B}_{n,k,1}(W)\to B({\mathcal{H}^W})$ is contained in $\bigcap_{i\in I}H_i$ by \cref{cyclic_arrangement}. Since ${\mathcal{H}^W}$ is generic (\cref{defn_generic_arrangement}), the codimension of $\bigcap_{i\in I}H_i$ equals $|I|$, so the codimension of $\mathcal{B}_D(W)$ is at least $|I|$.

Conversely, note that any $X\in\Gr_{l,n}$ contains a vector which changes sign at least $l-1$ times: put an $l\times n$ matrix whose rows span $X$ into reduced row echelon form, and take the alternating sum of the rows. Now fix any $V\in S_D$. The element $\alt(V^\perp)\in\Gr_{n-k,n}$ contains a vector 
$\alt(v)$ (for some $v\in V^\perp$)
which changes sign at least $n-k-1$ times.
Since for all $i\in [n]\setminus I$ there exists $w\in V^\perp$ with $w_i\neq 0$, we may perturb $v\in V^\perp$ to make the components $[n]\setminus I$ all nonzero, without changing the sign of any nonzero components of $v$. We obtain a vector $v'\in V^\perp$ satisfying $\var(\alt(v'))\ge n-k-1$ and $v'_i\neq 0$ for all $i\in [n]\setminus I$. By \cref{dual_translation}(i), we have $\overline{\var}(v')\le k$, and since $\overline{\var}(v')\ge k$ by \cref{gantmakher_krein}(i), we get $\overline{\var}(v') = k$. Letting $\sigma := \sign(v')$, we have $\mathcal{B}_\sigma(W)\subseteq\mathcal{B}_D(W)$ by \cref{image_sign_vectors}. The codimension of $\mathcal{B}_\sigma(W)$ equals the number of zero components of $\sigma$, which is at most $|I|$. Hence the codimension of $\mathcal{B}_D(W)$ is at most $|I|$.
\end{pf}

Next we determine which cells $S_D$ are mapped injectively to the $m=1$ amplituhedron, and explicitly describe the images of such cells.
\begin{defn}\label{L-condition}
Let $\overline{\mathcal{L}_{n,k,1}}$ denote 
the set of \Le -diagrams of type $(k,n)$ 
which have at most one $+$ in each row, 
and 
which satisfy the \emph{\textup{L}-condition}: there is no $0$ which has
a $+$ above it in the same column and a $+$ to the right in the same row.
(However, we do not require each $+$ to appear at the 
right end of its row.) We let $\mathcal{L}_{n,k,1}$ denote the subset of $\overline{\mathcal{L}_{n,k,1}}$ of \Le -diagrams with exactly $k$ $+$'s.
\end{defn}
Note that 
$\mathcal{D}_{n,k,1} 
\subseteq\mathcal{L}_{n,k,1}$ and 
$\overline{\mathcal{D}_{n,k,1}}
\subseteq\overline{\mathcal{L}_{n,k,1}}$. For example, we have
$$
\mathcal{D}_{4,2,1} = \left\{\;\begin{ytableau}0 & + \\ 0 & +\end{ytableau}\;, \;\begin{ytableau}0 & + \\ +\end{ytableau}\;, \;\begin{ytableau}+ \\ +\end{ytableau}\;\right\} \quad \text{ and } \quad \mathcal{L}_{4,2,1} = \mathcal{D}_{4,2,1}\sqcup\left\{\;\begin{ytableau}+ & 0 \\ + & 0\end{ytableau}\;, \;\begin{ytableau}+ & 0 \\ +\end{ytableau}\;\right\}.
$$

\begin{rmk}\label{rmk:L-bar}
It is not hard to see that 
$\mathcal{L}_{n,k,1}$ (respectively, $\overline{\mathcal{L}_{n,k,1}}$) consists of the \Le -diagrams we obtain from diagrams in $\mathcal{D}_{n',k,1}$ (respectively, $\overline{\mathcal{D}_{n',k,1}}$) by inserting $n-n'$ columns of all $0$'s, as $n'$ ranges over all $n'\le n$.
\end{rmk}

\begin{defn}\label{defn_slide}
Given a \Le -diagram $D\in\overline{\mathcal{L}_{n,k,1}}$, we define a set of \Le -diagrams $\Slide(D)$ as follows. If $D\in\mathcal{L}_{n,k,1}$ (i.e.\ $D$ has no zero rows), then we let $\Slide(D)$ be the set of \Le -diagrams which can be obtained from $D$ by doing the following for each $+$ of $D$.
\begin{itemize}
\item[(1)] Slide the $+$ weakly to the right somewhere in the same row, say from box $b$ to box $b'$, such that the southeast corner of box $b'$ lies on the southeast border of $D$.
\item[(2)] Remove all boxes to the right of box $b'$ in the same row.
\item[(3)] If $b'\neq b$ and the entire lower edge of $b'$ lies on the southeast border of $D$, we can choose to remove box $b'$ (or not).
\end{itemize}
More generally, if $D \in \overline{\mathcal{L}_{n,k,1}}$, we label the steps of the southeast border of $D$ by $1, \dots, n$, and let $I \subseteq [n]$ be the set of $i$ which label a vertical step whose row contains no $+$'s. Let $D'$ be the \Le -diagram obtained by deleting the rows corresponding to $I$ from $D$, so $D'$ has no zero rows. Then we define $\Slide(D)$ as a set of \Le -diagrams in bijection with $\Slide(D')$, where given a \Le -diagram in $\Slide(D')$ we obtain the corresponding \Le -diagram in $\Slide(D)$ by adding a row of all $0$'s for each $i \in I$, such that its vertical step on the southeast border gets labeled by $i$ when we label the southeast border by $1, \dots, n$.
Note that $\Slide(D)\subseteq\overline{\mathcal{D}_{n,k,1}}$.
\end{defn}

\begin{eg}\label{eg_slide}
Let
$$
D := \;\begin{ytableau}
0 & + & 0 & 0 \\
0 & 0 & 0 \\
0 & + & 0 \\
0 & + & 0 \\
+
\end{ytableau}\;\in\overline{\mathcal{L}_{5,9,1}}.
$$
Then $\Slide(D)$ equals 
the set of \Le -diagrams appearing below.
\begin{gather*}
\;\begin{ytableau}
0 & 0 & + \\
0 & 0 & 0 \\
0 & 0 & + \\
0 & + \\
+
\end{ytableau}\;\quad
\;\begin{ytableau}
0 & 0 & + \\
0 & 0 & 0 \\
0 & 0 & + \\
0 & 0 & + \\
+
\end{ytableau}\;\quad
\;\begin{ytableau}
0 & 0 & + \\
0 & 0 & 0 \\
0 & 0 & + \\
0 & 0 \\
+
\end{ytableau}\;\quad
\;\begin{ytableau}
0 & 0 & 0 & + \\
0 & 0 & 0 \\
0 & 0 & + \\
0 & + \\
+
\end{ytableau}\;\quad
\;\begin{ytableau}
0 & 0 & 0 \\
0 & 0 & 0 \\
0 & 0 & + \\
0 & + \\
+
\end{ytableau}\; \\[6pt]
\;\begin{ytableau}
0 & 0 & 0 & + \\
0 & 0 & 0 \\
0 & 0 & + \\
0 & 0 & + \\
+
\end{ytableau}\;\quad
\;\begin{ytableau}
0 & 0 & 0 & + \\
0 & 0 & 0 \\
0 & 0 & + \\
0 & 0 \\
+
\end{ytableau}\;\quad
\;\begin{ytableau}
0 & 0 & 0 \\
0 & 0 & 0 \\
0 & 0 & + \\
0 & 0 & + \\
+
\end{ytableau}\;\quad
\;\begin{ytableau}
0 & 0 & 0 \\
0 & 0 & 0 \\
0 & 0 & + \\
0 & 0 \\
+
\end{ytableau}\; \\[-34pt]
\end{gather*}
\end{eg}

\begin{eg}
If $D\in\overline{\mathcal{D}_{n,k,1}}$, then each $+$ of $D$ is in the rightmost box of its row, so we cannot slide it further to the right. Hence $\Slide(D) = \{D\}$.
\end{eg}

\begin{thm}\label{image_as_union}
Let $D$ be a \Le -diagram of type $(k,n)$. Then the map $S_D\to\mathcal{B}_{n,k,1}(W), V\mapsto V^\perp\cap W$
is injective on the cell $S_D$ if and only if $D\in\overline{\mathcal{L}_{n,k,1}}$. In this case, we have
\begin{align}\label{image_equation}
\mathcal{B}_D(W)= \bigsqcup_{D'\in\Slide(D)}\mathcal{B}_{D'}(W).
\end{align}
\end{thm}

We will prove \cref{image_as_union} over the remainder of the section, divided into several steps. First, we consider two examples. 
\begin{eg}\label{ex:image}
Let $D$ be the \Le -diagram from \cref{eg_slide}. Then
\cref{image_as_union} asserts that $S_D$ maps injectively to the $m=1$ amplituhedron, and that its image is the disjoint union of the images of the $9$ cells corresponding to the \Le -diagrams in $\Slide(D)$.
\end{eg}

\begin{eg}
\cref{image_as_union} implies that if two cells of $\Gr_{k,n}^{\ge 0}$ map injectively to the $m=1$ amplituhedron, then their images are distinct. This can fail to hold if we do not assume that the cells map injectively. For example, if $n = 3$ and $k = 1$, then the cells $\{(1 : 0 : a) : a  > 0\}$ and $\{(1 : b : c) : b,c > 0\}$ of $\Gr_{1,3}^{\ge 0}$ have the same image by \cref{image_sign_vectors}, namely $\mathcal{B}_{(+,-,-)}(W)\sqcup\mathcal{B}_{(+,0,-)}(W)\sqcup\mathcal{B}_{(+,+,-)}(W)$.
\end{eg}

Now we begin proving \cref{image_as_union}. Let us characterize the \Le -diagrams in $\overline{\mathcal{L}_{n,k,1}}$ by a matroid-theoretic condition.
\begin{defn}\label{defn_circuits}
Let $D$ be a \Le-diagram of type $(k,n)$. The set of sign vectors $\mathcal{V}(D)$ from \cref{defn_sign(D)} is an induced subposet of $\{0,+,-\}^n$ from \cref{defn_sign_vector_partial_order}. We call the minimal elements of $\mathcal{V}(D)\setminus\{0\}$ the {\itshape circuits}\footnote{In the language of oriented matroids, $\mathcal{C}(D)$ is the set of (signed) circuits of the oriented matroid represented by any $V\in S_D$. The supports of circuits are precisely the (unsigned) circuits of the (unoriented) matroid $M(V)$. What we call a circuit will always be a sign vector.} of $D$, and denote the set of circuits by $\mathcal{C}(D)$. A more concrete way to think about circuits is the following. Fix $V\in S_D$ and an $(n-k)\times n$ matrix $A$ whose rows span $V^\perp$. For a subset $I\in\binom{[n]}{n-k}$ such that columns $I$ of $A$ form a basis for $\mathbb{R}^{n-k}$, let $A(I)$ be the $(n-k)\times n$ matrix whose rows span $V^\perp$ which we obtain by row reducing $A$ so as to get an identity matrix in columns $I$. Then the sign vectors of rows appearing in some matrix $A(I)$ are precisely the circuits of $D$ (up to sign).

We call $i\in [n]$ a {\itshape loop} of $D$ if the unit vector $e^{(i)}$ is contained in $\mathcal{V}(D)$, and a {\itshape coloop} if $\sigma_i = 0$ for all $\sigma\in\mathcal{V}(D)$. These are precisely the loops and coloops (from \cref{defn_matroid}) of the positroid $M(D)$ from \cref{prop:perf}. Alternatively, loops label the columns of $D$ which contain only $0$'s, and coloops label the rows which contain only $0$'s.
\end{defn}

\begin{lem}\label{L_circuits}
Let $D$ be a \Le -diagram of type $(k,n)$. Then $D\in\overline{\mathcal{L}_{n,k,1}}$ if and only if there do not exist $\sigma\in\mathcal{C}(D)$ and $a < b < c$ in $[n]$ such that
\begin{itemize}
\item[(i)] $\sigma_a, \sigma_c\neq 0$;
\item[(ii)] $\sigma_b = 0$; and
\item[(iii)] $b$ is neither a loop nor a coloop of $D$.
\end{itemize}
\end{lem}
\begin{pf}
Let $D'$ be the \Le -diagram contained in a $k'\times (n'-k')$ rectangle obtained from $D$ by deleting all rows and columns which contain only $0$'s. Equivalently, the positroid $M(D')$ is obtained from $M(D)$ by restricting the ground set $[n]$ (as in \cref{defn_restriction}) to elements which are neither loops nor coloops. Then we obtain $\mathcal{C}(D')$ from $\mathcal{C}(D)$ by deleting the circuits $\pm e^{(i)}$ corresponding to loops $i$, restricting all circuits from $[n]$ to $\{i\in [n]: i \text{ is neither a loop nor a coloop of $D$}\}$, and multiplying certain components by $-1$. Note that $D'$ has no loops or coloops, and $D\in\overline{\mathcal{L}_{n,k,1}}$ if and only if $D'\in\mathcal{D}_{n',k',1}$. Hence it suffices to prove that $D'\in\mathcal{D}_{n',k',1}$ if and only if there do not exist $\sigma\in\mathcal{C}(D')$ and $a < b < c$ in $[n]$ such that $\sigma_a,\sigma_c\neq 0$ and $\sigma_b = 0$.

($\Rightarrow$): Suppose that $D'\in\mathcal{D}_{n',k',1'}$. Let us fix $V'\in S_{D'}$, so that the representable matroid $M(V'^\perp)$ (defined in \cref{eg_matroid}) equals $M(D')^*$. By \cref{positroid_downset_m=1}, $M(D')^*$ is the direct sum of uniform matroids of rank $1$ whose ground sets are all intervals. That is, we can write $V'^\perp$ as the row span of an $(n'-k')\times n'$ matrix of the form
\begin{align}\label{V'_matrix}
\setcounter{MaxMatrixCols}{20}\begin{bmatrix}
* & \cdots & * & 0 & \cdots & 0 & 0 & \cdots & 0 & & \cdots & \\
0 & \cdots & 0 & * & \cdots & * & 0 & \cdots & 0 & & \cdots & \\
0 & \cdots & 0 & 0 & \cdots & 0 & * & \cdots & * & & \cdots & \\
 & \vdots & & & \vdots & & & \vdots & & & \ddots &
\end{bmatrix},
\end{align}
where all the $*$'s are nonzero. The sign vectors of these rows, and their negations, are precisely the circuits of $D'$. We see that the circuits of $D'$ satisfy the required condition.

($\Leftarrow$): Suppose that there do not exist $\sigma\in\mathcal{C}(D')$ and $a < b < c$ in $[n]$ such that $\sigma_a,\sigma_c\neq 0$ and $\sigma_b = 0$. Let us again take some $V'\in S_{D'}$, and put an $(n'-k')\times n'$ matrix whose rows span $V'^\perp$ into reduced row echelon form. Then the sign vector of each row of this matrix is a circuit of $D'$, and so the nonzero entries in each row form a consecutive block. It follows that this matrix is of the form \eqref{V'_matrix}, where all the $*$'s are nonzero. (The matrix has no zero columns since $D'$ has no coloops.) Hence the matroid $M(D')^* = M(V'^\perp)$ is the direct sum of uniform matroids of rank $1$ whose ground sets are all intervals, whence $D'\in\mathcal{D}_{n',k',1}$ by \cref{positroid_downset_m=1}.
\end{pf}

\begin{lem}\label{big_fiber}
Let $D$ be a \Le -diagram of type $(k,n)$. Suppose that there exist a circuit $\sigma\in\mathcal{C}(D)$ and $a < b < c$ in $[n]$ such that 
\begin{itemize}
\item[(i)] $\sigma_a, \sigma_c\neq 0$;
\item[(ii)] $\sigma_b = 0$; and
\item[(iii)] $b$ is not a coloop of $D$.
\end{itemize}
Then there exists $\tau\in\mathcal{V}(D)$ with $\overline{\var}(\tau) = k$ and $\tau_b = 0$.
\end{lem}
We will use such a sign vector $\tau$ in \cref{map_not_injective} to provide a certificate of non-injectivity, if $D\notin\overline{\mathcal{L}_{n,k,1}}$. In proving \cref{map_not_injective}, we will want to require that $b$ (satisfying $\tau_b = 0$) is not a loop of $D$, as in \cref{L_circuits}. But in order to prove just \cref{big_fiber}, we do not need to assume that $b$ is not a loop.
\begin{pf}
Fix $V\in S_D$ and an $(n-k)\times n$ matrix $A$ whose rows span $U := \alt(V^\perp)$. By \cref{dual_translation}(i), it suffices to show that there exists $w\in U$ with $\var(w) = n-k-1$ and $w_b = 0$. We will use the fact that $U\in\Gr_{n-k,n}^{\ge 0}$, by \cref{dual_translation}(ii). The idea is to take $u\in U$ with sign vector $\sigma' := \alt(\sigma)$, and perturb it so that it has $n-k-1$ sign changes. The presence of both $a$ and $c$ will allow us to get a `free' sign change without `using' $b$.

Let $J := \{j\in [n] : \sigma_j = 0\}$. The fact that $\sigma$ is a circuit implies that columns $J$ of $A$ do not span $\mathbb{R}^{n-k}$, but columns $J\cup\{a\}$ do span $\mathbb{R}^{n-k}$. Also, since $b$ is not a coloop of $D$, column $b$ of $A$ is nonzero. Hence we can extend $\{a,b\}$ to $I\in\binom{J\cup\{a\}}{n-k}$ such that columns $I$ of $A$ form a basis of $\mathbb{R}^{n-k}$. Let $A(I)$ denote the $(n-k)\times n$ matrix whose rows span $U$ which we obtain by row reducing $A$ so as to get an identity matrix in columns $I$. We denote by $v^{(i)}$ the row of $A(I)$ whose pivot column is $i$, i.e.\ $v^{(i)}_i = 1$ and $v^{(i)}_j = 0$ for all $j\in I\setminus\{i\}$.
\begin{claim}
$\sigma'_c = (-1)^{|I\cap(a,c)|}\sigma'_a$, where $(a,c)$ denotes the open interval between $a$ and $c$.
\end{claim}

\begin{claimpf}
Write $I = \{i_1 < \cdots < i_{n-k}\}$, so that $a = i_r$ for some $r\in [n-k]$. Since $u_{i_s} = 0$ for all $s\in [n-k]\setminus\{r\}$, we can use $v^{(i_s)}$ to perturb $u$, without changing the sign of $u_c$, so that component $i_s$ gets the sign
$$
\begin{cases}
(-1)^{s-r}\sigma'_a, & \text{if $i_s < c$}, \\
(-1)^{s-r+1}\sigma'_a, & \text{if $i_s > c$}.
\end{cases}
$$
Since the resulting vector $u'$ lies in $U\in \Gr_{n-k,n}^{\geq 0}$, we get $\var(u') \leq n-k-1$ by \cref{gantmakher_krein}(i). Hence $u'$ does not alternate in sign on $I\cup\{c\}$, and so $\sigma'_c\neq -(-1)^{|I\cap(a,c)|}\sigma'_a$.
\end{claimpf}

We obtain our desired vector $w\in U$ as follows. Write $I\setminus\{b\} = \{i'_1 < \cdots < i'_{n-k-1}\}$, so that $a = i'_r$ for some $r\in [n-k]$. Then for $s\in [n-k-1]\setminus\{r\}$, we use $v^{(i'_s)}$ to perturb $u$, without changing the sign of $u_c$, so that component $i'_s$ gets sign 
$$
\begin{cases}
(-1)^{s-r}\sigma'_a, & \text{if $i'_s < c$}, \\
(-1)^{s-r+1}\sigma'_a, & \text{if $i'_s > c$}.
\end{cases}
$$
The resulting vector $w$ alternates in sign on $(I\setminus\{b\})\cup\{c\}$ by the claim, so $\var(w)\ge n-k-1$. We have $\var(w)\le n-k-1$ by \cref{gantmakher_krein}(i), and $w_b = 0$ by construction.
\end{pf}

\begin{cor}\label{map_not_injective}
Let $D$ be a \Le -diagram of type $(k,n)$ which is not in $\overline{\mathcal{L}_{n,k,1}}$. Then the map $S_D\to\mathcal{B}_{n,k,1}(W), V\mapsto V^\perp\cap W$ is not injective.
\end{cor}

\begin{pf}
Fix $V\in S_D$. By \cref{L_circuits} and \cref{big_fiber}, there exist $w\in V^\perp$ and $b\in [n]$ such that $\overline{\var}(w) = k$, $w_b = 0$, and $b$ is neither a loop nor a coloop of $D$. The strategy of the proof is to use the positive torus action (see \cref{positive_torus_remark}) so as to get from $w$ an element of $\mathcal{B}_{n,k,1}(W)$, and then use the positive torus action again in component $b$ to get many elements of $S_D$ which map to our chosen element of $\mathcal{B}_{n,k,1}(W)$.

Let $\sigma := \sign(w)$, so that by \cref{positive_covectors}(i) there exists $w'\in W$ with $\sign(w') = \sigma$. 
Letting $w = (w_1,\dots, w_n)$, we can therefore 
write $w' = (c_1w_1, \dots, c_nw_n)$ for some $c_1, \dots, c_n > 0$. Letting $V' := \{(\frac{v_1}{c_1}, \dots, \frac{v_n}{c_n}) : v\in V\}\in S_D$, we have $w'\in V'^\perp$. Then for $t > 0$, we define
$$
V_t' := \{(v_1, \dots, v_{b-1}, tv_b, v_{b+1}, \dots, v_n) : v\in V'\}\in S_D.
$$
\begin{claim}
The $V_t'$ ($t > 0$) are all distinct.
\end{claim}

\begin{claimpf}
Let $M(D)$ denote the matroid from \cref{prop:perf}. Since $b$ is not a loop or a coloop of $D$, there exist $I,J\in M(D)$ with $b\in I$ and $b\notin J$. Then $\frac{\Delta_I(V_t')}{\Delta_J(V_t')}\cdot\frac{\Delta_J(V')}{\Delta_I(V')} = t$ for $t > 0$.
\end{claimpf}
Since $w'\in V_t'^\perp\cap W$ and $\dim(V_t'^\perp\cap W) = 1$, the elements $V_t'$ ($t > 0$) all map to the line spanned by $w'$ in $\mathcal{B}_{n,k,1}(W)$.
\end{pf}

We now show that if $D\in\overline{\mathcal{L}_{n,k,1}}$, then $S_D$ maps injectively to the $m=1$ amplituhedron, and we also determine its image. We need to understand the sign vectors $\mathcal{V}(D)$ for such $D$.
This can be done, thanks to the fact
that when we remove enough loops (columns containing only $0$'s) from $D$, we obtain an element of  
$\overline{\mathcal{D}_{n,k,1}}$ (see \cref{rmk:L-bar}).
We first consider the case that $D$ is in $\mathcal{L}_{n,k,1}$ (i.e.\ $D$ has no coloops) and its $+$'s lie in a single column. 
\begin{lem}\label{image_basic_case}
Suppose that $D\in\mathcal{L}_{n,k,1}$ with $k\ge 1$ such that the $+$'s of $D$ lie in a single column. Then the map $S_D\to\mathcal{B}_{n,k,1}(W), V\mapsto V^\perp\cap W$ is injective on $S_D$, and $\mathcal{B}_D(W)= \bigsqcup_{D'\in\Slide(D)}\mathcal{B}_{D'}(W)$.
\end{lem}

\begin{pf}
Let $L\in\binom{[n]}{n-k-1}$ denote the set of loops of $D$ (from \cref{defn_circuits}). That is, if we label the steps of the southeast border of $D$ from $1$ to $n$, then $L$ is the set of labels of the horizontal steps, excluding the label corresponding to the column of $+$'s. Then $M(D)$ is a matroid of rank $k$, which is uniform when restricted to $[n]\setminus L$ and has set of loops $L$. It follows from \cref{positive_covectors}(ii) that $\sigma\in\mathcal{V}(D)$ if and only if
$$
\sigma|_{[n]\setminus L}\in\{0, (+,-,+,-,\cdots), (-,+,-,+,\cdots)\}.
$$
(The sign of $\sigma_i$ for $i\in L$ is arbitrary.)

To compute $\mathcal{B}_D(W)$, 
we use \cref{image_sign_vectors}.  We need to identify
those $\sigma\in\mathcal{V}(D)$ satisfying $\overline{\var}(\sigma) = k$. Let us first consider the example $D = \;\begin{ytableau}+ & 0 & 0\end{ytableau}\;$, so $n=4$, $k=1$, $L = \{2,3\}$. We have
$$
\mathcal{V}(D) = \{(+,*,*,-)\}\cup\{(-,*,*,+)\}\cup\{(0,*,*,0)\},
$$
where each $*$ can be $0$, $+$, or $-$. Since
$\overline{\var}(0, *, *, 0) \geq 2$, the sign vectors $\sigma\in\mathcal{V}(D)$ with $\overline{\var}(\sigma) = 1$ (modulo multiplication by $\pm 1$) are
$$
(+,+,+,-), (+,+,-,-), (+,-,-,-), (+,+,0,-), (+,0,-,-),
$$
that is, a sequence of $+$'s followed by $-$'s, possibly with one $0$ in between. The corresponding \Le -diagrams (from \cref{def:Lesign}) are
$$
\;\begin{ytableau}+\end{ytableau}\;,\quad
\;\begin{ytableau}0 & +\end{ytableau}\;,\quad
\;\begin{ytableau}0 & 0 & +\end{ytableau}\;,\quad
\;\begin{ytableau}0\end{ytableau}\;,\quad
\;\begin{ytableau}0 & 0\end{ytableau}\;,
$$
respectively.

Now we describe what happens in general. Write $[n]\setminus L = \{r_1 < \cdots < r_{k+1}\}$. Given $\sigma\in\{0,+,-\}^n$, we have $\sigma\in\mathcal{V}(D)$ and $\overline{\var}(\sigma) = k$ if and only if
\begin{itemize}
\item $\sigma$ alternates in sign on $[n]\setminus L$;
\item $\sigma_i = \sigma_{r_1}$ for all $i < r_1$;
\item $\sigma_i = \sigma_{r_{k+1}}$ for all $i > r_{k+1}$; and
\item for all $j\in [k]$, there exists an integer $s_j\in [r_j, r_{j+1}-1]$ such that $\sigma_i = \sigma_{r_j}$ for all $i\in [r_j, s_j]$, $\sigma_i = \sigma_{r_{j+1}}$ for all $i\in [s_j+2, r_{j+1}]$, and if $s_j+1\neq r_{j+1}$ then $\sigma_{s_j+1}$ equals either $0$ or $\sigma_{r_{j+1}}$.
\end{itemize}
We obtain the \Le -diagram $\Omega_{DS}^{-1}(\sigma)$ (see \cref{def:Lesign}) by a slide as in \cref{defn_slide}, as follows. In step (1), we slide the $+$ in the $j$th row of $D$ so that it is $s_j-r_j$ boxes to the left of the rightmost box in the row. In step (2), we remove the boxes to the right of this $+$.  
In step (3), we remove the box $b'$ containing this $+$ if and only if $\sigma_{s_j+1} = 0$. (Note that the case when $b'\neq b$ and the bottom edge of $b'$ lies on the southeast border of $D$ corresponds to $s_j+1\neq r_{j+1}$.) Thus $\Omega_{DS}^{-1}$ restricted to $\{\sigma\in\mathcal{V}(D) : \overline{\var}(\sigma) = k\}$ is a bijection $\{\sigma\in\mathcal{V}(D) : \overline{\var}(\sigma) = k\}\to\Slide(D)$, with inverse $\Omega_{DS}$. Then \cref{image_sign_vectors} implies \eqref{image_equation}.

To see that the map $S_D\to\mathcal{B}_{n,k,1}(W)$ is injective, suppose that $w\in W$ with $\sign(w)\in\mathcal{V}(D)$ and $\overline{\var}(w)=k$. We have just observed that $w$ is nonzero when restricted to $[n]\setminus L$, so the vectors $w$ and $e^{(i)}$ for $i\in L$ are linearly independent. Their span is an element of $\Gr_{n-k,n}$, whose orthogonal complement is the unique $V\in S_D$ with $V^\perp\cap W = \spn(w)$.
\end{pf}

\begin{pf}[of \cref{image_as_union}]
If $D\notin\overline{\mathcal{L}_{n,k,1}}$, then the map $S_D\to\mathcal{B}_{n,k,1}(W)$ is not injective by \cref{map_not_injective}. Now suppose that $D\in\overline{\mathcal{L}_{n,k,1}}$. We must show that the map $S_D\to\mathcal{B}_{n,k,1}(W)$ is injective, and that its image equals $\bigsqcup_{D'\in\Slide(D)}\mathcal{B}_{D'}(W)$. We first explain how to reduce to the case that $D\in\mathcal{L}_{n,k,1}$, i.e.\ $D$ has no coloops (see \cref{defn_circuits}). Let us suppose that $i\in [n]$ is a coloop of $D$, and see what happens when we delete $i$. In terms of \Le -diagrams, $i$ labels a row of $D$ which contains only $0$'s, and deleting this row gives a \Le -diagram $D'$ indexing a cell of $\Gr_{k-1,n-1}^{\geq 0}$. In terms of cells, we have the map
\begin{gather*}
\{V^\perp: V\in S_D\}\to\{V'^\perp : V'\in S_{D'}\}, \\
\{(v_1, \dots, v_{i-1}, 0, v_{i+1}, \dots, v_n)\}\mapsto \{(v_1, \dots, v_{i-1}, -v_{i+1}, \dots, -v_n)\}.
\end{gather*}
In terms of the subspace $W\in\Gr_{k+1,n}^{>0}$, we map it to
$$
W' := \{(w_1, \dots, w_{i-1}, -w_{i+1}, \dots, -w_n) : w\in W\text{ with }w_i = 0\},
$$
where $W'\in\Gr_{k,n-1}^{>0}$ by \cref{gantmakher_krein}(ii). Note that the map $V \mapsto V^\perp \cap W$ is injective on $S_D$ if and only if 
the map $V' \mapsto {V'}^\perp \cap W'$ is injective on $S_{D'}$. We also have a bijection $\Slide(D)\to\Slide(D')$ given by deleting the row labeled by $i$. By \cref{image_sign_vectors}, this shows that the result for $D'$ implies the result for $D$. Hence we may assume that $D$ has no coloops,
i.e.\ $D\in \mathcal{L}_{n,k,1}$.

In this case, we can decompose $D$ as follows: there exist \Le -diagrams $D_1, \dots, D_l$, where each $D_j$ has no coloops and all its $+$'s lie in a single column, such that we get $D$ by gluing together $D_1, \dots, D_l$ from northeast to southwest (so that the southwest corner of $D_j$ and the northeast corner of $D_{j+1}$ coincide, for all $j\in [l-1]$), and filling the space above and to the left of the resulting diagram with $0$'s. For example, if $D = \;\begin{ytableau} 0 & 0 & + \\ 0 & 0 & + \\ +\end{ytableau}\;$, then $D_1 = \;\begin{ytableau}0 & + \\ 0 & +\end{ytableau}\;$ and $D_2 = \;\begin{ytableau}+\end{ytableau}\;$. 
We have $M(D) = M(D_1)\oplus\cdots\oplus M(D_l)$, so any $V\in S_D$ can be written uniquely as $V_1\oplus\cdots\oplus V_l$, where $V_j \in S_{D_j}$. Now for $j\in [l]$, let $E_j$ be the ground set of $M(D_j)$, and let $U_j$ denote the orthogonal 
complement of $V_j$ in $\R^{E_j}$, so that $V^{\perp} = U_1 \oplus \dots \oplus U_l$.  
It follows that the vectors $w\in V^\perp$ satisfying $\overline{\var}(w) = k$ can be written precisely as $w = u_1 + \cdots + u_l$, where $u_j\in U_j$ with $\overline{\var}(u_j) = \dim(V_j)$ for all $j\in [l]$, and the signs of the $u_j$'s are such that we never get an extra sign change from $u_j$ to $u_{j+1}$ for $j\in [l-1]$. Thus it suffices to verify injectivity and \eqref{image_equation} locally for each $D_j$, which we can do thanks to \cref{image_basic_case}.
\end{pf}

\section{Relaxing to Grassmann polytopes}\label{sec_grassmann_polytopes}

\noindent In this section we discuss what happens when we relax the condition that $W\in\Gr_{k+m,n}$ is 
totally positive, in the sense of Lam's {\itshape Grassmann polytopes} \cite{lam}. 
\begin{defn}[{\cite[Definition 15.1]{lam}}]\label{defn_grassmann_polytope}
Let $Z$ be a real $r\times n$ matrix  with row span $W\in\Gr_{k+m,n}$, where $k+m \le r$ (we allow $Z$ to not have full row rank). Suppose that $W$ contains a totally positive $k$-dimensional subspace. Then by \cite[Proposition 15.2]{lam}, the map $\tilde{Z}:\Gr_{k,n}^{\ge 0}\to\Gr_{k,r}$ is well defined, i.e.\ $\dim(Z(V)) = k$ for all $V\in\Gr_{k,n}^{\ge 0}$. We call the image $\tilde{Z}(\overline{S_M})\subseteq\Gr_{k,r}$ of the closure of a cell $S_M$ of $\Gr_{k,n}^{\ge 0}$ a {\itshape Grassmann polytope}. When $\overline{S_M} = \Gr_{k,n}^{\ge 0}$, i.e.\ $M$ is the uniform matroid of rank $k$, we call the image $\tilde{Z}(\Gr_{k,n}^{\ge 0})$ a {\itshape full Grassmann polytope}.

Analogously, given a cell $S_M$ of $\Gr_{k,n}^{\ge 0}$, let us define a {\itshape Grassmann arrangement} as
$$
\{V^\perp\cap W : V\in\overline{S_M}\}\subseteq\Gr_m(W).
$$
We call this a {\itshape full Grassmann arrangement} in the case that $\overline{S_M} = \Gr_{k,n}^{\ge 0}$. A Grassmann arrangement is well defined, and homeomorphic to the corresponding Grassmann polytope, by the same arguments which appear in \cref{sec_B_amplituhedron}.
\end{defn}

The amplituhedron $\mathcal{A}_{n,k,1}(Z)$ is an example of a full Grassmann polytope (see \cref{TP_restriction}). In the case $k=1$, Grassmann polytopes are precisely polytopes in the projective space $\Gr_{1,r} = \mathbb{P}^{r-1}$. Therefore, Grassmann polytopes generalize polytopes into Grassmannians.

\cref{intrinsic_description} generalizes to any full Grassmann arrangement; the proof is the same. In particular, in the case that $m=1$, the analogue of \cref{B_m=1} holds.

\begin{prop}
For $m=1$, the full Grassmann arrangement can be described as follows:
$$
\{V^\perp\cap W : V\in\Gr_{k,n}^{\ge 0}\} = \{w\in\mathbb{P}(W) : \overline{\var}(w)\ge k\} \subseteq\mathbb{P}(W).
$$
\end{prop}
As in \cref{sec_cyclic_arrangement}, if $m=1$ we can define a hyperplane arrangement ${\mathcal{H}^W}$ depending on $W$, and show that the full Grassmann arrangement above is homeomorphic to $B({\mathcal{H}^W})$. In slightly more detail, we take a $k$-dimensional totally positive subspace $W'$ of $W$ as provided by \cref{defn_grassmann_polytope}, and let $w^{(1)}, \dots, w^{(k)}$ be its basis. We extend this to a basis $w^{(0)}, w^{(1)}, \dots, w^{(k)}$ of $W$, and define ${\mathcal{H}^W}$ as in \cref{defn_cyclic_arrangement} (we ignore the requirement that $w^{(0)}$ is positively oriented). Then \cref{unbounded_faces_lemma} holds for ${\mathcal{H}^W}$, with the same proof. This implies (as in the first paragraph of the proof of \cref{bounded_faces_lemma}) that the face labels of $B({\mathcal{H}^W})$ are precisely the sign vectors $\sigma\in\{0,+,-\}^n$ such that $\overline{\var}(\sigma)\ge k$ and $\psi_{{\mathcal{H}^W}}^{-1}(\sigma)\neq\emptyset$. We deduce that $\Psi_{{\mathcal{H}^W}}$ restricted to $B({\mathcal{H}^W})$ is a homeomorphism onto the full Grassmann arrangement of $W$, as claimed. If all Pl\"{u}cker coordinates of $W$ are nonzero, this implies that the full Grassmann arrangement of $W$ (for $m=1$) is homeomorphic to a ball, as in \cref{homeomorphic_to_ball}. We remark that Lam has conjectured that every Grassmann polytope is contractible \cite[p.\ 112]{lam}.

We observe that as $W$ varies, we do not recover all bounded complexes of hyperplane arrangements in this way, since the condition that $W\in\Gr_{k+1,n}$ contains a subspace in $\Gr_{k,n}^{>0}$ is very restrictive. Indeed, if $\sigma$ labels a face of $B({\mathcal{H}^W})$ for such $W$, then $\overline{\var}(\sigma)\ge k$.

\begin{rmk}\label{grassmann_polytope_remark}
By \cite[Theorem 4.2]{karp}, the map $\tilde{Z}:\Gr_{k,n}^{\ge 0} \to \Gr_{k,r}$ is well defined if and only if $\var(v)\ge k$ for all nonzero $v\in\ker(Z)$, so it makes sense to consider images $\tilde{Z}(\Gr_{k,n}^{\ge 0})$ for such $Z$. Galashin's negative answer to \cref{lam_problem}(ii) \cite{galashin} demonstrates the existence of such $Z$ not satisfying the assumptions of \cref{defn_grassmann_polytope}.\footnote{The existence of a `wild' matrix $Z$ of rank $d$ is equivalent to a negative answer to \cref{lam_problem}(ii) for $l = n-k$ and $m = n - d$, by a duality argument similar to that in the proof of \cref{intrinsic_description}.} Lam has proposed calling their images $\tilde{Z}(\Gr_{k,n}^{\ge 0})$ {\itshape wild Grassmann polytopes}. The connection to Grassmann polytopes was one motivation for posing \cref{lam_problem}\footnote{\cref{lam_problem}(ii) was raised in \cite{karp}; see Theorem 4.2 and the discussion which follows it.}, and Galashin's example is an exciting development.

Similarly, we may define {\itshape wild Grassmann arrangements}. Note that in our construction above, it is essential that $W$ (the row span of $Z$) has a totally positive $k$-dimensional subspace. Indeed, the proof of \cref{unbounded_faces_lemma} (which classifies the face labels of $\mathcal{H}^W$) relies on $\spn(w^{(1)}, \dots, w^{(k)})$ being totally positive. Therefore, our arguments do not extend to the wild case. It would be interesting to explore the properties of wild Grassmann arrangements, including in the case $m=1$, and to see whether their properties are less `tame' than Grassmann arrangements as defined in \cref{defn_grassmann_polytope}. We leave this to future work.
\end{rmk}

\bibliographystyle{alpha}
\bibliography{bibliography}

\end{document}